\documentclass[a4paper,11pt,BCOR1cm,parskip]{scrartcl}

\usepackage{german}
\usepackage[iso]{umlaute}
\usepackage{theorem}
\usepackage{amsmath}
\usepackage{amssymb}
\usepackage{makeidx}
\usepackage[matrix,arrow,curve]{xy}
\usepackage{epsfig}
\usepackage{boxedminipage}
\usepackage{epic}
\usepackage{longtable}
\usepackage{ifthen}
\usepackage{tabularx}
\usepackage{parskip}

\newcommand{\R}{\mathrm{I\!R}}
\newcommand{\N}{\mathrm{I\!N}}
\newcommand{\HH}{\mathrm{I\!H}}

\newcommand{\K}{\mathrm{I\!K}}
\newcommand{\PP}{\mathrm{I\!P}}

\newcommand{\Z}{\mathchoice {\hbox{$\sf\textstyle Z\kern-0.4em
Z$}}{\hbox{$\sf\textstyle Z\kern-0.4em Z$}}{\hbox{$\sf\scriptstyle
Z\kern-0.3em Z$}}{\hbox{$\sf\scriptscriptstyle Z\kern-0.2em Z$}}}

\newcommand{\Q}{\mathchoice {\setbox0=\hbox{$\displaystyle\rm
Q$}\hbox{\raise0.15\ht0\hbox to0pt{\kern0.4\wd0\vrule
height0.8\ht0\hss}\box0}}{\setbox0=\hbox{$\textstyle\rm
Q$}\hbox{\raise0.15\ht0\hbox to0pt{\kern0.4\wd0\vrule
height0.8\ht0\hss}\box0}}{\setbox0=\hbox{$\scriptstyle\rm
Q$}\hbox{\raise0.15\ht0\hbox to0pt{\kern0.4\wd0\vrule
height0.7\ht0\hss}\box0}}{\setbox0=\hbox{$\scriptscriptstyle\rm
Q$}\hbox{\raise0.15\ht0\hbox to0pt{\kern0.4\wd0\vrule
height0.7\ht0\hss}\box0}}}

\newcommand{\OO}{\mathchoice {\setbox0=\hbox{$\displaystyle\rm
O$}\hbox{\hbox to0pt{\kern0.4\wd0\vrule
height0.9\ht0\hss}\box0}}{\setbox0=\hbox{$\textstyle\rm O$}\hbox{\hbox
to0pt{\kern0.4\wd0\vrule
height0.9\ht0\hss}\box0}}{\setbox0=\hbox{$\scriptstyle\rm O$}\hbox{\hbox
to0pt{\kern0.4\wd0\vrule
height0.9\ht0\hss}\box0}}{\setbox0=\hbox{$\scriptscriptstyle\rm
O$}\hbox{\hbox to0pt{\kern0.4\wd0\vrule height0.9\ht0\hss}\box0}}}

\newcommand{\id}{\mathrm{id}}
\newcommand{\GL}{\mathrm{GL}}

\newcommand{\SO}{\mathrm{SO}}

\newcommand{\rk}{\mathrm{rk}}

\newcommand{\eps}{\varepsilon}
\newcommand{\vi}{\varphi}
\newcommand{\vkap}{\varkappa}

\newcommand{\qmq}[1]{\quad\mbox{#1}\quad}
\newcommand{\Menge}[2]{\{\,#1\,|\,#2\,\}}

\renewcommand{\bigoplus}{\mathop{\bigcirc
  \raisebox{-0.22em}{\hskip-0.53em\hbox{\vrule height2.08ex width0.04em}
  \raisebox{ 0.48em}{\hskip-0.75em\hbox{\vrule height0.04em width 0.8em}}
  \hskip- 0.2em}}}

\newcommand {\g}[2]{\langle #1,#2\rangle}

\newcommand {\operp}{\mathbin{\mbox{$\ominus\raisebox{2.9pt}
 {\hskip-0.42em\hbox{\vrule height0.7ex width0.02em}\hskip0.42em }$}}}

\newcommand{\Ug}{\mathrm{U}}
\newcommand{\SU}{\mathrm{SU}}
\newcommand{\Sp}{\mathrm{Sp}}

\newcommand{\Spin}{\mathrm{Spin}}

\newcommand{\RE}{\mathop{\mathrm{Re}}\nolimits}
\newcommand{\IM}{\mathop{\mathrm{Im}}\nolimits}

\newcommand{\Aut}{\mathop{\mathrm{Aut}}\nolimits}

\newcommand{\ad}{\mathop{\mathrm{ad}}\nolimits}
\newcommand{\Ad}{\mathop{\mathrm{Ad}}\nolimits}

\newcommand{\Fix}{\mathop{\mathrm{Fix}}\nolimits}
\newcommand{\spn}{\mathop{\mathrm{span}}\nolimits}
\newcommand{\tr}{\mathop{\mathrm{tr}}\nolimits}

\newcommand{\diag}{\mathop{\mathrm{diag}}\nolimits}

\newcommand{\EIII}{\mathrm{EIII}}
\newcommand{\EIV}{\mathrm{EIV}}
\newcommand{\Esix}{{E}_6}
\newcommand{\Ffour}{{F}_4}
\newcommand{\Gtwo}{{G}_2}

\newcommand{\liea}{\mathfrak{a}}
\newcommand{\lieb}{\mathfrak{b}}
\newcommand{\liec}{\mathfrak{c}}

\newcommand{\liee}{\mathfrak{e}}
\newcommand{\lief}{\mathfrak{f}}
\newcommand{\lieg}{\mathfrak{g}}

\newcommand{\liek}{\mathfrak{k}}

\newcommand{\liem}{\mathfrak{m}}
\newcommand{\lieo}{\mathfrak{o}}

\newcommand{\lieso}{\mathfrak{so}}
\newcommand{\liesp}{\mathfrak{sp}}
\newcommand{\liet}{\mathfrak{t}}
\newcommand{\lieu}{\mathfrak{u}}
\newcommand{\liesu}{\mathfrak{su}}
\newcommand{\liew}{\mathfrak{w}}
\newcommand{\liez}{\mathfrak{z}}
\newcommand{\RP}{\ensuremath{\R\mathrm{P}}}
\newcommand{\CP}{\ensuremath{\C\mathrm{P}}}
\newcommand{\HP}{\ensuremath{\HH\mathrm{P}}}
\newcommand{\OP}{\ensuremath{\OO\mathrm{P}}}
\newcommand{\KP}{\ensuremath{\K\mathrm{P}}}

\newcommand{\frakJ}{\mathfrak{J}}

\newcommand{\bbS}{\mathbb{S}}

\newcommand{\bbV}{\mathbb{V}}

\newcommand{\Sph}{\bbS}

\newcommand{\beweis}{\begingroup\footnotesize \emph{Proof. }}
\newcommand{\beweisende}{\strut\hfill $\Box$\par\medskip\endgroup}
\newcommand{\Mengegr}[2]{\{\,#1\,{\bigr |}\,#2\,\}}

\newcommand{\wt}{\widetilde}
\newcommand{\wh}{\widehat}

\newcommand{\C}{\mathchoice {\setbox0=\hbox{$\displaystyle\rm
C$}\hbox{\hbox to0pt{\kern0.4\wd0\vrule
height0.95\ht0\hss}\box0}}{\setbox0=\hbox{$\textstyle\rm C$}\hbox{\hbox
to0pt{\kern0.4\wd0\vrule
height0.95\ht0\hss}\box0}}{\setbox0=\hbox{$\scriptstyle\rm C$}\hbox{\hbox
to0pt{\kern0.4\wd0\vrule
height0.95\ht0\hss}\box0}}{\setbox0=\hbox{$\scriptscriptstyle\rm
C$}\hbox{\hbox to0pt{\kern0.4\wd0\vrule height0.95\ht0\hss}\box0}}}

\theoremstyle{plain} 

\theorembodyfont{\rmfamily\mdseries\itshape}
\newtheorem{Def}{Definition}[section]
\newtheorem{Prop}[Def]{Proposition}
\newtheorem{Theorem}[Def]{Theorem}

\theorembodyfont{\rmfamily\mdseries\upshape}

\newtheorem{Remark}[Def]{Remark}

\begin{document}
\selectlanguage{english}

\title{Totally geodesic submanifolds of the exceptional Riemannian symmetric spaces of rank 2}
\author{Sebastian Klein$^{1}$}
\date{8 September 2008}
\maketitle
\footnotetext[1]{This work was supported by a fellowship within the Postdoc-Programme of the German Academic Exchange Service (DAAD).}
\addtocounter{footnote}{1}

\abstract{\textbf{Abstract.} 
The present article is the final part of a series on the classification of the 
totally geodesic submanifolds of the irreducible Riemannian symmetric spaces of rank 2. After this problem has been solved for the 2-Grassmannians
in my papers \cite{Klein:2007-claQ} and \cite{Klein:2007-tgG2}, and for the space \,$\SU(3)/\SO(3)$\, in Section~6 of \cite{Klein:2007-Satake},
we now solve the classification for the remaining irreducible Riemannian symmetric spaces of rank 2 and compact type: \,$\SU(6)/\Sp(3)$\,, \,$\SO(10)/\Ug(5)$\,,
\,$E_6/(\Ug(1)\cdot\Spin(10))$\,, \,$E_6/F_4$\,, \,$G_2/\SO(4)$\,, \,$\SU(3)$\,, \,$\Sp(2)$\, and \,$G_2$\,.

Similarly as for the spaces already investigated in the earlier papers, it turns out that for many of the spaces investigated here, the earlier classification
of the maximal totally geodesic submanifolds of Riemannian symmetric spaces by \textsc{Chen} and \textsc{Nagano} (\cite{Chen/Nagano:totges2-1978}, \S 9) is incomplete.
In particular, in the spaces \,$\Sp(2)$\,, \,$G_2/\SO(4)$\, and \,$G_2$\,, there exist maximal totally geodesic submanifolds, isometric to 2- or 3-dimensional spheres,
which have a ``skew'' position in the ambient space in the sense that their geodesic diameter is strictly larger than the geodesic diameter of the
ambient space. They are all missing from \cite{Chen/Nagano:totges2-1978}. 

}
\bigskip

\textbf{Author's address.} \\
Sebastian Klein \\ Department of Mathematics \\ University College Cork \\ Cork \\ Ireland \\
\texttt{mail@sebastian-klein.de}

\bigskip


\bigskip

\textbf{MS classification numbers.} 53C35 (Primary); 53C17

\section{Introduction}
\label{Se:intro}

The classification of the totally geodesic submanifolds in Riemannian symmetric spaces is an interesting and significant problem
of Riemannian geometry. Presently, I solve this problem for the irreducible Riemannian symmetric spaces of rank \,$2$\,.

The totally geodesic submanifolds of the 2-Grassmannians \,$G_2^+(\R^n)$\,, \,$G_2(\C^n)$\, and \,$G_2(\HH^n)$\, have already
been classified in my papers \cite{Klein:2007-claQ} and \cite{Klein:2007-tgG2}; moreover the totally geodesic submanifolds of \,$\SU(3)/\SO(3)$\,
have been classified in Section~6 of my paper \cite{Klein:2007-Satake}.
In the present paper I complete the classification of the totally geodesic
submanifolds in the irreducible Riemannian symmetric spaces of rank \,$2$\, (simply connected and of compact type)
by considering the remaining spaces of this kind; they are the
spaces of type I
$$ \SO(10)/\Ug(5)\,, \quad E_6/(\Ug(1)\cdot \Spin(10))\,, \quad
\SU(6)/\Sp(3)\,, \quad E_6/F_4 \quad\textrm{and}\quad G_2/\SO(4) $$
as well as the spaces of Lie group type
$$ \SU(3)\,, \quad \Sp(2) \quad\textrm{and}\quad G_2 \; ; $$
herein \,$G_2$\,, \,$F_4$\, and \,$E_6$\, denote the respective exceptional, simply connected, compact, real Lie groups.

It should be mentioned that
already \textsc{Chen} and \textsc{Nagano} gave what they claimed to be a complete classification of the isometry types of maximal totally geodesic submanifolds in all
Riemannian symmetric spaces of rank \,$2$\, in \S 9 of their paper \cite{Chen/Nagano:totges2-1978} based on their \,$(M_+,M_-)$-method.
However, as it will turn out in the present paper,
their classification is faulty also for several of the spaces under consideration here. 
In particular, in the spaces \,$\Sp(2)$\,, \,$G_2$\, and \,$G_2/\SO(4)$\,, 
there exist maximal totally geodesic submanifolds, isometric to spheres of dimension \,$2$\, or \,$3$\,,
which have a ``skew'' position in the ambient space in the sense that their geodesic diameter is strictly larger than the geodesic diameter of the ambient space;
these submanifolds are missing from Chen's and Nagano's classification. 
Also in the spaces \,$\SO(10)/\Ug(5)$\, and \,$E_6/(\Ug(1)\cdot \Spin(10))$\,, such ``skew'' totally geodesic submanifolds exist, although they are not maximal.
Moreover several other details of Chen's and Nagano's classification are incorrect. For a detailed discussion with respect to the individual spaces studied, see
the following Remarks of the present paper:

{\footnotesize
\begin{center}
\begin{tabular}{|c|c|c|c|c|c|c|c|c|}
\hline
space & $\SO(10)/\Ug(5)$ & $E_6/(\Ug(1)\cdot\Spin(10))$ & $\SU(6)/\Sp(3)$ & $E_6/F_4$ & $G_2/\SO(4)$ & $\SU(3)$ & $\Sp(2)$ & $G_2$ \\
\hline
Remark & \ref{R:EIII:DIII:CN} & \ref{R:EIII:EIII:CN} & \ref{R:EIV:AII:CN} & \ref{R:EIV:EIV:CN} & \ref{R:G2:G:CN} & \ref{R:EIV:A2:CN} & \ref{R:EIII:B2:CN} & \ref{R:G2:CN} \\
\hline
\end{tabular}
\end{center}
}

Even apart from these problems, Chen's and Nagano's investigation
is not satisfactory, as they name only the isometry type of the totally geodesic submanifolds, without giving any description of their position
in the ambient space. (Such a description can, for example, be constituted by giving explicit totally geodesic, isometric embeddings for the various congruence classes
of totally geodesic submanifolds, or at least by describing the tangent spaces of the totally geodesic submanifolds (i.e.~the Lie triple systems)
as subspaces of the tangent space of the ambient symmetric space in an explicit way.)

The usual strategy for the classification of totally geodesic submanifolds in a Riemannian symmetric space \,$M=G/K$\,,
which is used also here, is as follows.
Let \,$\lieg = \liek \oplus \liem$\, be the decomposition of the Lie algebra of \,$G$\, induced by the
symmetric structure of \,$M$\,. As it is well-known, the Lie triple systems \,$\liem'$\, in \,$\liem$\, (i.e.~the linear
subspaces \,$\liem' \subset \liem$\, which satisfy \,$[[\liem',\liem'],\liem'] \subset \liem'$\,) are in one-to-one
correspondence with the (connected, complete) totally geodesic submanifolds \,$M_{\liem'}$\, of \,$M$\,
running through the ``origin point'' \,$p_0 = eK \in M$\,,
the correspondence being that \,$M_{\liem'}$\, is characterized by \,$p_0 \in M_{\liem'}$\, and
\,$T_{p_0}M_{\liem'} = \tau(\liem')$\,, where \,$\tau: \liem \to T_{p_0}M$\, is the canonical isomorphism.

Thus the task of classifying the totally geodesic submanifolds of \,$M$\, splits into two steps: (1) To classify the Lie triple systems
in \,$\liem$\,, and (2) for each of the Lie triple systems \,$\liem'$\, found in the first step,
to construct a (connected, complete) totally geodesic
submanifold \,$M_{\liem'}$\, of \,$M$\, so that \,$p_0 \in M_{\liem'}$\, and \,$\tau^{-1}(T_{p_0}M_{\liem'}) = \liem'$\, holds.

Herein, step (1) is the one which generally poses the more significant difficulties. As an approach to accomplishing this step, we describe
in Section~\ref{Se:generallts} for an arbitrary Riemannian symmetric space \,$M$\, of compact type
relations between the roots and root spaces of \,$M$\, and the
roots resp.~root spaces of its totally geodesic submanifolds (regarded as symmetric subspaces).
These relations provide conditions which are necessary for a linear subspace \,$\liem'$\, of \,$\liem$\, to be a Lie triple system. However,
these conditions are not generally sufficient, and therefore a specific investigation needs to be made to see which of the linear subspaces of \,$\liem$\,
satisfying the conditions are in fact Lie triple systems; this investigation is
the laborious part of the proof of the classification theorems.

It should be emphasized that to carry out this investigation for a given Riemannian symmetric space \,$M$\,, it does not suffice to know the
(restricted) root system (with multiplicities) of that space, or equivalently, the action of the Jacobi operators \,$R(\,\cdot\,,v)v$\,
on the various root spaces. Rather, a full description of the curvature tensor of \,$M$\, is needed.
The well-known formula \,$R(u,v)w = -[[u,v],w]$\, relating the curvature tensor \,$R$\, of \,$M$\, to the Lie bracket of the Lie algebra \,$\lieg$\,
of the transvection group \,$G$\, of \,$M$\, lets one calculate \,$R$\, relatively easily if \,$M$\, is a classical symmetric space (then \,$\lieg$\, is a matrix Lie algebra,
with the Lie bracket being simply the commutator of matrices), but not so easily if \,$M$\, is one of the exceptional symmetric spaces, because then the explicit
description of the exceptional Lie algebra \,$\lieg$\, as a matrix algebra is too unwieldy to be useful. 

In its place, we use the description of the curvature tensor based on the root space decomposition of \,$\lieg$\, which was described in
\cite{Klein:2007-Satake}, and which permits the reconstruction of \,$R$\, using only the Satake diagram of the Riemannian symmetric space \,$M$\,.
To actually carry out the computations involved in the application of the results from \cite{Klein:2007-Satake}, we use the example implementation
of the algorithms for \textsf{Maple} also presented in that paper; this implementation is found on \texttt{http://satake.sourceforge.net}.
Whenever in the present paper, a claim is made about the evaluation of the
Lie bracket of a Lie algebra
or the curvature tensor of a Riemannian symmetric space for specific input vectors, the result has been obtained in this way. \textsf{Maple} worksheets containing
all the calculations can also be found on \texttt{http://satake.sourceforge.net}.

\enlargethispage{2em}
Certain of the spaces under investigation here are locally isometric to totally geodesic submanifolds of others; more specifically, we have
the following inclusions of totally geodesic submanifolds:
\begin{gather*}
\Sp(2)/\Z_2 \subset \SO(10)/\Ug(5) \subset E_6/(\Ug(1)\cdot \Spin(10))\,, \\
\SU(3) \subset \SU(6)/\Sp(3)\,,\quad (\SU(6)/\Sp(3))/\Z_3 \subset E_6/F_4 \qmq{and} \\
G_2/\SO(4) \subset G_2 \; . 
\end{gather*}
If \,$M$\, is a Riemannian symmetric space and \,$M' \subset M$\, a totally geodesic submanifold, then the totally geodesic submanifolds of \,$M'$\,
are exactly those totally geodesic submanifolds of \,$M$\, which are contained in \,$M'$\,. For this reason, we can obtain a classification
of the totally geodesic submanifolds of \,$M'$\, from a classification of the totally geodesic submanifolds of \,$M$\,: We just need to determine
which of the totally geodesic submanifolds of \,$M$\, are contained in \,$M'$\,. 
Thus we do not need to carry out the classification of totally geodesic submanifolds for each space under investigation here individually 
by the approach described above. Rather it suffices to do the classification for the three spaces
\,$E_6/(\Ug(1)\cdot\Spin(10))$\,, \,$E_6/F_4$\, and \,$G_2$\,, by virtue of the mentioned inclusions we then also obtain classifications
for the remaining Riemannian symmetric spaces of rank 2.

The present paper is laid out as follows: Section~\ref{Se:generallts} contains general facts on Lie triple systems, in particular on the relationship
between their (restricted) roots resp.~root spaces, and the roots resp.~root spaces of the ambient space. Section~\ref{Se:EIII} is concerned
primarily with the investigation of the Riemannian symmetric space \,$E_6/(\Ug(1)\cdot\Spin(10))$\,: In Subsection~\ref{SSe:EIII:geometry} we make
general observations about the geometry of this space; using these results we then classify the Lie triple systems of \,$E_6/(\Ug(1)\cdot\Spin(10))$\,
in Subsection~\ref{SSe:EIII:lts}, corresponding to step (1) of the classification as described above. In Subsection~\ref{SSe:EIII:tgsub} we 
describe totally geodesic embeddings for each congruence class of Lie triple systems in \,$E_6/(\Ug(1)\cdot\Spin(10))$\,, thereby completing the classification
of totally geodesic submanifolds for that space. In Subsections \ref{SSe:EIII:B2} and \ref{SSe:EIII:DIII}, we use the inclusions of totally geodesic
submanifolds \,$\Sp(2) \subset G_2(\HH^4)$\, resp.~\,$\SO(10)/\Ug(5) \subset E_6/(\Ug(1)\cdot\Spin(10))$\, to derive the classification of totally
geodesic submanifolds in \,$\Sp(2)$\, resp.~in \,$\SO(10)/\Ug(5)$\, from previous results.

Section~\ref{Se:EIV} covers the investigation of \,$E_6/F_4$\, and is structured similarly: After the introduction of basic geometric facts
on that space in Subsection~\ref{SSe:EIV:geometry}, we classify its Lie triple systems in Subsection~\ref{SSe:EIV:lts}. As a consequence of the classification
it turns out that in \,$E_6/F_4$\,, all maximal totally geodesic submanifolds are reflective. Thus we can learn the global isometry type of
the corresponding totally geodesic submanifolds from the classification of reflective submanifolds in symmetric spaces by \textsc{Leung}, \cite{Leung:reflective-1979},
as is described in Subsection~\ref{SSe:EIV:tgsub},
and do not need to construct totally geodesic embeddings in this case explicitly. In Subsections~\ref{SSe:EIV:AII}
resp.~\ref{SSe:EIV:A2} we use the inclusion \,$(\SU(6)/\Sp(3))/\Z_3 \subset E_6/F_4$\, resp.~\,$\SU(3) \subset \SU(6)/\Sp(3)$\, to derive the
classification for the space \,$\SU(6)/\Sp(3)$\, resp.~\,$\SU(3)$\,. The space \,$\SU(3)/\SO(3)$\,, whose totally geodesic submanifolds
have already been classified in Section~6 of \cite{Klein:2007-Satake}, is contained in \,$\SU(3)$\,; therefore its Lie triple systems also occur
in the present paper. Subsection~\ref{SSe:EIV:AI} gives the relationship between the types of Lie triple systems of \,$\SU(3)/\SO(3)$\, 
as defined in Section~6 of \cite{Klein:2007-Satake} and types of Lie triple systems defined here.

Section~\ref{Se:G2} then investigates the Lie group \,$G_2$\, seen as a Riemannian symmetric space. In Subsection~\ref{SSe:G2:geometry} we investigate
the geometry of this space, then we proceed in Subsection~\ref{SSe:G2:lts} to the classification of its Lie triple systems, and describe
embeddings for (most of) its totally geodesic submanifolds in Subsection~\ref{SSe:G2:tgsub}. In Subsection~\ref{SSe:G2:G} we use the inclusion
\,$G_2/\SO(4) \subset G_2$\, to derive a classification of the totally geodesic submanifolds of \,$G_2/\SO(4)$\,.

Finally, in Section~\ref{Se:summary} we give a table of the isometry types of the maximal totally geodesic submanifolds of all irreducible
Riemannian symmetric spaces of rank 2 and compact type, thereby summarizing the results of my papers \cite{Klein:2007-claQ}, \cite{Klein:2007-tgG2}, 
\cite{Klein:2007-Satake} (Section~6), as well as of the present paper.

\bigskip

The results of the present paper were obtained by me while working at the University College Cork under the advisorship of Professor J.~Berndt. 
I would like to thank him for his dedicated support and guidance, as well as his generous hospitality.

\section{General facts on Lie triple systems}
\label{Se:generallts}

In this section we suppose that \,$M=G/K$\, is any Riemannian symmetric space of compact type. 
We consider the decomposition \,$\lieg = \liek \oplus \liem$\, of the Lie algebra \,$\lieg$\, of \,$G$\, induced by the symmetric
structure of \,$M$\,. Because \,$M$\, is of compact type, the Killing form \,$\vkap: \lieg \times \lieg \to \R,\;(X,Y) \mapsto \tr(\ad(X) \circ \ad(Y))$\,
is negative definite, and therefore \,$\g{\,\cdot\,}{\,\cdot\,} := -c \cdot \vkap$\, gives rise to a Riemannian metric on \,$M$\, for arbitrary \,$c \in \R_+$\,.
In the sequel we suppose that \,$M$\, is equipped with such a Riemannian metric.%
\footnote{The dependence of the sectional curvature of \,$M$\, on the choice of the Riemannian metric is as follows: If we multiply the
Riemannian metric with some factor \,$c>0$\,, then this causes the sectional curvature function to be multiplied with \,$\tfrac1c$\,.}

Let us fix notations concerning flat subspaces, roots and root spaces of \,$M$\, (for the corresponding theory, see for example \cite{Loos:1969-2}, Section~V.2):
A linear subspace \,$\liea \subset \liem$\, is called \emph{flat} if \,$[\liea,\liea] = \{0\}$\, holds. The maximal flat subspaces of \,$\liem$\,
are all of the same dimension, called the \emph{rank} of \,$M$\, (or \,$\liem$) and denoted by \,$\rk(M)$\, or \,$\rk(\liem)$\,; they are called the
\emph{Cartan subalgebras} of \,$\liem$\,. If a Cartan subalgebra 
\,$\liea \subset \liem$\, is fixed, we put for any linear form \,$\lambda \in \liea^*$\,
$$ \liem_\lambda := \Menge{\,X \in \liem\,}{\,\forall Z \in \liea: \ad(Z)^2 X = -\lambda(Z)^2X\,} $$
and consider the \emph{(restricted) root system}
$$ \Delta(\liem,\liea) := \Menge{\,\lambda \in \liea^*\setminus \{0\}\,}{\,\liem_\lambda \neq \{0\}\,} $$
of \,$\liem$\, with respect to \,$\liea$\,. The elements of \,$\Delta(\liem,\liea)$\, are called \emph{(restricted) roots} of \,$\liem$\, with respect to \,$\liea$\,,
for \,$\lambda \in \Delta(\liem,\liea)$\, the subspace \,$\liem_\lambda$\, is called the \emph{root space} corresponding to \,$\lambda$\,, and \,$n_{\lambda} := \dim(\liem_{\lambda})$\,
is called the \emph{multiplicity} of the root \,$\lambda$\,. If we fix a system of positive roots \,$\Delta_+ \subset \Delta(\liem,\liea)$\,
(i.e.~we have \,$\Delta_+ \dot{\cup} (-\Delta_+) = \Delta(\liem,\liea)$\,), 
we obtain the \emph{(restricted) root space decomposition} of \,$\liem$\,:
\begin{equation}
\label{eq:roots:mdecomp}
\liem = \liea \;\oplus\; \bigoplus_{\lambda \in \Delta_+} \liem_\lambda \;.
\end{equation}
The \emph{Weyl group} \,$W(\liem,\liea)$\, is the transformation group on \,$\liea$\, generated by the reflections in the hyperplanes \,$\Menge{v \in \liea}{\lambda(v)=0}$\,
(where \,$\lambda$\, runs through \,$\Delta(\liem,\liea)$\,);
it can be shown that the root system \,$\Delta(\liem,\liea)$\, is invariant under the action of \,$W(\liem,\liea)$\,. 

\bigskip

Let us now consider a Lie triple system \,$\liem' \subset \liem$\,, i.e.~\,$\liem'$\, is a linear subspace of \,$\liem$\, so that
\,$[\,[\liem',\liem'] \,,\, \liem'\,] \subset \liem'$\, holds. In spite of the fact that the symmetric space corresponding to \,$\liem'$\,
does not need to be of compact type (it can contain Euclidean factors), it is easily seen that the usual statements of the root space theory
for symmetric spaces of compact type carry over to \,$\liem'$\,, see \cite{Klein:2007-claQ}. 

More specifically, the maximal flat subspaces of \,$\liem'$\, are all of the same dimension
(again called the \emph{rank} of \,$\liem'$\,), and they are again called the \emph{Cartan subalgebras} of \,$\liem'$\,. For any Cartan subalgebra \,$\liea'$\,
of \,$\liem'$\,, there exists a Cartan subalgebra \,$\liea$\, of \,$\liem$\, so that \,$\liea' = \liea \cap \liem'$\, holds. With respect to any
Cartan subalgebra \,$\liea'$\, of \,$\liem'$\, we have a root system \,$\Delta(\liem',\liea')$\, (defined analogously as for \,$\liem$\,) 
and the corresponding root space decomposition
\begin{equation}
\label{eq:roots:m'decomp}
\liem' = \liea' \;\oplus\; \bigoplus_{\alpha \in \Delta_+(\liem',\liea')} \liem'_{\alpha}
\end{equation}
(with a system of positive roots \,$\Delta_+(\liem',\liea') \subset \Delta(\liem',\liea')$\,); we also again call \,$n_{\alpha}' := \dim(\liem_{\alpha}')$\,
the multiplicity of \,$\alpha \in \Delta(\liem',\liea')$\,. 
\,$\Delta(\liem',\liea')$\, is again invariant under the action of the corresponding Weyl group \,$W(\liem',\liea')$\,. 
It should be noted, however,
that in the case where a Euclidean factor is present in \,$\liem'$\,, \,$\Delta(\liem',\liea')$\, does not span \,$(\liea')^*$\,.

The following proposition describes the relation between the root space decompositions~\eqref{eq:roots:m'decomp} of \,$\liem'$\, and 
\eqref{eq:roots:mdecomp} of \,$\liem$\,. In particular, it shows the 
extent to which the the position of the individual root spaces \,$\liem'_\alpha$\, of \,$\liem'$\, is adapted to the root space
decomposition~\eqref{eq:roots:mdecomp} of the ambient space \,$\liem$\,. These relations will play a fundamental role in our
classification of the Lie triple systems in the Riemannian symmetric spaces of rank \,$2$\,. 

\begin{Prop}
\label{P:cla:subroots:subroots-neu}
Let \,$\liea'$\, be a Cartan subalgebra of \,$\liem'$\,, and let \,$\liea$\, be a Cartan subalgebra of \,$\liem$\, so that \,$\liea' = \liea \cap \liem'$\, holds.
\begin{enumerate}
\item
The roots resp.~root spaces of \,$\liem'$\, and of \,$\liem$\, are related in the following way:
\begin{gather}
\label{eq:cla:subroots:subroots-neu:toshow-Delta}
\Delta(\liem',\liea') \;\subset\; \;\bigr\{\;\lambda|\liea'\; \bigr| \;\lambda \in \Delta(\liem,\liea), \lambda|\liea' \neq 0\; \bigr\} \;.\\
\label{eq:cla:subroots:subroots-neu:toshow-liemalpha}
\textstyle
\forall\alpha \in \Delta(\liem',\liea')\;:\; \liem_\alpha' = \left( \bigoplus_{\substack{\lambda \in \Delta(\liem,\liea) \\ \lambda|\liea' = \alpha}} \liem_\lambda \right) \;\cap\; \liem' \; .
\end{gather}
In particular, if \,$\lambda\in\Delta(\liem,\liea)$\, satisfies \,$\lambda|\liea'=0$\,, then \,$\liem'$\, is orthogonal to \,$\liem_\lambda$\,. 
\item
We have \,$\rk(\liem') = \rk(\liem)$\, if and only if \,$\liea' = \liea$\, holds. If this is the case, then we have
\begin{equation}
\label{eq:cla:subroots:subroots-neu:c}
\Delta(\liem',\liea') \subset \Delta(\liem,\liea) \;,\quad \forall \alpha \in \Delta(\liem',\liea') : \liem'_\alpha = \liem_\alpha \cap \liem' \; . 
\end{equation}
\end{enumerate}
\end{Prop}

\beweis
See \cite{Klein:2007-claQ}, the proof of Proposition~2.1.
\beweisende

For the remainder of the section, we fix a Cartan subalgebra \,$\liea'$\, of \,$\liem'$\,, and let \,$\liea$\, be 
any Cartan subalgebra of \,$\liem$\, so that \,$\liea' = \liea \cap \liem'$\, holds.

\begin{Def}
\label{D:cla:subroots:Elemcomp}
Let \,$\alpha \in \Delta(\liem',\liea')$\, be given.
Recall that by Proposition~\ref{P:cla:subroots:subroots-neu}(a) there exists at least one root \,$\lambda \in \Delta(\liem,\liea)$\, with
\,$\lambda|\liea' = \alpha$\,. We call \,$\alpha$\, 
\begin{enumerate}
\item \emph{elementary}, if there exists only one root \,$\lambda \in \Delta(\liem,\liea)$\, with \,$\lambda|\liea' = \alpha$\,;
\item \emph{composite}, if there exist at least two different roots \,$\lambda, \mu \in \Delta(\liem,\liea)$\, with \,$\lambda|\liea' = \alpha = \mu|\liea'$\,.
\end{enumerate}
\end{Def}

Elementary roots play a special role:
If \,$\alpha \in \Delta(\liem',\liea')$\, is elementary, then the root space \,$\liem_\alpha'$\, 
is contained in the root space \,$\liem_\lambda$\,, where \,$\lambda \in \Delta(\liem,\liea)$\, is the unique root with \,$\lambda|\liea' = \alpha$\,.
As we will see in Proposition~\ref{P:cla:subroots:Comp} below, this property
causes restrictions for the possible positions (in relation to \,$\liea'$\,) of \,$\lambda$\,. The exploitation of these restrictions will play
an important role in the classification of the rank \,$1$\, Lie triple systems in the rank \,$2$\, spaces under investigation.

It should also be mentioned that in the case \,$\rk(\liem') = \rk(\liem)$\, we have \,$\liea' = \liea$\,, and therefore in that case
every \,$\alpha \in \Delta(\liem',\liea')$\, is elementary (compare Proposition~\ref{P:cla:subroots:subroots-neu}(b)).

For any linear form \,$\lambda \in \liea^*$\, we now denote by \,$\lambda^\sharp$\, the Riesz vector
corresponding to \,$\lambda$\,, i.e.~the vector \,$\lambda^\sharp \in \liea$\, characterized by \,$\g{\,\cdot\,}{\lambda^\sharp} = \lambda$\,. Here \,$\g{\,\cdot\,}{\,\cdot\,} = -c\cdot \vkap$\, is again the inner product obtained from the Killing form \,$\vkap$\,
of \,$\lieg$\,. 

\begin{Prop}
\label{P:cla:subroots:Comp}
Let \,$\alpha \in \Delta(\liem',\liea')$\, be given. 
\begin{enumerate}
\item
If \,$\alpha$\, is elementary 
and \,$\lambda \in \Delta(\liem,\liea)$\, is the unique root with \,$\lambda|\liea' = \alpha$\,, then we have \,$ \lambda^\sharp \in \liea'$\,. 
\item
If \,$\alpha$\, is composite and \,$\lambda, \mu \in \Delta(\liem,\liea)$\, are two different roots with
\,$\lambda|\liea' = \alpha = \mu|\liea'$\,, then \,$\lambda^\sharp - \mu^\sharp$\, is orthogonal to \,$\liea'$\,. 
\end{enumerate}
\end{Prop}

\beweis
For (a) see \cite{Klein:2007-claQ}, the proof of Proposition~2.3(a). (b) is obvious.
\beweisende

\begin{Prop}
\label{P:cla:skew}
Suppose that \,$\alpha \in \Delta(\liem',\liea')$\, is a composite root such that there exist precisely two roots \,$\lambda,\mu \in \Delta(\liem,\liea)$\, 
with \,$\lambda|\liea' = \alpha = \mu|\liea'$\,. Further suppose that \,$\alpha^\sharp$\, can be written as a linear combination \,$\alpha^\sharp = a\,\lambda^\sharp
+ b\,\mu^\sharp$\, with non-zero \,$a,b \in \R$\,. 

Then we have \,$a,b > 0$\,, and
there exists a linear subspace \,$\liem_{\lambda}' \subset \liem_\lambda$\, and an isometric linear map \,$\Phi: \liem_{\lambda}' \to \liem_\mu$\,
so that
\begin{equation}
\label{eq:cla:skew:skew}
\liem_{\alpha}' = \Menge{x+ \sqrt{\tfrac{b}{a}}\, \Phi(x)}{x \in \liem_{\lambda}'}
\end{equation}
holds. In particular we have \,$n_{\alpha}' \leq \min\{n_\lambda,n_\mu\}$\,. 
\end{Prop}

\beweis
See \cite{Klein:2007-tgG2}, the proof of Proposition~2.4.
\beweisende

We mention one important principle for the construction of Lie triple systems with only elementary roots.

\begin{Def}
\label{D:cla:closed}
A subset \,$\Delta' \subset \Delta(\liem,\liea)$\, is called a \emph{closed root subsystem} of \,$\Delta(\liem,\liea)$\, if for every \,$\lambda\in\Delta'$\,
we also have \,$-\lambda\in\Delta'$\,, and if for every \,$\lambda,\mu\in\Delta'$\, with \,$\lambda+\mu\in\Delta(\liem,\liea)$\, we have \,$\lambda+\mu\in\Delta'$\,.
\end{Def}

\begin{Prop}
\label{P:cla:assoclts}
Let \,$\Delta'$\, be a closed root subsystem of \,$\Delta(\liem,\liea)$\,, and let \,$\Delta_+'$\, be a positive root system of \,$\Delta'$\,.
Then \,$\liem' := \spn_{\R}\Menge{\lambda^\sharp}{\lambda\in\Delta'} \oplus \bigoplus_{\lambda\in\Delta'_+} \liem_\lambda$\, is a Lie triple system in \,$\liem$\,. \,$\liem'$\, is called
the \emph{Lie triple system associated to \,$\Delta'$\,}.
\end{Prop}

\beweis
This follows immediately from the fact that for any \,$\lambda,\mu\in\Delta(\liem,\liea)\cup\{0\}$\, we have
$$ [\liem_\lambda,\liem_\mu] \subset \liek_{\lambda+\mu} \oplus \liek_{\lambda-\mu} \qmq{and}
[\liek_\lambda,\liem_\mu] \subset \liem_{\lambda+\mu} \oplus \liem_{\lambda-\mu} \;, $$
see \cite{Loos:1969-2}, Proposition~VI.1.4c, p.~60. Here \,$\liek_\lambda$\, denotes the root space of \,$\liek$\, corresponding to \,$\lambda \in \Delta(\liem,\liea)$\,.
\beweisende

The isotropy group \,$K$\, of the symmetric space \,$M$\,
acts on \,$\liem$\, via the adjoint representation, i.e.~by
\,$K \times \liem \to \liem,\; (g,v) \mapsto \Ad(g)v$\,; this action is called the \emph{isotropy action}.
In the investigation of Riemannian symmetric spaces, the orbits of this action play an important role. In the case 
of spaces of rank 2, they form a 1-parameter family, which can be parametrized in the following way (generalizing the approach that was used for the 2-Grassmannians
in \cite{Klein:2007-claQ} and \cite{Klein:2007-tgG2}):

We suppose that \,$M$\, is of rank \,$2$\,, and fix a Weyl chamber \,$\liec$\, in \,$\liea$\,. We denote the two rays in \,$\liea$\, delineating this
Weyl chamber by \,$R_1$\, and \,$R_2$\,; in the case where \,$\Delta(\liem,\liea)$\, contains roots of different length (i.e.~the root system \,$\Delta(\liem,\liea)$\,
is of one of the types \,$B_2$\,, \,$BC_2$\, or \,$G_2$\,), we suppose that \,$R_1$\, points into the direction of one of the shorter roots. Let \,$\vi_{\max}$\,
be the angle between \,$R_1$\, and \,$R_2$\,; \,$\vi_{\max}$\, equals \,$\tfrac\pi3$\,, \,$\tfrac\pi4$\,, \,$\tfrac\pi4$\, or \,$\tfrac\pi6$\,, according to
whether \,$\Delta(\liem,\liea)$\, is of type \,$A_2$\,, \,$B_2$\,, \,$BC_2$\, or \,$G_2$\,, respectively.

Any given \,$v \in \liem\setminus\{0\}$\, is congruent under the isotropy action to one and only one vector \,$v_0 \in \overline{\liec}$\,,
and we denote the angle between \,$R_1$\, and \,$v_0$\, by \,$\vi(v)$\,. In this way we obtain a continuous function 
\,$\vi: \liem \setminus \{0\} \to [0,\vi_{\max}]$\,. Two vectors \,$v_1,v_2\in\liem'$\, with \,$\|v_1\| = \|v_2\| \neq 0$\, are congruent under
the isotropy action if and only if \,$\vi(v_1) = \vi(v_2)$\, holds. We call the value \,$\vi(v)$\, the \emph{isotropy angle} of a vector
\,$v\in\liem \setminus \{0\}$\,. 

Notice that if \,$\liem'$\, is a Lie triple system of \,$\liem$\, of rank \,$1$\,, then \,$\vi$\, is constant on \,$\liem' \setminus \{0\}$\,, and 
Proposition~\ref{P:cla:subroots:Comp} shows that there are only finitely many \,$t \in [0,\vi_{\max}]$\, so that there exists a Lie triple
system \,$\liem' \subset \liem$\, of rank \,$1$\, with \,$\vi|(\liem' \setminus \{0\}) = t$\, and \,$\dim(\liem') \geq 2$\,. We will call the value \,$t$\,
for such a Lie triple system \,$\liem'$\, the \emph{isotropy angle} of \,$\liem'$\,. On the other hand,
if \,$\liem'$\, is of rank \,$2$\,, then we have \,$\vi(\liem' \setminus \{0\}) = [0,\vi_{\max}]$\,. 

\section[The symmetric spaces \,$E_6/(\Ug(1)\cdot\Spin(10))$\,, \,$\Sp(2)$\, and \,$\SO(10)/\Ug(5)$\,]{The symmetric spaces \,$\boldsymbol{E_6/(\Ug(1)\cdot\Spin(10))}$\,, \,$\boldsymbol{\Sp(2)}$\, and \,$\boldsymbol{\SO(10)/\Ug(5)}$\,}
\label{Se:EIII}

\subsection[The geometry of \,$E_6/(\Ug(1)\cdot\Spin(10))$\,]{The geometry of \,$\boldsymbol{E_6/(\Ug(1)\cdot\Spin(10))}$\,}
\label{SSe:EIII:geometry}

In the present section we will study
the Hermitian symmetric space \,$\EIII:=\Esix/(\Ug(1)\cdot\Spin(10))$\,, which has the Satake diagram
\begin{center}
\begin{minipage}{5cm}
\medskip
\xymatrix@=.4cm{
\\ 
{\displaystyle \mathop{\circ}^{1}} \ar@{-}[r]<-.8ex> \ar@{<->}@/^1.5pc/[0,4]
        & \displaystyle \mathop{\bullet}^{3} \ar@{-}[r]<-.8ex>
        & \displaystyle \mathop{\bullet}^{4} \ar@{-}[r]<-.8ex> \ar@{-}[d]
        & \displaystyle \mathop{\bullet}^{5} \ar@{-}[r]<-.8ex>
        & {\displaystyle \mathop{\circ}^{6}} \\
&& {\displaystyle \mathop{\circ}_{2}} &&
}
\medskip
\end{minipage}
\end{center}

We consider the Lie algebra \,$\lieg := \liee_6$\,
of the transvection group \,$\Esix$\, of \,$\EIII$\,, and the splitting \,$\lieg = \liek \oplus \liem$\, induced by the symmetric structure of \,$\EIII$\,.
Here \,$\liek = \R \oplus \lieso(10)$\, is the Lie algebra of the isotropy group of \,$\EIII$\,, and \,$\liem$\, is isomorphic to the tangent space
of \,$\EIII$\, in the origin. The \,$E_6$-invariant Riemannian metric on \,$\EIII$\, induces an \,$\Ad(\Ug(1)\cdot\Spin(10))$-invariant Riemannian
metric on \,$\liem$\,. As was explained in Section~\ref{Se:generallts}, this metric is only unique up to a factor; we choose the factor in such a way
that the shortest restricted roots of \,$\EIII$\, (see below) have length \,$1$\,. 

\paragraph{The root space decomposition.} Let \,$\liet$\, be a Cartan subalgebra of \,$\lieg$\, which is maximally non-compact, i.e.~\,$\liet$\, is
chosen such that the flat subspace
\,$\liea := \liet \cap \liem$\, of \,$\liem$\, is of the maximal dimension \,$2$\,, and hence a Cartan subalgebra of \,$\liem$\,. 
Then we consider the root system \,$\Delta^\lieg \subset \liet^*$\, of \,$\lieg$\, with respect to \,$\liet$\,,
as well as the restricted root system \,$\Delta \subset \liea^*$\, of the symmetric space \,$\EIII$\, with respect to \,$\liea$\,. 
\,$\EIII$\, has the restricted Dynkin diagram with multiplicities  $\xymatrix@=.4cm{ \mathop{\bullet}^6 \ar@{<=>}[r] & \mathop{\bullet\hspace{-.29cm}\bigcirc}^{8[1]} }$,
in other words: its restricted root system \,$\Delta$\, is of type \,$BC_2$\,, i.e.~we have
\,$\Delta = \{\pm\lambda_1, \pm\lambda_2, \pm\lambda_3, \pm\lambda_4, \pm 2\lambda_1, \pm 2\lambda_2\}$\,, where
\,$(\lambda_1,\lambda_3)$\, is a system of simple roots of \,$\Delta$\,, these two roots are at an angle of \,$\tfrac{3}{4}\,\pi$\, 
with \,$\lambda_3$\, being the longer of the two, and we have \,$\lambda_2 = \lambda_1+\lambda_3$\,, \,$\lambda_4 = 2\lambda_1 + \lambda_3$\,. 
Moreover, the restricted roots have the following multiplicities: 
\,$n_{\lambda_1} = n_{\lambda_2} = 8$\,, \,$n_{\lambda_3} = n_{\lambda_4} = 6$\, and \,$n_{2\lambda_1} = n_{2\lambda_2} = 1$\,.
\,$\Delta$\, has the following graphical representation: 
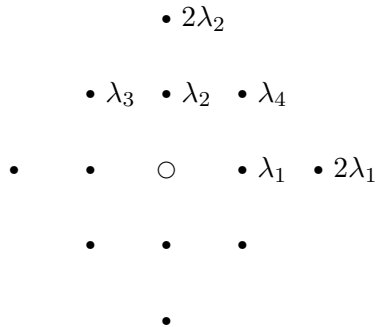
\begin{center}
\begin{minipage}{5cm}
\begin{center}
\strut \\[-.3cm]
\setlength{\unitlength}{1.0cm}
\begin{picture}(2,5)
\put(1,3){\circle{0.2}}
\put(2,3){\circle*{0.1}}        
\put(3,3){\circle*{0.1}}        
\put(2,2){\circle*{0.1}}        
\put(1,2){\circle*{0.1}}        
\put(1,1){\circle*{0.1}}        
\put(0,2){\circle*{0.1}}        
\put(0,3){\circle*{0.1}}        
\put(-1,3){\circle*{0.1}}        
\put(0,4){\circle*{0.1}}        
\put(1,4){\circle*{0.1}}        
\put(1,5){\circle*{0.1}}        
\put(2,4){\circle*{0.1}}        
\put(2.2,3.9){{$\lambda_4$}}
\put(2.2,2.9){{$\lambda_1$}}
\put(3.2,2.9){{$2\lambda_1$}}
\put(0.2,3.9){{$\lambda_3$}}
\put(1.2,3.9){{$\lambda_2$}}
\put(1.2,4.9){{$2\lambda_2$}}
\end{picture}
\end{center}
\end{minipage}
\end{center}
\vspace{-.3cm}

To be able to apply the results from \cite{Klein:2007-Satake} and the corresponding computer package for the calculation of the curvature tensor of \,$\EIII$\,,
we need to describe the relationship between the restricted roots of the symmetric space \,$\EIII$\, and the (non-restricted) roots of the Lie algebra \,$\liee_6$\,.
For this purpose, we order the simple roots of \,$\liee_6$\, as they are numbered in the Satake diagram of \,$\EIII$\, given above.
Then we label the 36 positive roots of  \,$\liee_6$\, by \,$\alpha_1,\dotsc,\alpha_{36}$\, 
in the order in which they are produced by Algorithm~(R) in Section~2 of \cite{Klein:2007-Satake} based on this ordering of the simple roots. 
It turns out that \,$\alpha_1,\dotsc,\alpha_{36}$\, have the following coordinates with respect to the simple roots of \,$\liee_6$\, ordered as before:

\medskip
\begin{tabular}{|c|c|}
\hline
\,$\alpha_1$\, & $(1,0,0,0,0,0)$ \\
\,$\alpha_2$\, & $(0,1,0,0,0,0)$ \\
\,$\alpha_3$\, & $(0,0,1,0,0,0)$ \\
\,$\alpha_4$\, & $(0,0,0,1,0,0)$ \\
\,$\alpha_5$\, & $(0,0,0,0,1,0)$ \\
\,$\alpha_6$\, & $(0,0,0,0,0,1)$ \\
\,$\alpha_7$\, & $(1,0,1,0,0,0)$ \\
\,$\alpha_8$\, & $(0,1,0,1,0,0)$ \\
\,$\alpha_9$\, & $(0,0,1,1,0,0)$ \\
\hline
\end{tabular}
\begin{tabular}{|c|c|}
\hline
\,$\alpha_{10}$\, & $(0,0,0,1,1,0)$ \\
\,$\alpha_{11}$\, & $(0,0,0,0,1,1)$ \\
\,$\alpha_{12}$\, & $(1,0,1,1,0,0)$ \\
\,$\alpha_{13}$\, & $(0,1,1,1,0,0)$ \\
\,$\alpha_{14}$\, & $(0,1,0,1,1,0)$ \\
\,$\alpha_{15}$\, & $(0,0,1,1,1,0)$ \\
\,$\alpha_{16}$\, & $(0,0,0,1,1,1)$ \\
\,$\alpha_{17}$\, & $(1,1,1,1,0,0)$ \\
\,$\alpha_{18}$\, & $(1,0,1,1,1,0)$ \\
\hline
\end{tabular}
\begin{tabular}{|c|c|}
\hline
\,$\alpha_{19}$\, & $(0,1,1,1,1,0)$ \\
\,$\alpha_{20}$\, & $(0,1,0,1,1,1)$ \\
\,$\alpha_{21}$\, & $(0,0,1,1,1,1)$ \\
\,$\alpha_{22}$\, & $(1,1,1,1,1,0)$ \\
\,$\alpha_{23}$\, & $(1,0,1,1,1,1)$ \\
\,$\alpha_{24}$\, & $(0,1,1,1,1,1)$ \\
\,$\alpha_{25}$\, & $(0,1,1,2,1,0)$ \\
\,$\alpha_{26}$\, & $(1,1,1,1,1,1)$ \\
\,$\alpha_{27}$\, & $(1,1,1,2,1,0)$ \\
\hline
\end{tabular}
\begin{tabular}{|c|c|}
\hline
\,$\alpha_{28}$\, & $(0,1,1,2,1,1)$ \\
\,$\alpha_{29}$\, & $(1,1,1,2,1,1)$ \\
\,$\alpha_{30}$\, & $(1,1,2,2,1,0)$ \\
\,$\alpha_{31}$\, & $(0,1,1,2,2,1)$ \\
\,$\alpha_{32}$\, & $(1,1,1,2,2,1)$ \\
\,$\alpha_{33}$\, & $(1,1,2,2,1,1)$ \\
\,$\alpha_{34}$\, & $(1,1,2,2,2,1)$ \\
\,$\alpha_{35}$\, & $(1,1,2,3,2,1)$ \\
\,$\alpha_{36}$\, & $(1,2,2,3,2,1)$ \\
\hline
\end{tabular}
\medskip

To find out which restricted root of \,$\EIII$\, corresponds to each root of \,$\liee_6$\,,
we tabulate the orbits of the action of \,$\sigma$\, on the root system \,$\Delta^\lieg$\,, and the restricted root of \,$\EIII$\, corresponding to each orbit
(compare Section~4 of \cite{Klein:2007-Satake}):

\medskip
{\footnotesize
\begin{center}
\begin{tabular}{|c||c|c|c|c|c|}
\hline
orbit & $\{\alpha_1,-\alpha_{21}\}$ & $\{\alpha_6,-\alpha_{18}\}$ & $\{\alpha_7,-\alpha_{16}\}$ & $\{\alpha_{11},-\alpha_{12}\}$ & $\{\alpha_{23},-\alpha_{23}\}$ \\
\hline
corresp.~restr.~root & $\lambda_1$ & $\lambda_1$ & $\lambda_1$ & $\lambda_1$ & $2\lambda_1$ \\
\hline
\end{tabular}

\medskip

\begin{center}
\begin{tabular}{|c||c|c|c|c|c|}
\hline
orbit & $\{\alpha_{17},-\alpha_{31}\}$ & $\{\alpha_{20},-\alpha_{30}\}$ & $\{\alpha_{22},-\alpha_{28}\}$ & $\{\alpha_{24},-\alpha_{27}\}$ & $\{\alpha_{36},-\alpha_{36}\}$ \\
\hline
corresp.~restr.~root & $\lambda_2$ & $\lambda_2$ & $\lambda_2$ & $\lambda_2$ & $2\lambda_2$ \\
\hline
\end{tabular}

\medskip

\begin{tabular}{|c||c|c|c|c|c|c|}
\hline
orbit & $\{\alpha_{2},-\alpha_{25}\}$ & $\{\alpha_{8},-\alpha_{19}\}$ & $\{\alpha_{13},-\alpha_{14}\}$ & $\{\alpha_{26},-\alpha_{35}\}$ & $\{\alpha_{29},-\alpha_{34}\}$ & $\{\alpha_{32},-\alpha_{33}\}$  \\ 
\hline
corresp.~restr.~root & $\lambda_3$ & $\lambda_3$ & $\lambda_3$ & $\lambda_4$ & $\lambda_4$ & $\lambda_4$ \\
\hline
\end{tabular}
\end{center}

\end{center}

}

\medskip

Moreover, we have \,$\sigma(\alpha_k)=\alpha_k$\, for \,$k\in\{3,4,5,9,10,15\}$\,.

Using the notations of \cite{Klein:2007-Satake}, Proposition~5.2(a) we now put for \,$c_1,\dotsc,c_4 \in \C$\, and \,$t\in\R$\,, where \,$A$\, denotes either
of the letters \,$K$\, and \,$M$\,:
\begin{align*}
A_{\lambda_1}(c_1,c_2,c_3,c_4) & := A_{\alpha_{1}}(c_1) + A_{\alpha_{6}}(c_2) + A_{\alpha_{7}}(c_3) + A_{\alpha_{11}}(c_4) \; , \\
A_{2\lambda_1}(t) & := \wt{A}_{\alpha_{23}}(t) \; , \\
A_{\lambda_2}(c_1,c_2,c_3,c_4) & := A_{\alpha_{17}}(c_1) + A_{\alpha_{20}}(c_2) + A_{\alpha_{22}}(c_3) + A_{\alpha_{24}}(c_4) \; , \\
A_{2\lambda_2}(t) & := \wt{A}_{\alpha_{36}}(t) \; , \\
A_{\lambda_3}(c_1,c_2,c_3) & := A_{\alpha_{2}}(c_1) + A_{\alpha_{8}}(c_2) + A_{\alpha_{13}}(c_3) \; , \\
A_{\lambda_4}(c_1,c_2,c_3) & := A_{\alpha_{26}}(c_1) + A_{\alpha_{29}}(c_2) + A_{\alpha_{32}}(c_3) \; .
\end{align*}
Then we have \,$\liem_{\lambda_k} = M_{\lambda_k}(\C,\C,\C,\C)$\, and \,$\liem_{2\lambda_k} = M_{2\lambda_k}(\R)$\, for \,$k\in\{1,2\}$\,, and
\,$\liem_{\lambda_k} = M_{\lambda_k}(\C,\C,\C)$\, for \,$k\in\{3,4\}$\,. 

\paragraph{The action of the isotropy group.} We next look at the isotropy action of \,$\EIII$\,. Regarding it, we use the notations introduced at the end
of Section~\ref{Se:generallts}, in particular we have the continuous function \,$\vi: \liem \setminus \{0\} \to [0,\tfrac\pi4]$\, parametrizing the
orbits of the isotropy action. For the elements of the closure \,$\overline{\liec}$\, of the positive Weyl chamber \,$\liec := \Menge{v\in\liea}{\lambda_1(v)\geq 0,
\lambda_3(v)\geq 0}$\, we can explicitly describe the relation to their isotropy angle: \,$(\lambda_2^\sharp,\lambda_1^\sharp)$\, is an orthonormal basis
of \,$\liea$\, so that with \,$v_t := \cos(t)\lambda_2^\sharp + \sin(t)\lambda_1^\sharp$\, we have
\begin{equation}
\label{eq:EIII:isotropy:liec}
\overline{\liec} = \Menge{s\cdot v_t}{t\in[0,\tfrac\pi4],s\in\R_{\geq 0}}\;, 
\end{equation}
and because the Weyl chamber \,$\liec$\, is bordered
by the two vectors \,$v_0 = \lambda_2^\sharp$\, with \,$\vi(v_0)=0$\, and \,$v_{\pi/4}=\tfrac{1}{\sqrt{2}}\,\lambda_4^\sharp$\, with \,$\vi(v_{\pi/4})=\tfrac\pi4$\,, 
we have
\begin{equation}
\label{eq:EIII:isotropy:vivt}
\vi(s\cdot v_t)=t \qmq{for all} t\in[0,\tfrac\pi4]\,, s\in\R_+ \; . 
\end{equation}

The action of
the subgroup \,$K_0$\, of \,$K$\, whose Lie algebra is the centralizer \,$\liek_0 := \Menge{X\in\liek}{[X,\liea]=0}$\, of \,$\liea$\, in \,$\liek$\,
leaves the restricted root spaces \,$\liem_\lambda$\, invariant. The Dynkin diagram of \,$\liek_0$\, is given by the black roots in the Satake diagram
of \,$\EIII$\, (see above), therefore we have \,$\liek_0 = (\liet\cap \liek) \oplus \bigoplus K_{\alpha_\ell}(\C)$\,, where the sum runs over all
those roots \,$\alpha_\ell$\, of \,$\liee_6$\, with \,$\sigma(\alpha_\ell) = \alpha_\ell$\,, i.e.~\,$\ell\in \{3,4,5,9,10,15\}$\,. 
Because of this and the fact that \,$\dim(\liet\cap\liek)=4$\, holds,
it follows that \,$\liek_0$\, is isomorphic to \,$\lieu(4)$\,, and hence \,$K_0$\, is locally isomorphic to \,$\Ug(4)$\,. 

By using the \textsf{Maple} implementation to look at the adjoint action of \,$\liek_0$\, on the root spaces \,$\liem_\lambda$\,, we can describe the action of \,$K_0$\, on the root spaces
in more detail:

\begin{Prop}
\label{P:EIII:isotropy}
For \,$k\in\{1,2\}$\, the action of \,$K_0$\, on \,$\liem_{\lambda_k}$\, is locally equivalent to the vector representation of \,$\Ug(4)$\,, this means that
if we denote by \,$\vi$\, the linear isometry
$$ \vi: \C^4 \to \liem_{\lambda_k},\; (c_1,c_2,c_3,c_4) \mapsto M_{\lambda_k}(c_1,c_2,c_3,c_4) \; , $$
there exists a local isomorphism of Lie groups \,$\Phi: \Ug(4) \to K_0$\, so that the following diagram commutes:
\begin{equation*}
\begin{minipage}{5cm}
\begin{xy}
\xymatrix{
\Ug(4)\times \C^4 \ar[r]^{\Phi\times\vi} \ar[d] & K_0 \times \liem_{\lambda_k} \ar[d]^{\Ad} \\
\C^4 \ar[r]_{\vi} & \liem_{\lambda_k} \;,
}
\end{xy}
\end{minipage}
\end{equation*}
where the left vertical arrow represents the canonical action of \,$\Ug(4)$\, on \,$\C^4$\,.

Moreover, if we fix \,$v\in\liem_{\lambda_k} \setminus \{0\}$\,, then the Lie subgroup \,$U' := \Menge{B\in \Ug(4)}{B(\vi^{-1}v) = \vi^{-1}v}$\, of \,$\Ug(4)$\,
is isomorphic to \,$\Ug(3)$\,, and hence the Lie subgroup \,$K_0' := \Menge{g\in K_0}{\Ad(g)v=v}$\, of \,$K_0$\, is locally isomorphic to \,$\Ug(3)$\,.
For \,$\ell \in \{3,4\}$\,, the action of \,$K_0'$\, on \,$\liem_{\lambda_\ell}$\, is locally equivalent to the vector representation of \,$\Ug(3)$\,, 
i.e.~with the linear isometry
$$ \psi: \C^3 \to \liem_{\ell},\; (c_1,c_2,c_3) \mapsto M_{\lambda_\ell}(c_1,c_2,c_3) $$
there exists a local isomorphism of Lie groups \,$\Psi: \Ug(3) \to K_0'$\, so that the following diagram commutes:
\begin{equation*}
\begin{minipage}{5cm}
\begin{xy}
\xymatrix{
\Ug(3)\times \C^3 \ar[r]^{\Psi\times\psi} \ar[d] & K_0' \times \liem_{\lambda_\ell} \ar[d]^{\Ad} \\
\C^3 \ar[r]_{\psi} & \liem_{\lambda_\ell} \;,
}
\end{xy}
\end{minipage}
\end{equation*}
where the left vertical arrow represents the canonical action of \,$\Ug(3)$\, on \,$\C^3$\,.

In particular we see that \,$\Ad(K_0)$\, acts ``jointly transitively'' on the unit spheres in \,$\liem_{\lambda_k}$\, and \,$\liem_{\lambda_\ell}$\,
in the sense that for any given \,$v_1,v_2\in \liem_{\lambda_k}$\, and \,$w_1,w_2\in\liem_{\lambda_\ell}$\, with \,$\|v_1\| = \|v_2\|$\, and \,$\|w_1\| = \|w_2\|$\,
there exists \,$g\in K_0$\, with \,$\Ad(g)v_1 = v_2$\, and \,$\Ad(g)w_1 = w_2$\,. 

Finally, we note that the linear isometries
$$ \liem_{\lambda_1} \to \liem_{\lambda_2}, \; M_{\lambda_1}(c_1,c_2,c_3,c_4) \mapsto M_{\lambda_2}(c_2,c_1,c_4,c_3) $$
and 
$$ \liem_{\lambda_3} \to \liem_{\lambda_4},\; M_{\lambda_3}(c_1,c_2,c_3) \mapsto M_{\lambda_4}(c_1,c_2,c_3) $$
commute with the action of \,$\Ad(K_0)$\, on the respective root spaces.
\end{Prop}

\paragraph{The complex structure of \,$\boldsymbol{\EIII}$\,.} \,$\EIII$\, is a Hermitian symmetric space; the action of its complex structure \,$J$\, on \,$\liem$\,
is given by \,$J|\liem = \ad(j)|\liem$\,, where \,$j$\, is a element of the center \,$\liez(\liek)$\, of \,$\liek$\, so that \,$(\ad(j)|\liem)^2 = -\id_\liem$\, holds. Because
\,$\liez(\liek)$\, is one-dimensional, this condition already determines \,$j$\, up to sign; we find via computations with the \textsf{Maple} package 
for computation of the Lie bracket of \,$\liee_6$\, that
$$ j = \tfrac23\,(\alpha_1^\sharp-\alpha_6^\sharp) + \tfrac13\,(\alpha_3^\sharp -\alpha_5^\sharp) + K_{2\lambda_1}(1) - K_{2\lambda_2}(1) $$
is one of the two possible choices; here we again denote for \,$\alpha\in\liet^*$\, by \,$\alpha^\sharp \in \liet$\, the dual of \,$\alpha$\, 
with respect to the Killing form \,$\vkap$\,
of \,$\lieg$\,, i.e.~the vector so that \,$\vkap(\alpha^\sharp,\,\cdot\,) = \alpha$\, holds.

Using this presentation of \,$j$\, and the formula \,$Jv = \ad(j)v$\, for \,$v\in\liem$\,, we can again use the \textsf{Maple} package to calculate the
action of \,$J$\, on \,$\liem$\,. In this way, we obtain for \,$c_1,\dotsc,c_4\in\C$\, and \,$t,s \in \R$\,:
\begin{align*}
J(t\,\lambda_1^\sharp + s\,\lambda_2^\sharp) & = \tfrac12\,(M_{2\lambda_1}(t) - M_{2\lambda_2}(s)) \\
J(M_{\lambda_1}(c_1,c_2,c_3,c_4)) & = M_{\lambda_1}(i\,c_1,-i\,c_2,i\,c_3,-i\,c_4) \\
J(M_{\lambda_2}(c_1,c_2,c_3,c_4)) & = M_{\lambda_2}(i\,c_1,-i\,c_2,i\,c_3,-i\,c_4) \\
J(M_{\lambda_3}(c_1,c_2,c_3)) & = M_{\lambda_4}(i\,c_1,-i\,c_2,-i\,\overline{c_3}) \\
J(M_{\lambda_4}(c_1,c_2,c_3)) & = M_{\lambda_3}(i\,c_1,-i\,c_2,i\,\overline{c_3}) \\
J(M_{2\lambda_1}(t)) & = -2t\,\lambda_1^\sharp \\
J(M_{2\lambda_2}(s)) & = 2s\,\lambda_2^\sharp \; . 
\end{align*}
In particular we see that \,$\liem_{\lambda_1}$\, and \,$\liem_{\lambda_2}$\, are complex linear subspaces of \,$\liem$\,, whereas
\,$\liea$\,, \,$\liem_{2\lambda_1} \oplus \liem_{2\lambda_2}$\,, \,$\liem_{\lambda_3}$\, and
\,$\liem_{\lambda_4}$\, are totally real linear subspaces with
\,$J(\liea) = \liem_{2\lambda_1} \oplus \liem_{2\lambda_2}$\, and \,$J(\liem_{\lambda_3}) = \liem_{\lambda_4}$\,.

\subsection[Lie triple systems in \,$E_6/(\Ug(1)\cdot\Spin(10))$\,]{Lie triple systems in \,$\boldsymbol{E_6/(\Ug(1)\cdot\Spin(10))}$\,}
\label{SSe:EIII:lts}

We are now ready to describe the Lie triple systems in \,$\EIII$\,. 

\begin{Def}
Let \,$V$\, be a unitary space. We say that an \,$\R$-linear subspace \,$U\subset V$\, is
\begin{enumerate}
\item of \CP-type \,$(\C,\dim_{\C}(U))$\, if it is a complex subspace of \,$V$\,,
\item of \CP-type \,$(\R,\dim_{\R}(U))$\, if it is a totally real subspace of \,$V$\,.
\end{enumerate}
\end{Def}

\begin{Theorem}
\label{EIII:EIII:cla}
The linear subspaces \,$\liem'$\, of \,$\liem$\, listed in the following are Lie triple systems, and every Lie triple system
\,$\{0\} \neq \liem' \subsetneq \liem$\, is congruent under the isotropy action to one of them.%
\footnote{Please read Remarks~\ref{R:EIII:EIII:extratypes} and \ref{R:EIII:EIII:CN} below before you suspect that there might be Lie triple systems missing from the list.}
\begin{itemize}
\item \,$\boldsymbol{(\mathrm{Geo},\vi=t)}$\, with \,$t \in [0,\tfrac{\pi}{4}]$\, \\
\,$\liem' = \R \,(\cos(t)\lambda_2^\sharp + \sin(t)\lambda_1^\sharp)$\,
(compare Equation~\eqref{eq:EIII:isotropy:liec}).
\item \,$\boldsymbol{(\PP,\vi=0,(\C,5))}$\, \\
\,$\liem' = \R\,\lambda_2^\sharp \oplus \liem_{\lambda_2} \oplus \liem_{2\lambda_2}$\,. 
\item \,$\boldsymbol{(\PP,\vi=\tfrac\pi4,\tau)}$\, with \,$\tau \in \{\Sph^5,\Sph^6,\Sph^7,\Sph^8,\OP^2\}$\,. \\
Put \,$H := \lambda_1^\sharp + \lambda_2^\sharp$\, and \,$\wt{H} := M_{2\lambda_1}(1) + M_{2\lambda_2}(1)$\,. \\
For \,$\tau = \Sph^k$\,: \,$\liem'$\, is a \,$k$-dimensional linear subspace of \,$\R\,H \oplus \liem_{\lambda_4} \oplus \R\,\wt{H}$\,. \\
For \,$\tau = \OP^2$\,: {\tiny \,$\liem' = \R\,H \oplus M_{\lambda_4}(\C,\C,\C) \oplus \Menge{M_{\lambda_1}(c_1,c_2,c_3,c_4) + M_{\lambda_2}(c_2,c_1,-c_4,-c_3)}{c_1,c_2,c_3,c_4\in\C} \oplus \R\,\wt{H}$\,.}
\item \,$\boldsymbol{(\PP\times\PP^1,(\K_1,\ell),\K_2)}$\, with \,$\ell\in\{4,5\}$\, and \,$\K_1,\K_2 \in\{\R,\C\}$\, \\
We have \,$\liem' = \liea \oplus \liem_{\lambda_1}' \oplus \liem_{2\lambda_1}' \oplus \liem_{2\lambda_2}'$\,, where
\,$\liem_{\lambda_1}'$\, is a complex-$(\ell-1)$-dimensional linear subspace of \,$\liem_{\lambda_1}$\,, and where we put for \,$k\in\{1,2\}$\,
\,$\liem_{2\lambda_k}' := \liem_{2\lambda_k}$\, if \,$\K_k = \C$\,, \,$\liem_{2\lambda_k}' := \{0\}$\, if \,$\K_k = \R$\,.
\item \,$\boldsymbol{(Q)}$\, \\
\,$\liem' = \liea \oplus M_{\lambda_3}(\C,\C,\C) \oplus M_{\lambda_4}(\C,\C,\C) \oplus M_{2\lambda_1}(\R) \oplus M_{2\lambda_2}(\R)$\,.
\item \,$\boldsymbol{(Q,\tau)}$\, where \,$\tau$\, is one of the types listed in \cite{Klein:2007-claQ}, Theorem~4.1 for \,$m=8$\,,
i.e.~\,$\tau$\, is one of \,$(\mathrm{G1},k)$\, with \,$k\leq 8$\,, \,$(\mathrm{G2},k_1,k_2)$\, with \,$k_1+k_2\leq 8$\,,
\,$(\mathrm{G3})$\,, \,$(\mathrm{P1},k)$\, with \,$k\leq 8$\,, \,$(\mathrm{P2})$\,, \,$(\mathrm{A})$\,, 
\,$(\mathrm{I1},k)$\, with \,$k\leq 4$\,, and \,$(\mathrm{I2},k)$\, with \,$k\leq 4$\,. \\ 
\,$\liem'$\, is contained in a Lie triple system \,$\wh{\liem}'$\, of type \,$(Q)$\,, corresponding to a complex quadric \,$Q^8$\,,
and regarded as a Lie triple system of \,$\wh{\liem}'$\,, \,$\liem'$\, is of type \,$\tau$\, according to the classification in
\cite{Klein:2007-claQ}, Theorem~4.1.
\item \,$\boldsymbol{(G_2\C^6)}$\, \\
\,$\liem' = \liea \oplus M_{\lambda_1}(\C,\C,0,0) \oplus M_{\lambda_2}(\C,\C,0,0) \oplus M_{\lambda_3}(0,0,\C) \oplus M_{\lambda_4}(0,0,\C) \oplus 
M_{2\lambda_1}(\R) \oplus M_{2\lambda_2}(\R)$\,.
\item \,$\boldsymbol{(G_2\C^6,\tau)}$\,, where \,$\tau$\, is one of the following types listed in \cite{Klein:2007-tgG2}, Theorem~7.1 for \,$n=4$\,:
\,$(\PP,\vi=\arctan(\tfrac12),(\K,k))$\, with \,$\K\in\{\R,\C\}$\, and \,$k \leq 2$\,, 
\,$(\PP,\vi=\tfrac\pi4,(\K,2))$\, with \,$\K\in\{\R,\C,\HH\}$\,, 
\,$(\mathrm{G}_2,(\K,k))$\, with \,$\K\in\{\R,\C\}$\, and \,$k\in\{3,4\}$\,, 
and \,$(\PP\times\PP,(\K,k),(\K',k'))$\, with \,$\K,\K'\in\{\R,\C\}$\, and \,$k+k'\leq 4$\,. \\
\,$\liem'$\, is contained in a Lie triple system \,$\wh{\liem}'$\, of type \,$(G_2\C^6)$\,, corresponding to a complex Grassmannian \,$G_2(\C^6)$\,,
and regarded as a Lie triple system of \,$\wh{\liem}'$\,, \,$\liem'$\, is of type \,$\tau$\, according to the classification in
\cite{Klein:2007-tgG2}, Theorem~7.1.
\item \,$\boldsymbol{(G_2\HH^4)}$\, \\
\,$\liem' = \liea \oplus M_{\lambda_1}(\R,\R,\R,\R) \oplus M_{\lambda_2}(i\R,i\R,i\R,i\R) \oplus M_{\lambda_3}(\R,\R,\R) \oplus M_{\lambda_4}(\R,\R,\R)$\,.
\item \,$\boldsymbol{(G_2\HH^4,\tau)}$\,, where \,$\tau$\, is one of the following types listed in \cite{Klein:2007-tgG2}, Theorem~5.3 for \,$n=2$\,:
\,$(\PP,\vi=0,(\K,2))$\, with \,$\K\in\{\R,\C,\HH\}$\,,
\,$(\Sph,\vi=\arctan(\tfrac13),3)$\,, 
\,$(\PP,\vi=\tfrac\pi4,(\Sph^3))$\,, \,$(\PP,\vi=\tfrac\pi4,(\HH,1))$\,,
\,$(\Sph^5,\vi=\tfrac\pi4)$\,, 
\,$(\mathrm{G}_2,(\HH,1))$\,, 
\,$(\Sph^1\times\Sph^5,k)$\, with \,$3 \leq k\leq 5$\,, 
and \,$(\mathrm{Sp}_2)$\,. \\
\,$\liem'$\, is contained in a Lie triple system \,$\wh{\liem}'$\, of type \,$(G_2\HH^4)$\,, corresponding locally to a quaternionic Grassmannian \,$G_2(\HH^4)$\,,
and regarded as a Lie triple system of \,$\wh{\liem}'$\,, \,$\liem'$\, is of type \,$\tau$\, according to the classification in
\cite{Klein:2007-tgG2}, Theorem~5.3.
\item \,$\boldsymbol{(\mathrm{DIII})}$\, \\
\,$\liem' = \liea \oplus M_{\lambda_1}(\C,0,\C,0) \oplus M_{\lambda_2}(\C,0,\C,0) \oplus M_{\lambda_3}(0,\C,\C) \oplus M_{\lambda_4}(0,\C,\C)
\oplus M_{2\lambda_1}(\R) \oplus M_{2\lambda_2}(\R)$.
\end{itemize}
We call the full name \,$(\mathrm{Geo},\vi=t)$\,, \,$(\PP,\vi=0,(\C,5))$\, etc.~given in the above table the \emph{type} of the 
Lie triple systems which are isotropy-congruent to the space given in that entry. Then every Lie triple system of \,$\liem$\, is of exactly one type.

In the type names of Lie triple systems of rank \,$1$\,, the value given in the form \,$\vi=t$\, is the isotropy angle 
(see the end of Section~\ref{Se:generallts}) of the Lie triple systems of that type.

The Lie triple systems \,$\liem'$\, of the various types have the properties given in the following table. 
The column ``isometry type'' gives the isometry type of the totally geodesic submanifolds corresponding to the Lie triple systems of the respective type
in abbreviated form (without specification of the scaling factors of the Riemannian metrics), for the details see Section~\ref{SSe:EIII:tgsub}.
\begin{center}
\begin{longtable}{|c|c|c|c|c|c|}
\hline
type of \,$\liem'$\, & $\dim(\liem')$ & $\rk(\liem')$ & \begin{minipage}{2cm} \begin{center} {\tiny \,$\liem'$\, complex or \\ totally real? \par} \end{center} \end{minipage}
& \,$\liem'$\, maximal & isometry type \\
\hline
\endhead
\hline
\endfoot
\,$(\mathrm{Geo},\vi=t)$\, & $1$ & $1$ & totally real & no & \,$\R$\, or \,$\Sph^1$\, \\
\,$(\PP,\vi=0,(\C,5))$\, & $10$ & $1$ & complex & no & \,$\CP^5$\, \\
\,$(\PP,\vi=\tfrac\pi4,\Sph^\ell)$\, & $\ell$ & $1$ & totally real & no & \,$\Sph^\ell$\, \\
\,$(\PP,\vi=\tfrac\pi4,\OP^2)$\, & $16$ & $1$ & totally real & yes & \,$\OP^2$\, \\
\hline
\,$(\PP\times\PP^1,(\R,\ell),\R)$\, & $\ell + 1$ & $2$ & totally real & no & \,$\RP^\ell\times\RP^1$\, \\
\,$(\PP\times\PP^1,(\R,\ell),\C)$\, & $\ell + 2$ & $2$ & neither & no & \,$\RP^\ell\times\CP^1$\,  \\
\,$(\PP\times\PP^1,(\C,\ell),\R)$\, & $2\ell + 1$ & $2$ & neither & no & \,$\CP^\ell \times \RP^1$\, \\
\,$(\PP\times\PP^1,(\C,\ell),\C)$\, & $2\ell + 2$ & $2$ & complex & for \,$\ell=5$\, & \,$\CP^\ell \times \CP^1$\, \\
\,$(Q)$\, & $16$ & $2$ & complex & yes & \,$Q^8$\, \\
\,$(Q,\tau)$\, & \multicolumn{3}{c|}{see \cite{Klein:2007-claQ}, Theorem~4.1} & no & \\
\,$(G_2\C^6)$\, & $16$ & $2$ & complex & yes & \,$G_2(\C^6)$\,  \\
\,$(G_2\C^6,\tau)$\, & \multicolumn{3}{c|}{see \cite{Klein:2007-tgG2}, Theorem~7.1} & no & \\
\,$(G_2\HH^4)$\, & $16$ & $2$ & totally real & yes & \,$G_2(\HH^4)/\Z_2$\, \\
\,$(G_2\HH^4,\tau)$\, & \multicolumn{2}{c|}{see \cite{Klein:2007-tgG2}, Theorem~5.3} & totally real & no & \\
\,$(\mathrm{DIII})$\, & $20$ & $2$ & complex & yes & \,$\SO(10)/\Ug(5)$\, \\
\end{longtable}
\end{center}
\end{Theorem}

\begin{Remark}
\label{R:EIII:EIII:phirel}
The Lie triple systems of type \,$(Q,\tau)$\,, \,$(G_2\C^6,\tau)$\, and \,$(G_2\HH^4,\tau)$\, are contained in Lie triple systems of type \,$(Q)$\,
(corresponding to a complex quadric \,$Q^8$\,), \,$(G_2\C^6)$\, (corresponding to \,$G_2(\C^6)$\,) and \,$(G_2\HH^4)$\, (corresponding to \,$G_2(\HH^4)/\Z_2$\,),
respectively. To obtain explicit descriptions of these types, 
one needs to apply the results in \cite{Klein:2007-claQ} and \cite{Klein:2007-tgG2} on the classification of Lie
triple systems in these spaces.

To be able to do so, it is important to know how the root systems of the Lie triple systems of type \,$(Q)$\,, \,$(G_2\C^6)$\, and \,$(G_2\HH^4)$\, are embedded
in the root system of \,$\EIII$\,, and also how the function \,$\vi$\, parametrizing the orbits of the isotropy action defined for \,$Q^m$\, and \,$G_2(\K^n)$\,
in \cite{Klein:2007-claQ} resp.~in \cite{Klein:2007-tgG2} relates to the corresponding function \,$\vi$\, defined for \,$\EIII$\, in the present paper.

Because the Lie triple systems of type \,$(Q)$\,, \,$(G_2\C^6)$\, and \,$(G_2\HH^4)$\, have maximal rank in \,$\EIII$\,, their 
respective root systems \,$\Delta_{(Q)}$\,,
\,$\Delta_{(G_2\C^6)}$\, and \,$\Delta_{(G_2\HH^4)}$\, are simply subsets of
the root system \,$\Delta$\, of \,$\EIII$\, (see Proposition~\ref{P:cla:subroots:subroots-neu}(b), and also see the proof of Theorem~\ref{EIII:EIII:cla} below). 
In fact, from the definition of these types in
Theorem~\ref{EIII:EIII:cla} it follows immediately that we have
\begin{gather*}
\Delta_{(Q)} = \{\pm \lambda_3,\pm\lambda_4,\pm2\lambda_1,\pm2\lambda_2\} \;, \\
\Delta_{(G_2\C^6)} = \Delta \;, \\
\Delta_{(G_2\HH^4)} = \{\pm \lambda_1,\pm\lambda_2,\pm \lambda_3,\pm\lambda_4\} \;. 
\end{gather*}

For each of the types \,$\hat\tau \in \{(Q), (G_2\C^6), (G_2\HH^4)\}$\, we now let \,$\liem_{\hat\tau}$\, be a Lie triple system of \,$\EIII$\, of type \,$\hat\tau$\,,
and let \,$\vi_{\hat\tau}: \liem_{\hat\tau}\setminus\{0\} \to [0,\tfrac\pi4]$\, be the function parametrizing the orbits of the isotropy action of the symmetric space
corresponding to \,$\liem_{\hat\tau}$\, (i.e.~\,$Q^8$\,, \,$G_2(\C^6)$\, or \,$G_2(\HH^4)/\Z_2$\,) as introduced in \cite{Klein:2007-claQ} at the beginning of
Section~4.2 resp.~in \cite{Klein:2007-tgG2}, Section~4. Note that in these cases, we always measured the angle \,$\vi(v)$\, from the vector corresponding to the
shortest root present in \,$Q^n$\, resp.~\,$G_2(\C^n)$\, for large \,$n$\,, even if this root vanishes for certain small values of \,$n$\, (as happens
for \,$G_2(\HH^4)$\,). Keeping this in mind, and considering the root systems \,$\Delta_{\hat\tau}$\, as given above, we see that the functions \,$\vi_{\hat\tau}$\,
is related to the function \,$\vi:\liem\setminus\{0\} \to [0,\tfrac\pi4]$\, parametrizing the isotropy orbits of \,$\EIII$\, by
\begin{gather*}
\vi_{(Q)}(v) = \tfrac\pi4 - \vi(v) \qmq{for \,$v\in\liem_{(Q)}\setminus\{0\}$\,,} \\
\vi_{(G_2\C^6)}(v) = \vi(v) \qmq{for \,$v\in\liem_{(G_2\C^6)}\setminus\{0\}$\,,} \\
\vi_{(G_2\HH^4)}(v) = \tfrac\pi4 - \vi(v) \qmq{for \,$v\in\liem_{(G_2\HH^4)}\setminus\{0\}$\,.} \\
\end{gather*}
\end{Remark}

\begin{Remark}
\label{R:EIII:EIII:extratypes} 
We now introduce alternative definitions for some types of Lie triple systems, to make it more intuitive that indeed all congruence classes of Lie triple systems
are covered in Theorem~\ref{EIII:EIII:cla}, and also to simplify the notations in what follows.

First, we consider the types \,$(G_2\C^6,\tau)$\, resp.~\,$(G_2\HH^4,\tau)$\, also for those types \,$\tau$\, listed in \cite{Klein:2007-tgG2}, Theorem~7.1 for \,$n=4$\,
resp.~in \cite{Klein:2007-tgG2}, Theorem~5.3 for \,$n=2$\, which have not been mentioned in Theorem~\ref{EIII:EIII:cla}. Then a Lie triple system of \,$\EIII$\,
is contained in a Lie triple system of type \,$(G_2\C^6)$\, resp.~\,$(G_2\HH^4)$\, if and only if it is of type \,$(G_2\C^6,\tau)$\, resp.~of type
\,$(G_2\HH^4,\tau)$\, with some \,$\tau$\,. 

Moreover, we define the types \,$(\PP,\vi=0,(\K,\ell))$\, for any \,$\K\in\{\R,\C\}$\, and \,$\ell \leq 5$\,: We say that a linear subspace
of \,$\liem$\, is of that type if and only if it is isotropy-congruent to \,$\liem' = \R\,\lambda_2^\sharp \oplus \liem_{\lambda_2}' \oplus \liem_{2\lambda_2}'$\,,
where \,$\liem_{\lambda_2}'\subset\liem_{\lambda_2}$\, is a linear subspace of \,$\CP$-type \,$(\K,\ell-1)$\, and we put \,$\liem_{2\lambda_2}' := \liem_{2\lambda_2}$\,
if \,$\K=\C$\,, \,$\liem_{2\lambda_2}' := \{0\}$\, if \,$\K=\R$\,. Any such space is a Lie triple system of \,$\liem$\,, and the Lie triple systems of these
types are exactly those which are contained in a Lie triple system of type \,$(\PP,\vi=0,(\C,5))$\,. 

Likewise, we can define the type \,$(\PP,\vi=\tfrac\pi4,\tau)$\, also for \,$\tau = \Sph^\ell$\, with \,$\ell \leq 4$\, and for \,$\tau = \KP^2$\,
with \,$\K\in\{\R,\C,\HH\}$\, in the following way:
We put \,$H := \lambda_1^\sharp + \lambda_2^\sharp$\, and \,$\wt{H} := M_{2\lambda_1}(1) + M_{2\lambda_2}(1)$\,. Then a Lie triple system is of type
\,$(\PP,\vi=\tfrac\pi4,\Sph^\ell)$\, if it is isotropy-congruent to a \,$\ell$-dimensional linear subspace of
\,$\R\,H \oplus \liem_{\lambda_4} \oplus \R\,\wt{H}$\,. A Lie triple system is of type \,$(\PP,\vi=\tfrac\pi4,\KP^2)$\,
if it is congruent to the Lie triple system \,$\liem'$\,, where we have 
\begin{itemize}
\item for \,$\K=\R$\,: \,$\liem' = \R\,H \oplus \Menge{M_{\lambda_1}(t,0,0,0) + M_{\lambda_2}(0,t,0,0)}{t\in\R}$\,. 
\item for \,$\K=\C$\,: \,$\liem' = \R\,H \oplus \Menge{M_{\lambda_1}(c,0,0,0) + M_{\lambda_2}(0,c,0,0)}{c\in\C} \oplus \R\,\wt{H} $\,. 
\item for \,$\K=\HH$\,: {\small \,$\liem' = \R\,H \oplus \Menge{M_{\lambda_1}(c_1,c_2,0,0) + M_{\lambda_2}(c_2,c_1,0,0)}{c_1,c_2\in\C} \oplus M_{\lambda_4}(0,0,\C) \oplus \R\,\wt{H}$\,.}
\end{itemize}
Then the Lie triple systems of \,$\EIII$\, which are contained in a Lie triple system of type \,$(\PP,\vi=\tfrac\pi4,\OP^2)$\, are exactly those which are
of a type of the form \,$(\PP,\vi=\tfrac\pi4,\tau)$\,. 

Finally, the type \,$(\PP\times\PP^1,(\K_1,\ell),\K_2)$\, can be defined also for \,$\ell \leq 3$\, by applying the same definition as in the Theorem.
Then the Lie triple systems
of \,$\EIII$\, which are contained in \,$(\PP\times\PP^1,(\C,5),\C)$\, are exactly those which are of the type
\,$(\PP\times\PP^1,(\K_1,\ell),\K_2)$\, with some \,$\K_1,\K_2\in\{\R,\C\}$\, and \,$\ell \leq 5$\,.

These ``newly defined'' types are identical, however, to types of the form \,$(Q,\tau)$\, or \,$(G_2\C^6,\tau)$\, defined in Theorem~\ref{EIII:EIII:cla}. 
This is detailed in the following table:
\begin{center}
\begin{longtable}{|c|c|}
\hline
The type ... defined here & is identical to the type ... from Theorem~\ref{EIII:EIII:cla}. \\
\hline
\endhead
\hline
\endfoot
$(G_2\C^6,(\PP,\vi=0,(\R,k)))$ & $(Q,(\mathrm{I2},k))$ \\
$(G_2\C^6,(\PP,\vi=0,(\C,k)))$ & $(Q,(\mathrm{I1},k))$ \\
$(G_2\C^6,(\Sph,\vi=\arctan(\tfrac13),2))$ & $(Q,(\mathrm{A}))$ \\
$(G_2\C^6,(\PP,\vi=\tfrac\pi4,(\R,1)))$ & $(Q,(\mathrm{P1},1))$ \\
$(G_2\C^6,(\PP,\vi=\tfrac\pi4,(\C,1)))$ & $(Q,(\mathrm{P1},2))$ \\
$(G_2\C^6,(\PP,\vi=\tfrac\pi4,(\Sph^3)))$ & $(Q,(\mathrm{P1},3))$ \\
$(G_2\C^6,(\PP,\vi=\tfrac\pi4,(\HH,1)))$ & $(Q,(\mathrm{P1},4))$ \\
$(G_2\C^6,(\mathrm{G}_2,(\R,1)))$ & $(Q,(\mathrm{I2},2))$ \\
$(G_2\C^6,(\mathrm{G}_2,(\C,1)))$ & $(Q,(\mathrm{I1},2))$ \\
$(G_2\C^6,(\mathrm{G}_2,(\R,2)))$ & $(Q,(\mathrm{P1},4))$ \\
$(G_2\C^6,(\mathrm{G}_2,(\C,2)))$ & $(Q,(\mathrm{G1},4))$ \\
$(G_2\C^6,(\Sph^1\times\Sph^5,k))$ & $(Q,(\mathrm{P2},1,k))$ \\
$(G_2\C^6,(Q_3))$ & $(Q,(\mathrm{G2},3))$ \\
\hline
$(G_2\HH^4,(\PP,\vi=0,\tau'))$ with \,$\dim(\tau')=1$ & $(Q,(\mathrm{P1},w(\tau')))$ \\
$(G_2\HH^4,(\PP,\vi=\arctan(\tfrac13),2))$ & $(G_2\C^6,(\PP,\vi=\arctan(\tfrac12),(\R,2)))$ \\
$(G_2\HH^4,(\PP,\vi=\tfrac\pi4,(\K,1)))$ with \,$\K\in\{\R,\C\}$\, & $(G_2\C^6,(\PP,\vi=0,(\K,1)))$ \\
$(G_2\HH^4,(\mathrm{G}_2,(\K,1)))$ with \,$\K\in\{\R,\C\}$\, & $(G_2\C^6,(\PP,\vi=\tfrac\pi4,(\K,2)))$ \\
$(G_2\HH^4,(\mathrm{G}_2,(\R,2)))$ & $(G_2\C^6,(\PP\times\PP,(\R,1),(\R,1)))$ \\
$(G_2\HH^4,(\mathrm{G}_2,(\C,2)))$ & $(G_2\C^6,(\mathrm{G}_2,(\R,4)))$ \\
$(G_2\HH^4,(\PP\times\PP,\tau',\tau''))$ with \,$\dim(\tau')=\dim(\tau'')=1$\, & $(Q,(\mathrm{G2},w(\tau'),w(\tau'')))$ \\
$(G_2\HH^4,(\Sph^1 \times \Sph^5,1))$ & $(G_2\C^6,(\PP\times\PP,(\R,1),(\R,1)))$\, \\
$(G_2\HH^4,(\Sph^1 \times \Sph^5,2))$ & $(G_2\C^6,(\PP\times\PP,(\C,1),(\R,1)))$\, \\
$(G_2\HH^4,(Q_3))$ & $(G_2\C^6,(\mathrm{G}_2,(\R,3)))$ \\
\hline
$(\PP,\vi=0,(\R,1))$ & $(\mathrm{Geo},\vi=0)$ \\
$(\PP,\vi=0,(\R,2))$ & $(G_2\C^6,(\PP,\vi=0,(\C,1)))$ \\
$(\PP,\vi=0,(\R,3))$ & $(G_2\C^6,(\PP,\vi=0,(\Sph^3)))$ \\
$(\PP,\vi=0,(\R,4))$ & $(G_2\C^6,(\PP,\vi=0,(\HH,1)))$ \\
$(\PP,\vi=0,(\R,5))$ & $(G_2\HH^4,(\Sph^5,\vi=\tfrac\pi4))$ \\
$(\PP,\vi=0,(\C,\ell))$ with \,$\ell \leq 4$ & $(Q,(\mathrm{I1},\ell))$ \\
\hline
$(\PP,\vi=\tfrac\pi4,\Sph^\ell)$ with \,$\ell\leq 4$ & $(Q,(\mathrm{P1},\ell))$ \\
$(\PP,\vi=\tfrac\pi4,\KP^2)$ with \,$\K\in\{\R,\C,\HH\}$ & $(G_2\C^6,(\PP,\vi=\tfrac\pi4,(\K,2)))$ \\
\hline
$(\PP\times\PP^1,(\K_1,\ell),\K_2)$\, with \,$\ell\leq 3$\, & $(G_2\C^6,(\PP\times\PP,(\K_1,\ell),(\K_2,1)))$ \\
\end{longtable}
\end{center}
As an example for proving these identities, we consider the type \,$(\PP,\vi=0,(\C,4))$\,.
To prove that this type is identical to the type \,$(Q,(\mathrm{I1},4))$\,, it suffices to show that the space \,$\liem' := \R\,\lambda_2^\sharp \oplus
M_{\lambda_2}(\C,\C,\C,0) \oplus M_{2\lambda_2}(\R)$\, of type \,$(\PP,\vi=0,(\C,4))$\, is isotropy-congruent to a space \,$\Ad(g)\liem'$\,
contained in the Lie triple system
\,$\wh{\liem}' := \liea \oplus M_{\lambda_3}(\C,\C,\C) \oplus M_{\lambda_4}(\C,\C,\C) \oplus M_{2\lambda_1}(\R) \oplus M_{2\lambda_2}(\R)$\, of type \,$(Q)$\,. 
Because \,$\Ad(g)\liem'$\, has the isotropy angle \,$0$\, with respect to \,$\EIII$\,, it has the isotropy angle \,$\tfrac\pi4$\, with respect to \,$Q^8$\,
(see Remark~\ref{R:EIII:EIII:phirel}); because it is also a complex subspace, it then must be of type \,$(Q,(\mathrm{I1},4))$\, by the
classification of Lie triple systems of the complex quadric given in \cite{Klein:2007-claQ}, Theorem~4.1.  --- To show that such an isotropy-congruence
indeed holds, notice that with \,$Z := K_{\alpha_7}(\sqrt{8}) \in \liek$\, we have \,$\ad(Z)\lambda_2^\sharp = \ad(Z)M_{2\lambda_2}(t) = 0$\, and
\,$\ad(Z)M_{\lambda_2}(c_1,c_2,c_3,0) = M_{\lambda_3}(c_2,c_1,-\overline{c_3}) + M_{\lambda_4}(-c_1,-c_2,c_3)$\, for any \,$t\in\R$\,, \,$c_1,c_2,c_3\in\C$\,. 
This shows that with \,$g := \exp(\tfrac\pi2\,Z)\in K$\,
we have \,$\Ad(g)\liem' = \R\lambda_2^\sharp \oplus \Menge{M_{\lambda_3}(c_2,c_1,-\overline{c_3}) + M_{\lambda_4}(-c_1,-c_2,c_3)}{c_1,c_2,c_3\in\C} \oplus
M_{2\lambda_2}(\R) \subset \wh{\liem}'$\,. 
\end{Remark}

\begin{Remark}
\label{R:EIII:EIII:CN} 
For the space \,$\EIII$\,, Table~VIII of \cite{Chen/Nagano:totges2-1978} 
correctly lists the \emph{local} isometry types of the \emph{maximal}
totally geodesic submanifolds. However, the totally geodesic submanifolds corresponding to the types
\,$(G_2\HH^4)$\, and \,$(Q)$\, are of isometry type \,$G_2(\HH^4)/\Z_2$\, resp.~\,$G_2^+(\R^{10}) \cong Q^8$\,
(see Section~\ref{SSe:EIII:tgsub}),
and not of isometry type \,$G_2(\HH^4)$\, resp.~\,$G_2(\R^{10})$\, (as \cite{Chen/Nagano:totges2-1978} claims).

It should be noted that \,$\EIII$\, contains spaces of rank \,$1$\, as totally geodesic submanifolds
in a ``skew'' position in the sense that their geodesic diameter
is strictly larger than the geodesic diameter of the ambient space \,$\EIII$\,.  However, none of them is maximal in \,$\EIII$\,. 
The ``skew'' totally geodesic submanifolds which are maximal among the totally geodesic submanifolds of \,$\EIII$\, of rank \,$1$\, are
those of the types \,$(G_2\HH^4,(\PP,\vi=\arctan(\tfrac13),3))$\, (isometric to an \,$\RP^3$\, of sectional curvature \,$\tfrac25$\,),
\,$(Q,(\mathrm{A}))$\, (isometric to a 2-sphere of radius \,$\tfrac12\sqrt{10}$\,) and 
\,$(G_2\C^6,(\PP,\vi=\arctan(\tfrac12),(\C,2))$\, (isometric to a \,$\CP^2$\, of holomorphic sectional curvature \,$\tfrac45$\,). The existence of
these ``skew'' totally geodesic submanifolds cannot be inferred from Table~VIII of \cite{Chen/Nagano:totges2-1978}. For explicit constructions
of these ``skew'' totally geodesic submanifolds in \,$G_2(\HH^4)$\,, \,$Q^3$\, resp.~\,$G_2(\C^6)$\,, see \cite{Klein:2007-tgG2}, Sections~6 and 7.%
\footnote{The most general of these constructions in \cite{Klein:2007-tgG2} is the construction of a ``skew'' \,$\HP^2$\, (of type \,$(\PP,\vi=\arctan(\tfrac12),(\HH,2))$\,)
in \,$G_2(\HH^7)$\, described in Section~6 of \cite{Klein:2007-tgG2}. It is based on the fundamental \,$14$-dimensional representation with quaternionic structure
of \,$\Sp(3)$\,, which is realized as a sub-representation of the representation of \,$\Sp(3)$\, on \,$\bigwedge^3 \C^6$\,, see also
\cite{Broecker-tom-Dieck:1985}, p.~269ff. I would like to remark that this representation is not equivalent to the representation of \,$\Sp(3)$\,
on \,$\frakJ(3,\HH)^{\C}$\, involved in Cartan's construction of isoparametric hypersurfaces in the sphere. This is easily seen, because the latter representation,
although it is also 14-dimensional and irreducible, admits a real structure, and thus cannot admit a quaternionic structure, see
\cite{Broecker-tom-Dieck:1985}, Proposition II.6.5, p.~98.} 
\end{Remark}

The remainder of the present section is concerned with the proof of Theorem~\ref{EIII:EIII:cla}.

We first mention that it is easily checked using the \textsf{Maple} implementation of the algorithms for the computation
of the curvature tensor that the spaces defined in the theorem, and therefore
also the linear subspaces \,$\liem' \subset \liem$\, which are congruent to one of them, are Lie triple systems. It is also
easily seen that the information in the table concerning the dimension, the rank, and the question if \,$\liem'$\, is complex
or totally real is correct (for the latter, use the description of the complex structure of \,$\EIII$\, given in 
Section~\ref{SSe:EIII:geometry}). The information on the isometry type of the corresponding totally geodesic submanifolds will
be proved in Section~\ref{SSe:EIII:tgsub}.

We next show that the information on the maximality of the Lie triple systems given in the table is correct.
For this purpose, we presume that the list of Lie triple systems given in the theorem is in fact complete; this will
be proved in the remainder of the present section.

That the Lie triple systems which are claimed to be maximal in the table indeed are: This is clear for the type
\,$(\mathrm{DIII})$\,, because it has the maximal dimension among all the Lie triple systems of \,$\EIII$\,. 
The Lie triple systems of the types \,$(\PP,\vi=\tfrac\pi4,\OP^2)$\,, \,$(Q)$\,, \,$(G_2\C^6)$\, and \,$(G_2\HH^4)$\,
all are of dimension \,$16$\,, therefore if they were not maximal, they could only be contained in a Lie triple
system of type \,$(\mathrm{DIII})$\,, because these are the only ones of greater dimension. The spaces of the
types \,$(\PP,\vi=\tfrac\pi4,\OP^2)$\, and \,$(G_2\HH^4)$\, are real forms of \,$\EIII$\,, and therefore
cannot be contained in a (complex) Lie triple systems of type \,$(\textrm{DIII})$\,. The restricted Dynkin diagrams
with multiplicities of the Lie triple systems of type \,$(Q)$\, and \,$(\textrm{DIII})$\, are 
$\xymatrix@=.4cm{ \mathop{\bullet}^1 \ar@{=>}[r] & \mathop{\bullet}^6 }$ and
$\xymatrix@=.4cm{ \mathop{\bullet}^4 \ar@{<=>}[r] & \mathop{\bullet\hspace{-.29cm}\bigcirc}^{4[1]} }$, respectively.
Thus the short roots in a Lie triple system of type \,$(Q)$\, have greater multiplicity than all the roots in a Lie triple
system of type \,$(\mathrm{DIII})$\,, and hence a Lie triple system of type \,$(Q)$\, cannot be contained in any
Lie triple system of type \,$(\mathrm{DIII})$\, either. Assume that the Lie triple system 
\,$\liem' := \liea \oplus M_{\lambda_1}(\C,\C,0,0) \oplus M_{\lambda_2}(\C,\C,0,0) \oplus M_{\lambda_3}(0,0,\C) \oplus M_{\lambda_4}(0,0,\C) \oplus M_{2\lambda_1}(\R) \oplus M_{2\lambda_2}(\R)$\, of type \,$(G_2\C^6)$\, were contained in a Lie triple system of type
\,$(\textrm{DIII})$\, i.e.~in a space isotropy-congruent to 
\,$\wh{\liem}' := \liea \oplus M_{\lambda_1}(\C,0,\C,0) \oplus M_{\lambda_2}(\C,0,\C,0) \oplus M_{\lambda_3}(0,\C,\C) \oplus M_{\lambda_4}(0,\C,\C) \oplus M_{2\lambda_1}(\R) \oplus M_{2\lambda_2}(\R)$\,.
Then there would exist \,$g\in K$\, so that \,$\Ad(g)$\, maps \,$M_{\lambda_k}(\C,\C,0,0)$\, onto \,$M_{\lambda_k}(\C,0,\C,0)$\,
for \,$k\in\{1,2\}$\,. But this is a contradiction to the fact that the action of \,$\Ad(g)$\, commutes with the map
\,$M_{\lambda_1}(c_1,c_2,c_3,c_4)\mapsto M_{\lambda_2}(c_2,c_1,c_4,c_3)$\, (Proposition~\ref{P:EIII:isotropy}), so also the Lie triple
systems of type \,$(G_2\C^6)$\, cannot be contained in a Lie triple system of type \,$(\mathrm{DIII})$\,. 
Finally, we note that the Lie triple systems
of type \,$(\PP\times\PP^1,(\C,5),\C)$\, are of rank \,$2$\, and have the Dynkin diagram 
$\xymatrix@=.4cm{ \mathop{\bullet\hspace{-.29cm}\bigcirc}^{8[1]} & \strut\!\!\!\!\!\!\!\mathop{\bullet}^1 }$.
They have a restricted root of multiplicity \,$8$\,, which is greater than the multiplicity of any root in any
other Lie triple system of \,$\EIII$\, of rank \,$2$\,. Therefore also this type is maximal.

That no Lie triple systems are maximal besides those mentioned in the theorem follows from the following table:

\begin{center}
\begin{tabular}{|c|c|}
\hline
Every Lie triple system of type ... & is contained in a Lie triple system of type ... \\
\hline
$(\mathrm{Geo},\vi=t)$ & $(\PP\times\PP^1,(\R,1),\R)$ \\
$(\PP,\vi=\tfrac\pi4,\Sph^k)$ & $(\PP,\vi=\tfrac\pi4,\OP^2)$ \\
$(\PP\times\PP^1,(\K,\ell),\R)$ & $(\PP\times\PP^1,(\K,\ell),\C)$ \\
$(\PP\times\PP^1,(\K,\ell),\C)$ with \,$(\K,\ell) \neq (\C,5)$\, & $(\PP\times\PP^1,(\C,5),\C)$ \\
$(Q,\tau)$ & $(Q)$ \\
$(G_2\C^6,\tau)$ & $(G_2\C^6)$ \\
$(G_2\HH^4,\tau)$ & $(G_2\HH^4)$ \\
\hline
\end{tabular}
\end{center}

We now turn to the proof that the list of Lie triple systems of \,$\EIII$\, given in Theorem~\ref{EIII:EIII:cla} is indeed
complete. For this purpose, we let an arbitrary Lie triple system \,$\liem'$\, of \,$\liem$\,, \,$\{0\} \neq \liem' \subsetneq \liem$\,,
be given. In the sequel, we will also use the additional names for types of Lie triple systems introduced in Remark~\ref{R:EIII:EIII:extratypes};
it has been shown there that these types are equivalent to other types defined in Theorem~\ref{EIII:EIII:cla}.

Because the symmetric space \,$\EIII$\, is of rank \,$2$\,, the rank of \,$\liem'$\, is either \,$1$\, or \,$2$\,. We will handle
these two cases separately in the sequel.

We first suppose that \,$\liem'$\, is a Lie triple system of rank \,$2$\,. Let us fix a Cartan subalgebra \,$\liea$\, of \,$\liem'$\,;
because of \,$\rk(\liem') = \rk(\liem)$\,, \,$\liea$\, is then also a Cartan subalgebra of \,$\liem$\,. In relation to this situation, we use
the notations introduced in Sections~\ref{Se:generallts} and \ref{SSe:EIII:geometry}. In particular,
we consider the positive root system
\,$\Delta_+ := \{\lambda_1,\lambda_2,\lambda_3,\lambda_4,2\lambda_1,2\lambda_2\}$\, of the root system \,$\Delta := \Delta(\liem,\liea)$\, of \,$\liem$\,,
and also the root system \,$\Delta' := \Delta(\liem',\liea)$\, of \,$\liem'$\,. By Proposition~\ref{P:cla:subroots:subroots-neu}(b), \,$\Delta'$\,
is a root subsystem of \,$\Delta$\,, and therefore \,$\Delta_+' := \Delta' \cap \Delta_+$\, is a positive system of roots for \,$\Delta'$\,.
Moreover, in the root space decompositions of \,$\liem$\, and \,$\liem'$\,
\begin{equation}
\label{eq:EIII:rk2:decomp}
\liem = \liea \;\oplus\; \bigoplus_{\lambda\in\Delta_+} \liem_\lambda \qmq{and}
\liem' = \liea \;\oplus\; \bigoplus_{\lambda\in\Delta_+'} \liem_\lambda'
\end{equation}
the root space \,$\liem_\lambda'$\, of \,$\liem'$\, with respect to \,$\lambda \in \Delta_+'$\, is related to the corresponding root space \,$\liem_\lambda$\, of \,$\liem$\,
by \,$\liem_\lambda' = \liem_\lambda \cap \liem'$\,. 

As was noted in Section~\ref{SSe:EIII:geometry}, \,$\liem_{\lambda_k}$\, is a complex subspace of \,$\liem$\, for \,$k\in\{1,2\}$\,. 
The following proposition describes how the position of \,$\liem_{\lambda_k}'$\, in \,$\liem_{\lambda_k}$\, with respect to the complex structure
is controlled by the presence of the root \,$2\lambda_k$\, in \,$\Delta'$\,.

\begin{Prop}
\label{P:EIII:rk2:ctr}
For \,$k\in\{1,2\}$\,, \,$\liem_{\lambda_k}'$\, is either a complex or a totally real subspace of \,$\liem_{\lambda_k}$\,; it is a complex subspace
if and only if \,$2\lambda_k \in \Delta'$\, holds.
\end{Prop}

\beweis
First suppose \,$2\lambda_k \in \Delta'$\,. Because of \,$n_{2\lambda_k}=1$\, we then have \,$\liem_{2\lambda_k}' = \liem_{2\lambda_k} = M_{2\lambda_k}(\R)$\,.
For any given \,$v\in\liem_{\lambda_1}'$\, we have%
\footnote{The evaluation of \,$R$\, is done here, as in all the following situations, using the \textsf{Maple} package described in 
\cite{Klein:2007-Satake}, as explained in the Introduction.}
\,$R(\lambda_k^\sharp,v)M_{2\lambda_k}(1) = -\tfrac18\,Jv$\,, and this vector is a member of \,$\liem'$\, by the fact that \,$\liem'$\, is a Lie triple system.
Thus \,$Jv \in \liem_{\lambda_k} \cap \liem'
= \liem_{\lambda_k}'$\, holds, and hence \,$\liem_{\lambda_k}'$\, is a complex subspace of \,$\liem_{\lambda_k}$\,. 

Now suppose \,$2\lambda_k \not\in \Delta'$\,. For any given \,$v,w \in \liem_{\lambda_k}'$\, we have \,$\liem' \ni R(\lambda_k^\sharp,v)w = \tfrac18\,\g{v}{w}\,\lambda_1^\sharp
+ \tfrac18\,\g{v}{Jw}\,M_{2\lambda_k}(1)$\,; because of \,$2\lambda_k\not\in\Delta'$\, it follows that \,$\g{v}{Jw}=0$\, holds. Hence \,$\liem_{\lambda_k}'$\,
is a totally real subspace of \,$\liem_{\lambda_k}$\,.
\beweisende

We now distinguish three cases depending on the structure of \,$\Delta'$\,, which we will treat separately in the sequel:
\begin{enumerate}
\item \,$\lambda_3,\lambda_4 \in \Delta'$\,
\item either, but not both, of \,$\lambda_3$\, and \,$\lambda_4$\, are members of \,$\Delta'$\,
\item \,$\lambda_3,\lambda_4 \not\in\Delta'$\,
\end{enumerate}

\textbf{Case (a).} Because of \,$\lambda_3 \in \Delta'$\,, there exists \,$v\in\liem_{\lambda_3}'$\, with \,$\|v\|=1$\,. By Proposition~\ref{P:EIII:isotropy}, 
there exists \,$g\in K_0 \subset K$\, so that \,$\Ad(g)$\, maps \,$v$\, into \,$M_{\lambda_3}(0,0,1)$\,, and therefore \,$\liem'$\, into another Lie triple
system \,$\liem'' := \Ad(g)\liem'$\,, so that we have \,$M_{\lambda_3}(0,0,1) \in \liem_{\lambda_3}''$\,. 
This argument shows that we may suppose without loss of generality that \,$M_{\lambda_3}(0,0,1) \in \liem_{\lambda_3}'$\, holds.

We have for any \,$v = M_{\lambda_1}(c_1,c_2,c_3,c_4)$\,: \,$R(\lambda_1^\sharp,v)M_{\lambda_3}(0,0,1) = \tfrac{\sqrt{2}}{16}\,M_{\lambda_2}(c_1i, c_2i,
-c_3i,-c_4i)$\,, and for any \,$v = M_{\lambda_2}(c_1,c_2,c_3,c_4)$\,: \,$R(\lambda_2^\sharp,v)M_{\lambda_3}(0,0,1) = -\tfrac{\sqrt{2}}{16}\,
M_{\lambda_1}(c_1i, c_2i, -c_3i,-c_4i)$\,. Because of the fact that \,$\liem'$\, is a Lie triple system, it follows that we have for any \,$c_1,\dotsc,c_4 \in \C$\,
\begin{equation}
\label{eq:EIII:rk2:m1m2}
M_{\lambda_1}(c_1,c_2,c_3,c_4) \in \liem_{\lambda_1}' \;\Longleftrightarrow\; M_{\lambda_2}(c_1i,c_2i,-c_3i,-c_4i) \in \liem_{\lambda_2}' \; .
\end{equation}
This equivalence in particular implies \,$\bigr(\lambda_1\in\Delta' \Leftrightarrow \lambda_2\in\Delta'\bigr)$\, and \,$n_{\lambda_1}' = n_{\lambda_2}'$\,.

Because \,$\Delta'$\, is invariant under the Weyl transformation given by the reflection in \,$(\lambda_3^\sharp)^{\perp,\liea}$\, (note \,$\lambda_3 \in \Delta'$\,),
we also have \,$\bigr(2\lambda_1\in\Delta' \Leftrightarrow 2\lambda_2\in\Delta'\bigr)$\,.

Let us first suppose \,$2\lambda_1,2\lambda_2 \in \Delta'$\,. Then
\,$\liem_{\lambda_k}'$\, is a complex subspace of \,$\liem_{\lambda_k}$\, for \,$k\in\{1,2\}$\, by Proposition~\ref{P:EIII:rk2:ctr}. Hence \,$n := n_{\lambda_1}' = n_{\lambda_2}'$\,
is an even number, and we consider the possible values \,$0,2,4,6,8$\, for \,$n$\, individually in the sequel.

If \,$n=0$\, holds, we have \,$\Delta' = \{\pm \lambda_3, \pm \lambda_4, \pm 2\lambda_1, \pm 2\lambda_2\}$\,; this is a closed root subsystem of \,$\Delta$\,. 
Therefore the maximal linear subspace \,$\wh{\liem}' := \liea \oplus \bigoplus_{\lambda \in \Delta'} \liem_\lambda$\, of \,$\liem$\, corresponding to \,$\Delta'$\, is a 
Lie triple system (see Proposition~\ref{P:cla:assoclts}); 
its corresponding Dynkin diagram with multiplicities is $\xymatrix@=.4cm{ \mathop{\bullet}^1 \ar@{=>}[r] & \mathop{\bullet}^6 }$.
Therefore the totally geodesic submanifold corresponding to \,$\wh{\liem}'$\, is locally isometric to the complex quadric \,$Q^8$\,. 
\,$\liem'$\, is also regarded as a subspace of \,$\wh{\liem}'$\, a Lie triple system; therefore \,$\liem'$\, is of one of the types 
described in the classification of the Lie triple systems of \,$Q^m$\, in \cite{Klein:2007-claQ}.  
It follows that if \,$\liem' = \wh{\liem}'$\, holds, then
\,$\liem'$\, is of type \,$(Q)$\,; otherwise it is of type \,$(Q,\tau)$\,, where \,$\tau$\, is one of the types of Lie triple systems
of \,$\wh{\liem}'$\, as described in \cite{Klein:2007-claQ}, Theorem~4.1 for \,$m=8$\,. 

For \,$n \neq 0$\,, an argument based on Proposition~\ref{P:EIII:isotropy} similar to the previous one shows that we may suppose without loss
of generality besides the earlier condition \,$M_{\lambda_3}(0,0,1) \in \liem_{\lambda_3}'$\, also \,$M_{\lambda_1}(1,0,0,0) \in \liem_{\lambda_1}'$\,.
Because \,$\liem_{\lambda_1}'$\, is a complex subspace of \,$\liem_{\lambda_1}$\,, we then in fact have
\,$M_{\lambda_1}(\C,0,0,0) \subset \liem_{\lambda_1}'$\,. This fact induces further relations between the root spaces of \,$\liem'$\,
besides \eqref{eq:EIII:rk2:m1m2}, which we now explore.

For any \,$v := M_{\lambda_3}(d_1,d_2,d_3)$\, we have
\,$u:= R(\lambda_2^\sharp, v)M_{\lambda_1}(1,0,0,0)=-\tfrac{\sqrt{2}}{16}\,M_{\lambda_2}(id_3,0,i\,\overline{d_2},-id_1)$\, and 
\,$R(u,M_{\lambda_1}(1,0,0,0))\lambda_1^\sharp = \tfrac{1}{128} v + \tfrac{1}{128} M_{\lambda_4}(-d_1,d_2,\overline{d_3})$\,. 
An analogous calculation applies starting with \,$v = M_{\lambda_4}(d_1,d_2,d_3)$\,,
and in this way we see via the fact that \,$\liem'$\, is a Lie triple system:
\begin{gather}
\label{eq:EIII:rk2:m3m2}
M_{\lambda_3}(d_1,d_2,d_3) \in \liem_{\lambda_3}' \;\Longrightarrow\; M_{\lambda_2}(id_3,0,i\,\overline{d_2},-id_1) \in \liem_{\lambda_2}' \;, \\
\label{eq:EIII:rk2:m3m4}
M_{\lambda_3}(d_1,d_2,d_3) \in \liem_{\lambda_3}' \;\Longleftrightarrow\; M_{\lambda_4}(-d_1,d_2,\overline{d_3}) \in \liem_{\lambda_4}' \; .
\end{gather}
Moreover for any \,$v := M_{\lambda_2}(c_1,c_2,c_3,c_4)$\, we have \,$R(M_{\lambda_1}(1,0,0,0),v)\lambda_3^\sharp = \tfrac{\sqrt{2}}{8}
M_{\lambda_3}(-c_4i,-\overline{c_3}\,i,c_1i)$\, and therefore, again by the fact that \,$\liem'$\, is a Lie triple system
\begin{equation}
\label{eq:EIII:rk2:m2m3}
M_{\lambda_2}(c_1,c_2,c_3,c_4) \in \liem_{\lambda_2}' \;\Longrightarrow\; M_{\lambda_3}(-c_4i,-\overline{c_3}\,i,c_1i) \in \liem_{\lambda_3}' \; . 
\end{equation}

We can use these relations to draw the following consequences from the fact \,$M_{\lambda_1}(\C,0,0,0) \subset \liem_{\lambda_1}'$\,: First, from 
\eqref{eq:EIII:rk2:m1m2} we obtain \,$M_{\lambda_2}(\C,0,0,0) \subset \liem_{\lambda_2}'$\,. By \eqref{eq:EIII:rk2:m2m3} therefrom
\,$M_{\lambda_3}(0,0,\C) \subset \liem_{\lambda_3}'$\, follows, and therefrom we finally obtain by \eqref{eq:EIII:rk2:m3m4}:
\,$M_{\lambda_4}(0,0,\C) \subset \liem_{\lambda_4}'$\,. Remember for the sequel also that we have
\,$\liem_{2\lambda_k} = M_{2\lambda_k}(\R)$\, for \,$k \in \{1,2\}$\,.

If \,$n=2$\, holds, then we in fact have \,$\liem_{\lambda_1}' = M_{\lambda_1}(\C,0,0,0)$\, and \,$\liem_{\lambda_2}' = M_{\lambda_2}(\C,0,0,0)$\,;
because of \eqref{eq:EIII:rk2:m3m2} we then have \,$\liem_{\lambda_3}' = M_{\lambda_3}(0,0,\C)$\,, and therefore \,$\liem_{\lambda_4}' = M_{\lambda_4}(0,0,\C)$\,
by \eqref{eq:EIII:rk2:m3m4}.
Thus we see by the root space decomposition~\eqref{eq:EIII:rk2:decomp} that 
$$ \liem' = \liea \oplus M_{\lambda_1}(\C,0,0,0) \oplus M_{\lambda_2}(\C,0,0,0) \oplus M_{\lambda_3}(0,0,\C) \oplus M_{\lambda_4}(0,0,\C) \oplus M_{2\lambda_1}(\R)
\oplus M_{2\lambda_2}(\R) $$
holds, and thus \,$\liem'$\, is of type \,$(\mathrm{G}_2\C^6, (\mathrm{G}_2,(\C,3)))$\,. 

If \,$n=4$\, holds, then the Dynkin diagram with multiplicities corresponding to \,$\liem'$\, is 
$\xymatrix@=.4cm{ \mathop{\bullet}^\ell \ar@{<=>}[r] & \mathop{\bullet\hspace{-.29cm}\bigcirc}^{4[1]} }$ with some \,$1\leq \ell\leq 6$\,;
from the classification of irreducible Riemannian symmetric spaces (see, for example, \cite{Loos:1969-2}, p.~119, 146)
we see that \,$\ell = 2$\, and \,$\ell=4$\, are the only possibilities. If \,$\ell=2$\, holds, we have \,$\liem_{\lambda_3}' = M_{\lambda_3}(0,0,\C)$\,
and \,$\liem_{\lambda_4}' = M_{\lambda_4}(0,0,\C)$\,. Because of \,$n=4$\, we see from
\eqref{eq:EIII:rk2:m2m3} that \,$\liem_{\lambda_2}' = M_{\lambda_2}(\C,\C,0,0)$\, and therefore by \eqref{eq:EIII:rk2:m1m2} also
\,$\liem_{\lambda_1}' = M_{\lambda_1}(\C,\C,0,0)$\, has to hold. 
Thus we see that
$$ \liem' = \liea \oplus M_{\lambda_1}(\C,\C,0,0) \oplus M_{\lambda_2}(\C,\C,0,0) \oplus M_{\lambda_3}(0,0,\C) \oplus M_{\lambda_4}(0,0,\C) \oplus 
M_{2\lambda_1}(\R) \oplus M_{2\lambda_2}(\R) $$
holds, and hence \,$\liem'$\, is of type \,$(\mathrm{G}_2\C^6)$\,. 
On the other hand, if \,$\ell=4$\, holds, we let \,$v \in \liem_{\lambda_1}'$\, be a unit vector which is orthogonal to \,$M_{\lambda_1}(\C,0,0,0) \subset \liem_{\lambda_1}'$\,;
we have \,$v=M_{\lambda_1}(0,c_2,c_3,c_4)$\, with some \,$c_2,c_3,c_4\in\C$\,, and 
\begin{equation}
\label{eq:EIII:rk2:DIII5}
\liem_{\lambda_1}' = M_{\lambda_1}(\C,0,0,0)\oplus\R v\oplus \R Jv
\end{equation}
holds. Because of \,$v\in\liem_{\lambda_1}'$\,, we have \,$M_{\lambda_2}(0,c_2 i,-c_3 i,-c_4 i) \in \liem_{\lambda_2}'$\, by \eqref{eq:EIII:rk2:m1m2}, 
therefore \,$M_{\lambda_3}(-c_4, \overline{c_3}, 0)\in \liem_{\lambda_3}'$\, by \eqref{eq:EIII:rk2:m2m3}, thus
\,$M_{\lambda_2}(0,0,-i c_3,-i c_4) \in \liem_{\lambda_2}'$\, by \eqref{eq:EIII:rk2:m3m2}, and hence finally 
\,$M_{\lambda_1}(0,0,c_3,c_4) \in \liem_{\lambda_1}'$\, by \eqref{eq:EIII:rk2:m1m2}.
From \eqref{eq:EIII:rk2:DIII5} and the explicit description of \,$J$\, in Section~\ref{SSe:EIII:geometry} we see that 
$$ (0,0,c_3,c_4) \in \R(0,c_2,c_3,c_4) \oplus \R(0,-ic_2,ic_3,-ic_4) \; $$
holds; this implies that we have either \,$c_2=0$\, or \,$c_3=c_4=0$\,. In fact, \,$c_3=c_4=0$\, is impossible, because then we would have
\,$\liem_{\lambda_1}' = M_{\lambda_1}(\C,\C,0,0)$\,, therefore by \eqref{eq:EIII:rk2:m1m2} also \,$\liem_{\lambda_2}' = M_{\lambda_2}(\C,\C,0,0)$\,,
and therefore by \eqref{eq:EIII:rk2:m3m2} \,$\liem_{\lambda_3}' \subset M_{\lambda_3}(0,0,\C)$\,, in contradiction to \,$\ell=4$\,. Therefore we have \,$c_2=0$\,
and thus \,$v = M_{\lambda_1}(0,0,c_3,c_4)$\,. We have \,$\ad(K_{\alpha_4}(2))M_{\lambda_1}(0,0,c_3,c_4) = M_{\lambda_1}(0,0,-\overline{c_4},\overline{c_3})$\,, therefore by 
the application of a rotation \,$\Ad(\exp(H))$\, with suitable \,$H \in \R\,\alpha_4^\sharp \oplus K_{\alpha_4}(\C) \cong \liesu(2)$\, to \,$\liem'$\,, we can arrange
\,$c_4 = 0$\,, and thus \,$\liem_{\lambda_1}' = M_{\lambda_1}(\C,0,\C,0)$\,. Then we have \,$\liem_{\lambda_2}' = M_{\lambda_2}(\C,0,\C,0)$\, by 
\eqref{eq:EIII:rk2:m1m2}, \,$\liem_{\lambda_3}' = M_{\lambda_3}(0,\C,\C)$\, by \eqref{eq:EIII:rk2:m2m3} and the fact that \,$\ell=4$\,, and
\,$\liem_{\lambda_4}' = M_{\lambda_4}(0,\C,\C)$\, by \eqref{eq:EIII:rk2:m3m4}. Therefore
$$ \liem' = \liea \oplus M_{\lambda_1}(\C,0,\C,0) \oplus M_{\lambda_2}(\C,0,\C,0) \oplus M_{\lambda_3}(0,\C,\C) \oplus M_{\lambda_4}(0,\C,\C) \oplus M_{2\lambda_1}(\R)
\oplus M_{2\lambda_2}(\R) $$
is of type \,$(\mathrm{DIII})$\,. 

The case \,$n=6$\, cannot occur, because the Dynkin diagram with multiplicities corresponding to \,$\liem'$\, would then be
$\xymatrix@=.4cm{ \mathop{\bullet}^\ell \ar@{<=>}[r] & \mathop{\bullet\hspace{-.29cm}\bigcirc}^{6[1]} }$ with some \,$1\leq \ell\leq 6$\,;
\eqref{eq:EIII:rk2:m2m3} shows \,$\ell \geq 4$\,. But the classification of irreducible Riemannian symmetric spaces
(see \cite{Loos:1969-2}, p.~119, 146) shows that no symmetric space with such a diagram exists.

Finally, if \,$n=8$\, holds, we have \,$\liem_{\lambda_k}' = \liem_{\lambda_k}$\, for \,$k\in\{1,2\}$\,, from \eqref{eq:EIII:rk2:m2m3} we obtain
\,$\liem_{\lambda_3}' = \liem_{\lambda_3}$\,, from \eqref{eq:EIII:rk2:m3m4} we then obtain \,$\liem_{\lambda_4}' = \liem_{\lambda_4}$\,, and we also have
\,$\liem_{2\lambda_k}' = \liem_{2\lambda_k}$\, for \,$k\in\{1,2\}$\,. Thus we have \,$\liem' = \liem$\,. 

\medskip

Let us now consider the case where \,$2\lambda_1,2\lambda_2 \not\in \Delta'$\,. Then \,$\liem_{\lambda_1}'$\, and \,$\liem_{\lambda_2}'$\, are totally real subspaces
of \,$\liem_{\lambda_1}$\, resp.~\,$\liem_{\lambda_2}$\, by Proposition~\ref{P:EIII:rk2:ctr}. We either have \,$\lambda_1,\lambda_2\in\Delta'$\, or
\,$\lambda_1,\lambda_2 \not\in\Delta'$\, because of the invariance of \,$\Delta'$\, under the Weyl transformation induced by \,$\lambda_3 \in \Delta'$\,.
If \,$\lambda_1,\lambda_2\not\in\Delta'$\,, i.e.~\,$\Delta' = \{\pm \lambda_3,\pm \lambda_4\}$\, holds, then \,$\liem'$\, is again contained in a Lie triple
system \,$\wh{\liem}'$\, of type \,$(Q)$\,, and therefore, by the classification of Lie triple systems in \,$\wh{\liem}'$\, given in 
\cite{Klein:2007-claQ}, \,$\liem'$\, is of type \,$(Q,\tau)$\,, where \,$\tau$\, is one of the types listed in Theorem~4.1 of \cite{Klein:2007-claQ} for
\,$m=8$\,. 

So we now suppose \,$\lambda_1,\lambda_2 \in \Delta'$\,. 
Once again using Proposition~\ref{P:EIII:isotropy}, we may suppose without loss of generality that \,$\liem_{\lambda_1}' \subset M_{\lambda_1}(\R,\R,\R,\R)$\,
holds because \,$\liem_{\lambda_1}'$\, is totally real, and also that \,$M_{\lambda_1}(1,0,0,0) \in \liem_{\lambda_1}'$\, holds. 
The proof of Equations~\eqref{eq:EIII:rk2:m3m2}--\eqref{eq:EIII:rk2:m2m3} was based only on the fact that \,$M_{\lambda_1}(1,0,0,0) \in \liem_{\lambda_1}'$\, holds, 
and therefore these equations will again be valid in the present situation.
Therefore we have \,$\liem_{\lambda_2}' \subset M_{\lambda_2}(i\R,i\R,i\R,i\R)$\, by \eqref{eq:EIII:rk2:m1m2}, 
then \,$\liem_{\lambda_3}' \subset M_{\lambda_3}(\R,\R,\R)$\, by \eqref{eq:EIII:rk2:m2m3},
and then \,$\liem_{\lambda_4}' \subset M_{\lambda_4}(\R,\R,\R)$\, by \eqref{eq:EIII:rk2:m3m4}. 

Therefore \,$\liem'$\, is contained in the Lie triple system
$$ \wh{\liem}' := \liea \oplus M_{\lambda_1}(\R,\R,\R,\R) \oplus M_{\lambda_2}(i\R,i\R,i\R,i\R) \oplus M_{\lambda_3}(\R,\R,\R) \oplus M_{\lambda_4}(\R,\R,\R) \;, $$
of type \,$(\mathrm{G}_2\HH^4)$\,. \,$\wh{\liem}'$\, has the Dynkin diagram 
$\xymatrix@=.4cm{ \mathop{\bullet}^3 \ar@{=>}[r] & \mathop{\bullet}^4 }$, and therefore the corresponding totally geodesic submanifold is locally isometric
to \,$G_2(\HH^4)$\,. \,$\liem'$\, is also a Lie triple system of \,$\wh{\liem}'$\,, therefore we have either \,$\liem' = \wh{\liem}'$\,, or \,$\liem'$\,
is of one of the types described in the classification of Lie triple systems of \,$G_2(\HH^{n+2})$\, given in Theorem~5.3 of \cite{Klein:2007-tgG2}
for \,$n=2$\,. It follows that \,$\liem'$\, is either of type \,$(G_2\HH^4)$\, (if \,$\liem' = \wh{\liem}'$\, holds), or of type \,$(G_2\HH^4,\tau)$\,,
where \,$\tau$\, is one of the types of Lie triple systems of \,$\wh{\liem}'$\, as described in Theorem~5.3 of \cite{Klein:2007-tgG2} for \,$n=2$\,.

\textbf{Case (b).} Here we suppose that either, but not both, of \,$\lambda_3$\, and \,$\lambda_4$\, are in \,$\Delta'$\,. Without loss of generality
we may suppose \,$\lambda_3\in\Delta'$\,, \,$\lambda_4\not\in\Delta'$\,. Because \,$\Delta'$\, is invariant under its Weyl transformation group, we then
have \,$\Delta' = \{\pm \lambda_3\}$\, and therefore \,$\liem' = \liea \oplus \liem_{\lambda_3}'$\, is of type \,$(Q,(\mathrm{G2},\ell,1))$\, with
\,$\ell := 1+n_{\lambda_3}'$\,, \,$2\leq\ell\leq 7$\,. 

\textbf{Case (c).} 
So now we have \,$\lambda_3,\lambda_4\not\in\Delta'$\,. Let us first consider the case where at least one of the roots \,$\lambda_1$\, and \,$\lambda_2$\,
is not in \,$\Delta'$\,. Without loss of generality we suppose \,$\lambda_2\not\in\Delta'$\,, so that we have \,$\Delta' \subset \{\pm \lambda_1,\pm 2\lambda_1,
\pm 2\lambda_2\}$\,. For \,$k\in\{1,2\}$\, we put \,$\K_k := \C$\, if \,$2\lambda_k\in\Delta'$\,, \,$\K_k := \R$\, if \,$2\lambda_k\not\in\Delta'$\,.
Proposition~\ref{P:EIII:rk2:ctr} then shows that \,$\liem_{\lambda_1}'$\, is a linear subspace of \,$\liem_{\lambda_1}$\, of type \,$(\K_1,\dim_{\K_1}(\liem_{\lambda_1}'))$\,.
It follows that \,$\liem'$\, is of type \,$(\PP\times\PP^1,(\K_1,1+\dim_{\K_1}(\liem_{\lambda_1}')),\K_2)$\,.

Now consider the case \,$\lambda_1,\lambda_2\in\Delta'$\,. As before, we may use Proposition~\ref{P:EIII:isotropy} to suppose without 
loss of generality that \,$M_{\lambda_1}(1,0,0,0) \in \liem_{\lambda_1}'$\, holds.
Let \,$v\in\liem_{\lambda_2}'$\, be given, say \,$v=M_{\lambda_2}(c_1,c_2,c_3,c_4)$\, with \,$c_1,\dotsc,c_4\in\C$\,, then we have 
\,$\liem' \ni R(M_{\lambda_1}(1,0,0,0),v)\lambda_3^\sharp = -\tfrac{\sqrt{2}}{8}\,M_{\lambda_3}(ic_4,i\,\overline{c_3}, -ic_1)$\,. Because of \,$\lambda_3\not\in\Delta'$\,
it follows that we have \,$c_1=c_3=c_4=0$\, and thus we have \,$\liem_{\lambda_2}' \subset M_{\lambda_2}(0,\C,0,0)$\,. Without loss of generality
we may suppose \,$M_{\lambda_2}(0,1,0,0) \in \liem_{\lambda_2}'$\,. Now let \,$v\in \liem_{\lambda_1}'$\, be given, say 
\,$v=M_{\lambda_1}(c_1,c_2,c_3,c_4)$\, with \,$c_1,\dotsc,c_4\in\C$\,, then we have \,$\liem' \ni R(M_{\lambda_2}(0,1,0,0),v)\lambda_3^\sharp 
= \tfrac{\sqrt{2}}{8}\,M_{\lambda_3}(ic_3,i\,\overline{c_4},ic_2)$\,. Because of \,$\lambda_3\not\in\Delta'$\, we obtain \,$c_2=c_3=c_4=0$\,, and thus
\,$\liem_{\lambda_1}' \subset M_{\lambda_1}(\C,0,0,0)$\,. Because \,$\liem_{\lambda_k}$\, is either complex or totally real according to whether
\,$2\lambda_k$\, is or is not a member of \,$\Delta'$\, by Proposition~\ref{P:EIII:rk2:ctr}, we see that \,$\liem'$\, is of the type 
\,$(\mathrm{G}_2\C^6,(\PP\times\PP, (\K_1,2), (\K_2,2)))$\,, where
for \,$k\in\{1,2\}$\, we put \,$\K_k := \C$\, if \,$2\lambda_k\in\Delta'$\,, \,$\K_k := \R$\, if \,$2\lambda_k\not\in\Delta'$\,. 

This completes the classification of the Lie triple systems of \,$\EIII$\, of rank \,$2$\,. 

\bigskip

We now turn our attention to the case where \,$\liem'$\, is a Lie triple system of rank \,$1$\,. Via the application of the isotropy action of \,$\EIII$\,, we may
suppose without loss of generality 
that \,$\liem'$\, contains a unit vector \,$H$\, from the closure \,$\overline{\liec}$\, 
of the positive Weyl chamber \,$\liec$\, of \,$\liem$\, (with respect to \,$\liea$\, and our
choice of positive roots). 
By Equations~\eqref{eq:EIII:isotropy:liec} and \eqref{eq:EIII:isotropy:vivt} we then have with \,$\vi_0 := \vi(H) \in [0,\tfrac{\pi}{4}]$\, 
\begin{equation}
\label{eq:EIII:rk1:H}
H = \cos(\vi_0)\,\lambda_2^\sharp + \sin(\vi_0)\,\lambda_1^\sharp \; .
\end{equation}

Because of \,$\rk(\liem')=1$\,, \,$\liea' := \R\,H$\, is a Cartan subalgebra of \,$\liem'$\,, and we have \,$\liea' = \liea \cap \liem'$\,. 
It follows from Proposition~\ref{P:cla:subroots:subroots-neu}(a) that the root systems \,$\Delta'$\,
and \,$\Delta$\, of \,$\liem'$\, resp.~\,$\liem$\, with respect to \,$\liea'$\, resp.~to \,$\liea$\, are related by
\begin{equation}
\label{eq:EIII:rk1:Delta'Delta}
\Delta' \;\subset\; \Mengegr{\lambda(H)\,\alpha_0}{\lambda\in\Delta,\,\lambda(H)\neq 0}
\end{equation}
with the linear form \,$\alpha_0: \liea' \to \R,\; tH\mapsto t$\,; moreover for \,$\liem'$\, we have the root space decomposition
\begin{equation}
\label{eq:EIII:rk1:m'decomp}
\liem' = \liea' \oplus \bigoplus_{\alpha\in\Delta_+'} \liem_\alpha' 
\end{equation}
where for any root \,$\alpha\in\Delta'$\,, the corresponding root space \,$\liem_\alpha'$\, is given by
\begin{equation}
\label{eq:EIII:rk1:malpha'}
\liem_\alpha' = \left( \bigoplus_{\substack{\lambda \in \Delta \\ \lambda(H) = \alpha(H)}} \liem_\lambda \right) \;\cap\; \liem' \; . 
\end{equation}

If \,$\Delta' = \varnothing$\, holds, then we have \,$\liem' = \R H$\,, and therefore \,$\liem'$\, is then of type \,$(\mathrm{Geo},\vi=\vi_0)$\,. Otherwise
by the same consideration as in my classification of the Lie triple systems in \,$G_2(\HH^n)$\, (\cite{Klein:2007-tgG2}, the beginning of Section~5.2),
we see that 
$$ \vi_0 \;\in\; \{0,\arctan(\tfrac13),\arctan(\tfrac12),\tfrac\pi4\} $$
holds; moreover in the cases \,$\vi_0 = \arctan(\tfrac13)$\, and \,$\vi_0 = \arctan(\tfrac12)$\,, \,$\Delta'$\, cannot have elementary roots in the sense
of Definition~\ref{D:cla:subroots:Elemcomp}. 

In the sequel we consider the four possible values for \,$\vi_0$\, individually.

\textbf{The case \,$\boldsymbol{\vi_0=0}$\,.} In this case we have \,$H = \lambda_2^\sharp$\, by Equation~\eqref{eq:EIII:rk1:H} and therefore
$$ \lambda_1(H) = 2\lambda_1(H) = 0,\; \lambda_2(H) = \lambda_3(H) = \lambda_4(H) = 1,\; 2\lambda_2(H) = 2 \; . $$
Thus we have
\,$\Delta' \subset \{\pm \alpha,\pm 2\alpha\}$\, with \,$\alpha := \lambda_2|\liea' = \lambda_3|\liea' = \lambda_4|\liea'$\, by Equation~\eqref{eq:EIII:rk1:Delta'Delta},
\,$\liem' = \R\,H \oplus \liem_{\alpha}' \oplus \liem_{2\alpha}'$\, by Equation~\eqref{eq:EIII:rk1:m'decomp}, and
\,$\liem_{\alpha}' \subset \liem_{\lambda_2} \oplus \liem_{\lambda_3} \oplus \liem_{\lambda_4}$\, and
\,$\liem_{2\alpha}' \subset \liem_{2\lambda_2}$\, by Equation~\eqref{eq:EIII:rk1:malpha'}.

We first note that if in fact \,$\liem_{\alpha}' \subset \liem_{\lambda_2}$\, holds (this is in particular the case for \,$\alpha\not\in\Delta'$\,),
then by the same argument as in the proof of Proposition~\ref{P:EIII:rk2:ctr}, \,$\liem_\alpha'$\, is either a complex or a totally real linear subspace 
of \,$\liem_{\lambda_2}$\,, depending on whether \,$2\alpha$\, is or is not
a member of \,$\Delta'$\,. Therefore \,$\liem'$\, then is of type 
\,$(\PP,\vi=0,(\K,\ell))$\, with \,$\K\in\{\R,\C\}$\, and \,$\ell := \dim_{\K}(\liem_\alpha') + 1$\,. 

Also, if \,$\liem_{\alpha}' \subset \liem_{\lambda_3} \oplus \liem_{\lambda_4}$\, holds,
then \,$\liem'$\, is contained in a Lie triple system \,$\wh{\liem}'$\, of type \,$(Q)$\,, and therefore \,$\liem'$\, is of type 
\,$(Q,\tau)$\,, where \,$\tau$\, is a type given in \cite{Klein:2007-claQ}, Theorem~4.1 for \,$m=4$\,.

Thus we now suppose \,$\liem_\alpha' \not\subset \liem_{\lambda_2}$\, and \,$\liem_{\alpha}' \not\subset \liem_{\lambda_3} \oplus \liem_{\lambda_4}$\,, 
in particular we have \,$\alpha\in\Delta'$\,. We will show that in this situation, \,$\liem'$\, is conjugate under the isotropy action to another
Lie triple system whose corresponding root space decomposition satisfies \,$\liem_\alpha' \subset \liem_{\lambda_2}'$\, or \,$\liem_\alpha' \subset
\liem_{\lambda_3}' \oplus \liem_{\lambda_4}'$\,. It then follows by the above discussion that \,$\liem'$\, is of one of the types of the Theorem.

It follows from our hypotheses \,$\liem_\alpha' \not\subset \liem_{\lambda_2}$\, and \,$\liem_{\alpha}' \not\subset \liem_{\lambda_3} \oplus \liem_{\lambda_4}$\,
that there exists a unit vector \,$v_0 \in \liem'$\,, say \,$v_0 = M_{\lambda_2}(c_1,c_2,c_3,c_4) + M_{\lambda_3}(d_1,d_2,d_3)
+ M_{\lambda_4}(e_1,e_2,e_3)$\,, with \,$(c_1,c_2,c_3,c_4) \neq (0,0,0,0)$\, and \,$(d_1,d_2,d_3,e_1,e_2,e_3) \neq (0,\dotsc,0)$\,. 
By virtue of Proposition~\ref{P:EIII:isotropy} we may suppose without loss of generality
that \,$(c_1,c_2,c_3,c_4) = (t,0,0,0)$\, holds with some \,$t\in \R\setminus \{0\}$\,, and furthermore that
\,$(d_1,d_2,d_3) = (s,0,0)$\, holds with \,$s\in \R$\,. Then we have
\begin{equation}
\label{eq:EIII:rk1:vi0:v0}
v_0 = M_{\lambda_2}(t,0,0,0) + M_{\lambda_3}(s,0,0) + M_{\lambda_4}(e_1,e_2,e_3) \; . 
\end{equation}

Because \,$R(H,v_0)v_0$\, is a member of \,$\liem'$\,, the \,$\liem_{\lambda_1}$-component of this vector, which equals
$$ t\cdot \tfrac{\sqrt{2}}{8}\,M_{\lambda_1}\bigr(\;i\,e_3 \;,\; 0 \;,\; -i\,\overline{e_2} \;,\; i(e_1-s) \; \bigr) \;, $$
must vanish (because of \,$\lambda_1(H)=0$\,), and therefore we have \,$e_1=s$\,, \,$e_2 = e_3 = 0$\,, and therefore
$$ v_0 = M_{\lambda_2}(t,0,0,0) + M_{\lambda_3}(s,0,0) + M_{\lambda_4}(s,0,0) \; . $$
We have
\begin{align*}
\ad(K_{\alpha_{11}}(\sqrt{8}))M_{\lambda_2}(1,0,0,0) & = M_{\lambda_3}(1,0,0) + M_{\lambda_4}(1,0,0) \\
\mathrm{and}\quad
\ad(K_{\alpha_{11}}(\sqrt{8}))(M_{\lambda_3}(1,0,0) + M_{\lambda_4}(1,0,0)) & = -M_{\lambda_2}(1,0,0,0) \;, 
\end{align*}
therefore the 1-parameter subgroup \,$\{\exp(K_{\alpha_{11}}(t))\}_{t\in\R}$\, of the isotropy group acts as a rotation group on the plane
\,$\R\,M_{\lambda_2}(1,0,0,0) \oplus \R(M_{\lambda_3}(1,0,0) + M_{\lambda_4}(1,0,0))$\,; it follows that a suitable member of this
1-parameter group maps (via the isotropy action) \,$v_0$\, onto \,$M_{\lambda_2}(1,0,0,0)$\,. By replacing \,$\liem'$\, with its image
under the action of that element, we may therefore suppose that \,$M_{\lambda_2}(1,0,0,0) \in \liem'$\, holds.

If this replacement causes either \,$\liem_{\alpha}' \subset \liem_{\lambda_2}$\, or \,$\liem_{\alpha}' \subset \liem_{\lambda_3} \oplus \liem_{\lambda_4}$\,
to hold, then we are done. Otherwise, there exists another vector 
\,$v_1 \in \liem_\alpha$\,, say \,$v_1 = M_{\lambda_2}(c_1,c_2,c_3,c_4) + M_{\lambda_3}(d_1,d_2,d_3)
+ M_{\lambda_4}(e_1,e_2,e_3)$\,, with \,$(c_1,c_2,c_3,c_4) \neq (0,0,0,0)$\, and \,$(d_1,d_2,d_3,e_1,e_2,e_3) \neq (0,\dotsc,0)$\,,
and which is orthogonal to \,$M_{\lambda_2}(1,0,0,0) \in \liem_{\alpha}'$\,, whence we have \,$\RE(c_1) = 0$\,. By Proposition~\ref{P:EIII:isotropy} we may
suppose without loss of generality \,$(d_1,d_2,d_3)=(1,0,0)$\, (whilst maintaining the condition \,$M_{\lambda_2}(1,0,0,0) \in \liem_\alpha'$\,).
Then we have 
$$ R(H,M_{\lambda_2}(1,0,0,0))v_1 = \tfrac{\sqrt{2}}{16}\,M_{\lambda_1}(i\,e_3, 0, -i\,\overline{e_2}, i\,(e_1-1)) \;-\; \tfrac18\,M_{2\lambda_2}(\IM(c_1)) \; . $$
Because this vector is a member of \,$\liem'$\,, its \,$\liem_{\lambda_1}$-component must vanish. Thus we have \,$e_1=1$\, and \,$e_2 = e_3 = 0$\,.
Moreover: If \,$2\alpha\not\in\Delta$\,, also the \,$\liem_{2\lambda_2}$-component vanishes, and thus we have 
\,$\IM(c_1)=0$\,, hence \,$c_1=0$\,. On the other hand, if \,$2\alpha\in\Delta$\,, we also have \,$[K_{2\lambda_2}(1),v_0] = -\tfrac12\,M_{\lambda_2}(i,0,0,0) \in \liem'$\,,
and therefore we can replace \,$v_1$\, by \,$v_1-\IM(c_1)\,M_{\lambda_2}(i,0,0,0)$\,. Hence we can suppose \,$c_1=0$\, in any case. Thus \,$v_1$\, is of the form
$$ v_1 = M_{\lambda_2}(0,c_2,c_3,c_4) + M_{\lambda_3}(1,0,0) + M_{\lambda_4}(1,0,0) \; . $$
We now calculate \,$R(H,v_1)v_1 \in \liem'$\,:
$$ R(H,v_1)v_1 = (\tfrac{\|c\|^2}4 + \tfrac{1}2)\cdot H - \tfrac{\sqrt{2}}{4}\,M_{\lambda_1}(i\,\overline{c_4}, 0, i\,\overline{c_2}, 0) \; . $$ 
The \,$\liem_{\lambda_1}$-component of this vector again vanishes, and thus we obtain \,$c_2=c_4=0$\,. Thus we have
$$ v_1 = M_{\lambda_2}(0,0,c_3,0) + M_{\lambda_3}(1,0,0) + M_{\lambda_4}(1,0,0) \; . $$
We now consider the Lie subalgebra \,$\lieb := \R\,\alpha_6^\sharp \oplus K_{\alpha_6}(\C)$\, of \,$\liek$\,, which is isomorphic to \,$\liesu(2)$\,. 
For \,$c\in\C$\, we have
\begin{align*}
\ad(K_{\alpha_6}(2))M_{\lambda_2}(0,0,c,0) & = \tfrac{1}{\sqrt{2}}(M_{\lambda_3}(c,0,0) + M_{\lambda_4}(c,0,0)) \\
\qmq{and} \ad(K_{\alpha_6}(2))\tfrac{1}{\sqrt{2}}(M_{\lambda_3}(c,0,0) + M_{\lambda_4}(c,0,0)) & = -M_{\lambda_2}(0,0,c,0) \; ,
\end{align*}
therefore the connected Lie subgroup \,$B$\, of \,$K$\, with Lie algebra \,$\lieb$\, acts on the complex 2-plane 
\,$M_{\lambda_2}(0,0,\C,0) \oplus \Menge{M_{\lambda_3}(c,0,0)+M_{\lambda_4}(c,0,0)}{c\in\C}$\, as \,$\SU(2)$\,, and further
$$ \ad(K_{\alpha_6}(2))M_{\lambda_2}(1,0,0,0) = 0 \;, $$
therefore the action of \,$B$\, leaves \,$M_{\lambda_2}(1,0,0,0)$\, invariant. Hence, by replacing
\,$\liem'$\, with \,$\Ad(g)\liem'$\, for an appropriate \,$g\in B$\,, we can transform \,$v_1$\, into \,$M_{\lambda_2}(0,0,1,0)$\,, 
while leaving \,$M_{\lambda_2}(1,0,0,0)$\, invariant. 
By replacing \,$\liem'$\, with \,$\Ad(g)\liem'$\,, we can thus ensure besides \,$M_{\lambda_2}(1,0,0,0) \in \liem_\alpha'$\, also \,$M_{\lambda_2}(0,0,1,0) \in 
\liem_\alpha'$\,. 

If this replacement causes \,$\liem_{\alpha}' \subset \liem_{\lambda_2}$\, or \,$\liem_{\alpha}' \subset \liem_{\lambda_3}\oplus\liem_{\lambda_4}$\,
to hold, then we are done. Otherwise, there exists yet another vector 
\,$v_2 \in \liem_\alpha$\,, say \,$v_2 = M_{\lambda_2}(c_1,c_2,c_3,c_4) + M_{\lambda_3}(d_1,d_2,d_3)
+ M_{\lambda_4}(e_1,e_2,e_3)$\, with \,$(d_1,d_2,d_3,e_1,e_2,e_3) \neq (0,\dotsc,0)$\,, 
which is orthogonal to \,$M_{\lambda_2}(1,0,0,0), M_{\lambda_2}(0,0,1,0) \in \liem_{\alpha}'$\,, 
whence we have \,$\RE(c_1) = \RE(c_3) = 0$\,. By an analogous argument as previously, we in fact obtain \,$c_1=c_3=0$\,. 
Then we calculate 
$$ R(H,v_2)M_{\lambda_2}(1,0,0,0) = \tfrac{\sqrt{2}}{16}\,M_{\lambda_1}(\,(e_3-\overline{d_3})i\,,\, 0\,,\, -(\overline{d_2}+\overline{e_2})i\,,\, (e_1-d_1)i\,)
+ \tfrac18 M_{2\lambda_2}(\IM(c_1)) $$
and 
$$ R(H,v_2)M_{\lambda_2}(0,0,1,0) = \tfrac{\sqrt{2}}{16}\,M_{\lambda_1}(\,(e_2-d_2)i\,,\, (e_1-d_1)i\,,\, (d_3+\overline{e_3})i\,,\, 0 \,)
+ \tfrac18 M_{2\lambda_2}(\IM(c_3)) $$
Because these vectors are elements of \,$\liem'$\,, their \,$\liem_{\lambda_1}$-components vanish. From this fact, we derive the equations
\,$e_1=d_1$\, and \,$d_2=d_3=e_2=e_3=0$\,. Using the fact that these equations hold, we now calculate
$$ R(H,v_2)v_2 = (\tfrac{\|c\|^2}4 + \tfrac{|d_1|^2}2)H + \tfrac{\sqrt{2}}{4}\,M_{\lambda_1}(\,i\,\overline{c_4}\,d_1\,,\, 0\,,\, i\,\overline{c_2}\,d_1\,,\, 0\,) \; . $$
Also this vector is an element of \,$\liem'$\,, and thus its \,$\liem_{\lambda_1}$-component once again vanishes, whence it follows (because \,$d_1\neq 0$\,)
that we have \,$c_2=c_4=0$\,, hence \,$v_2 = M_{\lambda_3}(d_1,0,0) + M_{\lambda_4}(d_1,0,0)$\,. 

Now let \,$Z := K_{\alpha_1}(2)$\,. Then we have \,$\ad(Z)v_2 = 0$\, as well as
\begin{align*}
\ad(Z)M_{\lambda_2}(1,0,0,0) & = \tfrac{1}{\sqrt{2}}(M_{\lambda_3}(0,0,1)+M_{\lambda_4}(0,0,1)) =: v_0'\;,\\
\ad(Z)v_0' & = -M_{\lambda_2}(1,0,0,0) 
\end{align*}
and
\begin{align*}
\ad(Z)M_{\lambda_2}(0,0,1,0) & = \tfrac{1}{\sqrt{2}}(M_{\lambda_3}(0,1,0)+M_{\lambda_4}(0,1,0)) =: v_1'\;, \\
\ad(Z)v_1' & = -M_{\lambda_2}(0,0,1,0) \; .
\end{align*}
These equations show that the adjoint action of the one-parameter subgroup \,$B$\, of \,$K$\, tangential to \,$Z$\, leaves the element \,$v_2$\, of \,$\liem'$\, invariant,
whereas it acts as a rotation on the 2-planes spanned by \,$M_{\lambda_2}(1,0,0,0)$\, and \,$M_{\lambda_3}(0,0,1)+M_{\lambda_4}(0,0,1)$\,, resp.~by
\,$M_{\lambda_2}(0,0,1,0)$\, and \,$M_{\lambda_3}(0,1,0)+M_{\lambda_4}(0,1,0)$\,. It follows that there exists \,$g\in B$\, so that we have
\,$\Ad(g)M_{\lambda_2}(1,0,0,0) = v_0'$\,, \,$\Ad(g)M_{\lambda_2}(0,0,1,0) = v_1'$\, and
\,$\Ad(g)v_2 = v_2$\, holds. We replace \,$\liem'$\, by the Lie triple system \,$\Ad(g)\liem'$\,. Then we have \,$v_0',v_1',v_2 \in \liem_\alpha'$\,,
and it turns out that any vector \,$v\in\liem_\alpha'$\, which is orthogonal to the \,$\C$-span of \,$v_0',v_1',v_2$\, is necessarily zero.
Therefore we now have \,$\liem_\alpha' \subset \liem_{\lambda_3}\oplus\liem_{\lambda_4}$\,. 

This completes the treatment of the case \,$\vi_0=0$\,. 

\textbf{The case \,$\boldsymbol{\vi_0=\arctan(\tfrac13)}$\,.} We have
by Equation~\eqref{eq:EIII:rk1:H}: \,$H = \tfrac{3}{\sqrt{10}}\,\lambda_2^\sharp + \tfrac{1}{\sqrt{10}}\,\lambda_1^\sharp$\, and therefore
$$ \lambda_1(H) = \tfrac{1}{\sqrt{10}},\; \lambda_2(H) = \tfrac{3}{\sqrt{10}},\; \lambda_3(H) = \tfrac{2}{\sqrt{10}},\; \lambda_4(H) = \tfrac{4}{\sqrt{10}} \; ,
2\lambda_1(H) = \tfrac{2}{\sqrt{10}},\; 2\lambda_2(H) = \tfrac{6}{\sqrt{10}}\;. $$
Because there are no elementary roots (Definition~\ref{D:cla:subroots:Elemcomp}) in the present case, it follows 
by Equation~\eqref{eq:EIII:rk1:Delta'Delta} that we have \,$\Delta' \subset \{\pm \alpha\}$\, with \,$\alpha := \lambda_3|\liea' = (2\lambda_1)|\liea'$\,,
and by Equations~\eqref{eq:EIII:rk1:m'decomp},\eqref{eq:EIII:rk1:malpha'} we have \,$\liem' = \R H \oplus \liem_\alpha'$\, with \,$\liem_\alpha' \subset \liem_{\lambda_3}
\oplus \liem_{2\lambda_1}$\,.

It follows that \,$\liem'$\, is contained in the Lie triple system \,$\wh{\liem}' := \liea \oplus \liem_{\lambda_3} \oplus \liem_{\lambda_4} \oplus \liem_{2\lambda_1}
\oplus \liem_{2\lambda_2}$\, of type \,$(Q)$\,. \,$\wh{\liem}'$\, corresponds to a complex quadric of complex dimension \,$8$\,, and therefore the Lie triple
systems contained in \,$\wh{\liem}'$\, have been classified in \cite{Klein:2007-claQ}. \,$\liem'$\, is a Lie triple system of rank \,$1$\,, and its isotropy
angle \,$\arctan(\tfrac13)$\, corresponds to the isotropy angle \,$\tfrac\pi4-\arctan(\tfrac13)=\arctan(\tfrac12)$\, in \,$\wh{\liem}'$\,, as has been explained
in Remark~\ref{R:EIII:EIII:phirel}.  It therefore follows from the classification in \cite{Klein:2007-claQ},
Theorem~4.1 that \,$\liem'$\, is, as Lie triple system of \,$\wh{\liem}'$\,, of type \,$(\mathrm{A})$\,. Thus \,$\liem'$\, is as Lie triple system of \,$\liem$\,
of type \,$(Q,(\mathrm{A}))$\,.

\textbf{The case \,$\boldsymbol{\vi_0=\arctan(\tfrac12)}$\,.}
In this case we have
by Equation~\eqref{eq:EIII:rk1:H}: \,$H = \tfrac{2}{\sqrt{5}}\,\lambda_2^\sharp + \tfrac{1}{\sqrt{5}}\,\lambda_1^\sharp$\, and therefore
$$ \lambda_1(H) = \tfrac{1}{\sqrt{5}},\; \lambda_2(H) = \tfrac{2}{\sqrt{5}},\; \lambda_3(H) = \tfrac{1}{\sqrt{5}},\; \lambda_4(H) = \tfrac{3}{\sqrt{5}},\;
2\lambda_1(H) = \tfrac{2}{\sqrt{5}},\; 2\lambda_2(H) = \tfrac{4}{\sqrt{5}}\;. $$
Because there are no elementary roots (Definition~\ref{D:cla:subroots:Elemcomp}) in the present case, it follows 
by Equation~\eqref{eq:EIII:rk1:Delta'Delta} that we have \,$\Delta' \subset \{\pm \alpha,\pm 2\alpha\}$\, with \,$\alpha := \lambda_1|\liea' = \lambda_3|\liea'$\,,
\,$2\alpha = \lambda_2|\liea' = (2\lambda_1)|\liea'$\,,
and by Equations~\eqref{eq:EIII:rk1:m'decomp},\eqref{eq:EIII:rk1:malpha'} we have
\begin{equation}
\label{eq:EIII:rk1:arctan12:m'decomp}
\liem' = \R H \oplus \liem_\alpha' \oplus \liem_{2\alpha}'
\end{equation}
with \,$\liem_\alpha' \subset \liem_{\lambda_1} \oplus \liem_{\lambda_3}$\, and \,$\liem_{2\alpha}' \subset \liem_{\lambda_2} \oplus \liem_{2\lambda_1}$\,. 

We have \,$\alpha^\sharp = \tfrac35\,\lambda_1^\sharp + \tfrac25\,\lambda_3^\sharp$\, and \,$(2\alpha)^\sharp = \tfrac15\,(2\lambda_1)^\sharp + \tfrac45\,\lambda_2^\sharp$\,.
By Proposition~\ref{P:cla:skew} it follows that there exist linear subspaces \,$\liem_{\lambda_3}' \subset \liem_{\lambda_3}$\,, \,$\liem_{2\lambda_1}' \subset
\liem_{2\lambda_1}$\, and isometric linear maps \,$\Phi_\alpha: \liem_{\lambda_3}' \to \liem_{\lambda_1}$\,, \,$\Phi_{2\alpha}: \liem_{2\lambda_1}' \to \liem_{\lambda_2}$\, so that
\begin{equation}
\label{eq:EIII:rk1:arctan12:malpha'skew}
\liem_{\alpha}' = \Mengegr{x + \sqrt{\tfrac32}\,\Phi_\alpha(x)}{x \in \liem_{\lambda_3}'}
\qmq{and}
\liem_{2\alpha}' = \Mengegr{x + 2\,\Phi_{2\alpha}(x)}{x \in \liem_{2\lambda_1}'}
\end{equation}
holds; in particular we have for the multiplicities of the roots of \,$\liem'$\,: \,$n_{\alpha}' \leq 6$\, and \,$n_{2\alpha}' \leq 1$\,. We now consider the cases
\,$2\alpha\in\Delta'$\, and \,$2\alpha\not\in\Delta'$\, separately.

First suppose \,$2\alpha\in\Delta'$\,. Then we have \,$n_{2\alpha}' = 1$\, and \,$\liem_{2\lambda_1}' = \liem_{2\lambda_1} = M_{2\lambda_1}(\R)$\,. 
\,$\Phi_{2\alpha}(M_{2\lambda_1}(1))$\, is a unit vector in \,$\liem_{\lambda_2}$\,, and via Proposition~\ref{P:EIII:isotropy}, we may suppose without
loss of generality that \,$\Phi_{2\alpha}(M_{2\lambda_1}(1)) = M_{\lambda_2}(1,0,0,0)$\, holds. By Equation~\eqref{eq:EIII:rk1:arctan12:malpha'skew} 
\,$\liem_{2\alpha}'$\, is then spanned by the vector
\begin{equation}
\label{eq:EIII:rk1:arctan12:v2a'}
v_{2\alpha} := M_{2\lambda_1}(1) + 2\,M_{\lambda_2}(1,0,0,0)\;.
\end{equation}
We now let \,$v \in \liem_\alpha'$\, be given, say \,$v = M_{\lambda_1}(a_1,a_2,a_3,a_4) + M_{\lambda_3}(b_1,b_2,b_3)$\, with \,$a_k,b_\ell \in \C$\,. Because \,$\liem'$\,
is a Lie triple system, we have \,$v_R := R(H,v_{2\alpha})v \in \liem'$\,. The root space decomposition \eqref{eq:EIII:rk1:arctan12:m'decomp} together with
Equations~\eqref{eq:EIII:rk1:arctan12:malpha'skew} shows that therefore the \,$\liem_{\lambda_4}$-component of \,$v_R$\,, which is equal to 
$$ \tfrac{\sqrt{5}}{20} \, M_{\lambda_4}\bigr(\; i\,(-\,b_1 + \sqrt{2}\,a_4) \;,\;  i\,(b_2 + \sqrt{2}\,\overline{a_3}) \;,\; i\,(\overline{b_3} + \sqrt{2}\,a_1) \;\bigr) $$
must vanish, and thus we have
\begin{equation}
\label{eq:EIII:rk1:arctan12:ab}
b_1 = \sqrt{2}\,a_4\;,\quad b_2 = -\sqrt{2}\,\overline{a_3}\;,\quad b_3 = -\sqrt{2}\,\overline{a_1} \; .
\end{equation}
By Equation~\eqref{eq:EIII:rk1:arctan12:malpha'skew} we have \,$\Phi_\alpha(M_{\lambda_3}(b_1,b_2,b_3)) = \sqrt{2/3}\,M_{\lambda_1}(a_1,a_2,a_3,a_4)$\,;
because \,$\Phi_\alpha$\, is isometric, it follows that
$$ \tfrac23\,\sum_k |a_k|^2 = \sum_k |b_k|^2 \overset{\eqref{eq:EIII:rk1:arctan12:ab}} = 2\,(|a_4|^2 + |a_3|^2 + |a_1|^2) $$
and hence
$$ |a_2|^2 = 2\,(|a_1|^2 + |a_3|^2 + |a_4|^2) $$
holds. It follows that the projection map \,$\liem_{\alpha}' \to \C,\; v = M_{\lambda_1}(a_1,a_2,a_3,a_4) + M_{\lambda_3}(b_1,b_2,b_3) \;\mapsto\; a_2$\, 
is injective, and hence we have \,$n_\alpha' \leq 2$\,. We now give \,$v_R = R(H,v_{2\alpha})v$\, explicitly for the situation where \,$v$\, satisfies
Equations~\eqref{eq:EIII:rk1:arctan12:ab}:
$$ v_R = \tfrac{\sqrt{5}}{20}\,\bigr(\,M_{\lambda_1}(i\,a_1, i\,a_2, i\,a_3, -i\,a_4) + \sqrt{2}\,M_{\lambda_3}(-i\,a_4, i\,\overline{a_3}, i\,\overline{a_1})\,\bigr) \; . $$
Because \,$v_R \in \liem_{\alpha}'$\, is therefore orthogonal to \,$v$\,, we see that \,$n_\alpha' \in \{0,2\}$\, holds. 
If \,$n_\alpha' = 0$\,, then we have \,$\liem' = \R\,H \oplus \R\,v_{2\alpha}$\,, and thus we see that \,$\liem'$\, is of type \,$(G_2\C^6,(\PP,\vi=\arctan(\tfrac12),(\C,1)))$\,. 
If \,$n_\alpha' = 2$\, then by Proposition~\ref{P:EIII:isotropy} 
we may suppose without loss of generality that \,$v=M_{\lambda_1}(1,\sqrt{2},0,0) + M_{\lambda_3}(0,0,-\sqrt{2})$\, holds; then
we have \,$\liem' = \R\,H \oplus \R\,v \oplus \R\,v_R \oplus \R\,v_{2\alpha}$\, with \,$4\,v_R=M_{\lambda_1}(i,i\,\sqrt{2},0,0) + M_{\lambda_3}(0,0,\sqrt{2}\,i)$\,.
Therefore \,$\liem'$\, is then of type \,$(G_2\C^6,(\PP,\vi=\arctan(\tfrac12),(\C,2)))$\,. 

Let us now consider the case \,$2\alpha\not\in\Delta'$\,, so that \,$\liem' = \R H \oplus \liem_{\alpha}'$\, holds.
If \,$\alpha\not\in\Delta$\,, then \,$\liem' = \R H$\, is of type \,$(\mathrm{Geo},\vi=\arctan(\tfrac12))$\,, otherwise 
because of Proposition~\ref{P:EIII:isotropy}, we may suppose without loss of generality that
$$ v_\alpha := M_{\lambda_1}(1,\sqrt{2},0,0) + M_{\lambda_3}(0,0,-\sqrt{2}) \;\in\; \liem_\alpha' \; . $$
If \,$n_\alpha'=1$\, holds, then we have \,$\liem' = \R H \oplus \R v_\alpha$\,, and therefore \,$\liem'$\, is then of type \,$(G_2\HH^4,(\Sph,\vi=\arctan(\tfrac13),2))$\,
(note that the isotropy angle \,$\vi=\arctan(\tfrac12)$\, of \,$\liem'$\, corresponds to the isotropy angle \,$\tfrac\pi4-\arctan(\tfrac12)=\arctan(\tfrac13)$\,
within the type \,$(G_2\HH^4)$\, by Remark~\ref{R:EIII:EIII:phirel}).
Otherwise, we let another vector \,$v\in\liem_\alpha'$\, which is orthogonal to \,$v_\alpha$\, be given, say \,$v = M_{\lambda_1}(a_1,a_2,a_3,a_4) + M_{\lambda_3}(b_1,b_2,b_3)$\,
with \,$a_k,b_\ell\in\C$\,, and consider \,$v_R := R(H,v_\alpha)v \in \liem'$\,. Both the \,$\liem_{\lambda_2}$-component and the \,$\liem_{2\lambda_1}$-component must
vanish (because of \,$2\alpha\not\in\Delta'$\,). The \,$\liem_{\lambda_2}$-component of \,$v_R$\, is proportional to 
$$ M_{\lambda_2}\bigr(\; i\,(\sqrt{2}\,b_3+2\,a_1) \;,\; i\,(2\,\overline{b_3}+2\,a_2) \;,\; 
i\,(\sqrt{2}\,\overline{b_2}-2\,a_3-2\,b_1) \;,\; i\,(-2\,\overline{b_2} - 2\,a_4 - \sqrt{2}\,b_1) \;\bigr) $$
and so we have
$$ b_3 = -\sqrt{2}\,a_1\,,\quad b_3 = -\overline{a_2}\,,\quad -2\,b_1 + \sqrt{2}\,\overline{b_2} = 2\,a_3 \qmq{and} -\sqrt{2}\,b_1 - 2\,\overline{b_2} = 2\,a_4 \;, $$
hence
\begin{equation}
\label{eq:EIII:arctan12:cond1}
b_1 = - \tfrac23\,a_3 -\tfrac{\sqrt{2}}{3}\,a_4 \,,\quad b_2 = \tfrac{\sqrt{2}}{3}\,\overline{a_3} - \tfrac23\,\overline{a_4} \,,\quad
b_3 = -\overline{a_2} \qmq{and} a_2 = \sqrt{2}\,\overline{a_1} \; . 
\end{equation}
Moreover the \,$\liem_{2\lambda_1}$-component of \,$v_R$\, is proportional to \,$M_{2\lambda_1}(\IM(a_1-\sqrt{2}\,a_2))$\, and so we have
\,$\IM(a_1-\sqrt{2}\,a_2)=0$\,, hence \,$\IM(a_1) = \sqrt{2}\,\IM(a_2) \overset{\eqref{eq:EIII:arctan12:cond1}}= -2\,\IM(a_1)$\,, and thus
\begin{equation}
\label{eq:EIII:arctan12:cond2}
\IM(a_1) = \IM(a_2) = 0 \; .
\end{equation}
Further, the condition that \,$v$\, is orthogonal to \,$v_\alpha$\, gives
$$ 0 = \g{v}{v_\alpha} = \RE(a_1) + \sqrt{2}\,\RE(a_2) -\sqrt{2}\,\RE(b_3) \overset{\eqref{eq:EIII:arctan12:cond1}}{=} \RE(a_1-\sqrt{2}\,a_2) $$
and therefore \,$\RE(a_1) = \sqrt{2}\,\RE(a_2) \overset{\eqref{eq:EIII:arctan12:cond1}}= 2\,\RE(a_1)$\,, hence
\begin{equation}
\label{eq:EIII:arctan12:cond3}
\RE(a_1) = \RE(a_2) = 0 \; .
\end{equation}
From Equations~\eqref{eq:EIII:arctan12:cond2} and \eqref{eq:EIII:arctan12:cond3} we obtain \,$a_1=a_2=0$\,.
By the remaining equations
from \eqref{eq:EIII:arctan12:cond1} we now see that
$$ v = M_{\lambda_1}\bigr(\,0,0,c,d\,\bigr) 
\;+\; M_{\lambda_3}\bigr(\; -\tfrac23\,c-\tfrac{\sqrt{2}}3\,d \,,\, \tfrac{\sqrt{2}}{3}\,\overline{c} - \tfrac23\,\overline{d} \,,\, 0 \;\bigr) $$
holds with some constants \,$c,d\in\C$\,. 

We now consider the Lie subgroup \,$\lieb := \R\,\alpha_4^\sharp \oplus K_{\alpha_4}(\C)$\, of \,$\liek_0$\,, which is isomorphic to \,$\liesu(2)$\,. 
For \,$z\in\C$\,, we have \,$\ad(K_{\alpha_4}(z))H = \ad(K_{\alpha_4}(z))v_\alpha=0$\,, whereas \,$\ad(K_{\alpha_4}(z))$\, acts on the complex plane
$$ \liew := \Mengegr{\; M_{\lambda_1}(0,0,c,d)
+M_{\lambda_3}(\, -\tfrac23\,c-\tfrac{\sqrt{2}}3\,d \,,\, \tfrac{\sqrt{2}}{3}\,\overline{c} - \tfrac23\,\overline{d} \,,\, 0 \,)\;}{c,d\in\C} $$
as a skew-adjoint, invertible endomorphism for \,$z\neq 0$\,. It follows that the adjoint action of the connected Lie group \,$B\subset K$\,
with Lie algebra \,$\lieb$\, on \,$\liem$\, leaves \,$H$\, and \,$v_\alpha$\, invariant, whereas it acts on \,$\liew$\, as \,$\SU(2)$\, does. 
Therefore there exists \,$g\in B$\, so that \,$\Ad(g)$\, leaves \,$H$\, and \,$v_\alpha$\, invariant, and satisfies
\,$\Ad(g)v = M_{\lambda_1}(0,0,3,0) + M_{\lambda_3}(-2,\sqrt{2},0)$\,. By replacing \,$\liem'$\, with the Lie triple system \,$\Ad(g)\liem'$\,
from the same congruence class, we can thus arrange that
$$ v = M_{\lambda_1}(0,0,3,0) + M_{\lambda_3}(-2,\sqrt{2},0) $$
holds. Hence we see that in the case \,$n_\alpha'=2$\,, \,$\liem' = \R H \oplus \R v_\alpha \oplus \R v$\, is of type \,$(G_2\HH^4,(\Sph,\vi=\arctan(\tfrac13),2))$\,.

Finally we show that the case \,$n_\alpha' \geq 3$\, cannot happen: Let \,$v' \in \liem_\alpha'$\, be orthogonal to both \,$v_\alpha$\, and \,$v$\,. Then, 
as above, the \,$\liem_{\lambda_2}$-component and the \,$\liem_{2\lambda_1}$-component of both \,$R(H,v_\alpha)v'$\, and \,$R(H,v)v'$\, have to vanish, 
and these conditions yield \,$v'=0$\,. 

\textbf{The case \,$\boldsymbol{\vi_0=\tfrac\pi4}$\,.}
In this case we have
by Equation~\eqref{eq:EIII:rk1:H}: \,$H = \tfrac{1}{\sqrt{2}}\,\lambda_2^\sharp + \tfrac{1}{\sqrt{2}}\,\lambda_1^\sharp$\, and therefore
$$ \lambda_1(H) = \tfrac{1}{\sqrt{2}},\; \lambda_2(H) = \tfrac{1}{\sqrt{2}},\; \lambda_3(H) = 0,\; \lambda_4(H) = \tfrac{2}{\sqrt{2}},\;
2\lambda_1(H) = \tfrac{2}{\sqrt{2}},\; 2\lambda_2(H) = \tfrac{2}{\sqrt{2}}\;. $$
It follows
by Equation~\eqref{eq:EIII:rk1:Delta'Delta} that we have \,$\Delta' \subset \{\pm \alpha,\pm 2\alpha\}$\, with \,$\alpha := \lambda_1|\liea' = \lambda_2|\liea'$\,,
\,$2\alpha = \lambda_4|\liea' = (2\lambda_1)|\liea' = (2\lambda_2)|\liea'$\,,
and by Equations~\eqref{eq:EIII:rk1:m'decomp},\eqref{eq:EIII:rk1:malpha'} we have
\begin{equation}
\label{eq:EIII:rk1:pi4:m'decomp}
\liem' = \R H \oplus \liem_\alpha' \oplus \liem_{2\alpha}'
\end{equation}
with \,$\liem_\alpha' \subset \liem_{\lambda_1} \oplus \liem_{\lambda_2}$\, and \,$\liem_{2\alpha}' \subset \liem_{2\lambda_1} \oplus \liem_{2\lambda_2} \oplus \liem_{\lambda_4}$\,. 

Further information on the structure of elements of \,$\liem_\alpha'$\, resp.~of \,$\liem_{2\alpha}'$\, can be obtained: First let \,$v \in \liem_\alpha'$\, be given,
say \,$v = M_{\lambda_1}(a_1,a_2,a_3,a_4) + M_{\lambda_2}(b_1,b_2,b_3,b_4)$\, with \,$a_k,b_k \in \C$\,. By Proposition~\ref{P:cla:skew} and the fact that
\,$\alpha^\sharp = \tfrac12\,(\lambda_1^\sharp + \lambda_2^\sharp)$\, holds, we see that we have \,$\|(a_1,\dotsc,a_4)\| = \|(b_1,\dotsc,b_4)\|$\,, in particular
\,$n_\alpha' \leq 8$\,. 

Similarly, for any \,$v \in \liem_{2\alpha}'$\,, say \,$v = M_{\lambda_4}(c_1,c_2,c_3) + M_{2\lambda_1}(t) + M_{2\lambda_2}(s)$\, with \,$c_k\in\C$\, and \,$t,s\in\R$\,,
we consider the vector \,$v_R := R(H,v)v \in \liem'$\,. The \,$\liea$-component of \,$v_R$\, must be proportional to \,$H$\,, and this condition
yields \,$|t|=|s|$\,, hence \,$t = \pm s$\,. Moreover, because of \,$\lambda_3(H)=0$\,, the \,$\lambda_3$-component of \,$v_R$\,, which is proportional to
$$ M_{\lambda_3}(\,ic_1(s-t) \,,\, ic_2(t-s) \,,\, i\overline{c_3}(s-t)\,) \; , $$
has to vanish, and thus we have either \,$c_1=c_2=c_3=0$\, or \,$t=s$\,. If we have \,$t=-s$\,, and hence \,$c_1=c_2=c_3=0$\,, 
we put \,$Y := K_{2\lambda_1}(\sqrt{8}) - K_{2\lambda_2}(\sqrt{8})$\,, 
then we have \,$\ad(Y)H = M_{2\lambda_1}(1) - M_{2\lambda_2}(1)$\, and \,$\ad(Y)(M_{2\lambda_1}(1)+M_{2\lambda_2}(1))
= 4\sqrt{2}\,\lambda_3^\sharp$\,. 
These equations show that a Lie triple system \,$\liem'$\, where the case \,$t=-s$\, occurs is congruent under the adjoint action of a member of the
1-parameter subgroup of \,$K$\, induced by \,$Y$\, to a Lie triple system corresponding to the case \,$t=s$\,. By replacing \,$\liem'$\, with the latter
Lie triple system, we may suppose without loss of generality that in any case \,$\liem'_{2\alpha} \subset \liem_{\lambda_4} \oplus \R(M_{2\lambda_1}(1)+M_{2\lambda_2}(1))
=: \wh{\liem}_{2\alpha}'$\, holds.

In the case \,$\alpha\not\in\Delta'$\, it now follows immediately that \,$\liem'$\, is
of type \,$(\PP,\vi=\tfrac\pi4, \Sph^{1+n_{2\alpha}'})$\,.

So let us now turn our attention to the case \,$\alpha\in\Delta'$\,. \,$\liem'$\, corresponds to a Riemannian symmetric space of rank \,$1$\,; the classification
of these spaces gives that we have \,$n_{2\alpha}' \in \{0,1,3,7\}$\, (corresponding to the projective spaces over the reals, the complex numbers, the quaternions,
and the octonions, respectively), and that \,$n_{2\alpha}'+1$\, divides \,$n_{\alpha}'$\,. 

We continue our investigation of the structure of \,$\liem_\alpha'$\,:
Let \,$v \in \liem_{\alpha}'$\, be given, say \,$v = M_{\lambda_1}(a_1,\dotsc,a_4)
+ M_{\lambda_2}(b_1,\dotsc,b_4)$\, with \,$a_k,b_k \in \C$\,. Then the \,$\liem_{\lambda_3}$-component of \,$R(H,v)v \in \liem'$\, equals
$$ M_{\lambda_3}\bigr( \; -\tfrac{i}{8}(a_4\,b_1 + a_1\,b_4+b_2\,a_3+a_2\,b_3) \;,\;
-\tfrac{i}{8}(\overline{a_4}\,b_2 - \overline{a_3}\,b_1 + a_1\,\overline{b_3} - a_2\,\overline{b_4}) \;,\;
\tfrac{i}{8}(b_1\,\overline{a_1} - \overline{a_4}\,b_4 + \overline{b_3}\,a_3 - \overline{b_2}a_2) \; \bigr) \; . $$
Because of \,$\lambda_3(H)=0$\,, this has to vanish. In this way it follows that 
$$ \liem_\alpha' \;\subset\; \Mengegr{ M_{\lambda_1}(a_1,a_2,a_3,a_4) + M_{\lambda_2}(a_2,a_1,-a_4,-a_3)}{a_1,\dotsc,a_4 \in \C} \;=:\; \wh{\liem}_\alpha' $$
holds.

Therefore in any case \,$\liem'$\, is contained in the Lie triple system \,$\wh{\liem}' := \R\,H \oplus \wh{\liem}_\alpha' \oplus
\wh{\liem}_{2\alpha}'$\, of type \,$(\PP,\vi=\tfrac\pi4,\OP^2)$\,. 
The totally geodesic submanifold corresponding to \,$\wh{\liem}'$\, is a Cayley plane \,$\OP^2$\,, and \,$\liem'$\,
also is a Lie triple system of \,$\wh{\liem}'$\,. Therefore it follows from the classification of the Lie triple systems of \,$\OP^2$\, (see \cite{Wolf:1963-elliptic},
Section~3), that \,$\liem'$\, is of one of the types \,$(\PP,\vi=\tfrac\pi4,\KP^2)$\, with \,$\K\in\{\R,\C,\HH,\OO\}$\,. 

This completes the proof of the classification of Lie triple systems in \,$\EIII$\,. 

\subsection[Totally geodesic submanifolds in \,$E_6/(\Ug(1)\cdot\Spin(10))$\,]{Totally geodesic submanifolds in \,$\boldsymbol{E_6/(\Ug(1)\cdot\Spin(10))}$\,}
\label{SSe:EIII:tgsub}

We now study the geometry of the totally geodesic submanifolds of \,$\EIII$\, associated to the Lie triple systems found in Theorem~\ref{EIII:EIII:cla}.
Of course, the local isometry type of the submanifolds can easily be obtained by determining the restricted root system (with multiplicities) of the
corresponding Lie triple systems as they are given in Theorem~\ref{EIII:EIII:cla}. But to obtain the global isometry type, and also to understand the
position of the submanifolds in \,$\EIII$\, better, we describe the totally geodesic submanifolds of \,$\EIII$\, explicitly. We will also want information
on how the transvection groups of the respective submanifolds are embedded in \,$\Esix$\,, the transvection group of \,$\EIII$\,. 

In this way, we obtain the results of the following table. Herein, we ascribe the type of a Lie triple system also to the corresponding totally geodesic
submanifold (or to a corresponding totally geodesic embedding). For \,$\ell\in\N$\, and \,$r > 0$\, we denote by \,$\Sph^\ell_r$\, the \,$\ell$-dimensional
sphere of radius \,$r$\,, and for \,$\vkap > 0$\, we denote by \,$\RP^\ell_\vkap$\, the \,$\ell$-dimensional real projective space
of sectional curvature \,$\vkap$\,, and by \,$\CP^\ell_\vkap$\, the \,$\ell$-dimensional complex projective space of constant holomorphic curvature \,$4\vkap$\,.
Note that with these notations, \,$\RP^\ell_\vkap$\, is a real form of \,$\CP^\ell_\vkap$\,. 
Moreover, with \,$\HP^\ell_\vkap$\, resp.~\,$\OP^2_\vkap$\, we denote the \,$\ell$-dimensional quaternionic projective space resp.~the Cayley projective plane,
with their invariant Riemannian metric scaled in such a way that the \emph{minimal} sectional curvature equals \,$\vkap$\,. 
Also for the irreducible Riemannian symmetric spaces of rank \,$2$\,, their invariant Riemannian
metric is a priori only defined up to a positive constant; in the table below we describe the appropriate metric of these spaces by giving the
length \,$a$\, of the shortest restricted root of the space in the index \,${}_{\mathrm{srr}=a}$\,. We continue to use also the additional names of types
introduced in Remark~\ref{R:EIII:EIII:extratypes}.

\begin{center}
\begin{longtable}{|c|c|c|}
\hline
type of Lie triple system & isometry type & properties\footnotemark \\
\hline
\endhead
\hline
\endfoot
$(\mathrm{Geo},\vi=t)$ & \,$\R$\, or \,$\Sph^1$ & \\
$(\PP,\vi=0,(\K,\ell))$ & $\KP^\ell_{\vkap=1}$ & \,$(\K,\ell)=(\C,1)$\,: Helgason sphere \\
$(\PP,\vi=\tfrac\pi4,\Sph^k)$ & $\Sph^k_{r=1/\sqrt{2}}$ & \\
$(\PP,\vi=\tfrac\pi4,\KP^2)$ & $\KP^2_{\vkap=1/2}$ & \,$\K=\OO$\,: reflective, real form, maximal \\
\hline
$(\PP\times\PP^1,(\K_1,\ell),\K_2)$ & $\K_1\mathrm{P}^\ell_{\vkap=1} \times \K_2\mathrm{P}^1_{\vkap=1}$ & \,$(\K_1,\ell,\K_2)=(\C,5,\C)$\,: meridian for \,$(\mathrm{DIII})$\,, maximal \\
$(Q)$ & $Q^8_{\mathrm{srr}=\sqrt{2}}$ & polar, meridian for itself, maximal \\
$(Q,\tau)$ & see \cite{Klein:2007-claQ}, Section~5 & \\
$(G_2\C^6)$ & \,$G_2(\C^6)_{\mathrm{srr}=1}$\, & reflective, maximal \\
$(G_2\C^6,\tau)$ & see \cite{Klein:2007-tgG2}, Section~7 & \\
$(G_2\HH^4)$ & $(G_2(\HH^4)/\Z_2)_{\mathrm{srr}=1}$\, & reflective, real form, maximal \\
$(G_2\HH^4,\tau)$ & see \cite{Klein:2007-tgG2}, Section~6 & \\
$(\mathrm{DIII})$ & $\SO(10)/\Ug(5)_{\mathrm{srr}=1}$\, & polar, maximal \\
\end{longtable}
\end{center}
\footnotetext{The polars and meridians are also reflective, without this fact being noted explicitly in the table.}

For the application of the information from \cite{Klein:2007-claQ} and \cite{Klein:2007-tgG2} it should be noted that these two papers use different conventions
regarding the metric used on the spaces under investigation: In the investigation of the complex quadrics in \cite{Klein:2007-claQ}, the metric is normalized
such that the shortest restricted roots of \,$Q^m$\, have length \,$\sqrt{2}$\,, whereas in the investigation of \,$G_2(\C^n)$\, and \,$G_2(\HH^n)$\,
in \cite{Klein:2007-tgG2}, the metric on these spaces is normalized such that those linear forms which are the shortest roots in 
\,$G_2(\K^n)$\, for \,$n\geq 5$\, have length \,$1$\,
(notice that they are not actually roots of \,$G_2(\K^4)$\,, because their multiplicities then degenerate to zero).
Also for the investigation of \,$\EIII$\, in the present paper, we normalize the metric such that the shortest roots of this space have length \,$1$\,. 

By looking at the root systems of the totally geodesic submanifolds of type \,$(Q)$\,, \,$(G_2\C^6)$\, and \,$(G_2\HH^4)$\, of \,$\EIII$\, 
(see Remark~\ref{R:EIII:EIII:phirel}),
it follows that the data given in the cited papers on the metric properties of totally geodesic submanifolds can be carried over without any change to the present situation
for the totally geodesic submanifolds \,$Q^8$\, and \,$G_2(\C^6)$\, of \,$\EIII$\,. However, for \,$G_2(\HH^4)/\Z_2$\, it is necessary to scale the data given
in \cite{Klein:2007-tgG2}, as this manifold is considered with \,$\mathrm{srr}=\sqrt{2}$\, in \cite{Klein:2007-tgG2}, whereas it has \,$\mathrm{srr}=1$\, here.

For the proof of the data in the table, and to obtain the desired information on the position of the totally geodesic submanifolds of \,$\EIII$\,,
it is sufficient to consider the maximal totally geodesic submanifolds. In the case of \,$\EIII$\, every maximal totally geodesic submanifold is reflective
(see~\cite{Leung:reflective-1979}), and therefore a connected component of the fixed point set of an involutive isometry of \,$\EIII$\,. We will
describe these submanifolds in this way in the first instance.

To prove that the fixed point sets of the involutive isometries of \,$\EIII$\, we investigate below are indeed of the isometry type claimed above,
we will then construct totally geodesic, equivariant embeddings of the appropriate manifolds onto these fixed point sets for many of the
types of maximal totally geodesic submanifolds of \,$\EIII$\,. We will also describe the subgroups of the transvection group \,$\Esix$\, of \,$\EIII$\,
which correspond to the transvection groups of these totally geodesic submanifolds.

For these investigations, we need a model of \,$\EIII$\, in which we can carry out calculations explicitly. For this purpose, we use the 
explicit presentations of \,$\EIII$\, and of the exceptional Lie group \,$\Esix$\,
given by \textsc{Yokota} in \cite{Yokota:invol1-1990} and by \textsc{Atsuyama} in \cite{Atsuyama:projspaces2-1997}.

To describe these presentations,
we denote by \,$\R$\,, \,$\C = \R\oplus \R i$\,, \,$\HH = \C \oplus \C j$\, and \,$\OO = \HH \oplus \HH e$\, 
the four normed real division algebras: the field of real numbers, the field of complex numbers,
the skew-field of quaternions, and the division algebra of octonions. For \,$\K\in\{\C,\HH,\OO\}$\, and \,$x\in\K$\,,
we have the conjugate \,$\overline{x}$\, of \,$x$\,. We will also consider the complexification \,$\K^{\C} := \K \otimes_{\R} \R^{\C}$\,
of \,$\K$\, with respect to another ``copy'' \,$\R^{\C} = \R \oplus \R I$\, of the field of complex numbers;
we linearly extend the conjugation map \,$x\mapsto \overline{x}$\, of \,$\K$\, to \,$\K^{\C}$\,. 
Notice that the algebras \,$\C^{\C}$\,, \,$\HH^{\C}$\, and \,$\OO^{\C}$\, have zero divisors.

Let \,$M(n\times m,\K)$\, be the linear space of \,$(n\times m)$-matrices over \,$\K$\,, abbreviate \,$M(n,\K) := M(n\times n,\K)$\,, and let
\,$\frakJ(n,\K) := \Menge{X\in M(n,\K)}{X^* = X}$\, be the subspace of Hermitian matrices; via the multiplication map \,$\frakJ(n,\K) \times \frakJ(n,\K)
\to \frakJ(n,\K),\; (X,Y) \mapsto X\circ Y := \tfrac12(XY + YX)$\,, \,$\frakJ(n,\K)$\, becomes a real Jordan algebra for \,$\K \in \{\R,\C,\HH\}$\, or
\,$\K=\OO$\,,\,$n=3$\,; it becomes a complex Jordan algebra for \,$\K \in \{\R^{\C},\C^{\C},\HH^{\C}\}$\, or \,$\K=\OO^{\C}$\,,\,$n=3$\,. 
\,$\frakJ(3,\OO)$\, resp.~\,$\frakJ := \frakJ(3,\OO^{\C})$\, is the real resp.~complex exceptional Jordan algebra. 

We now consider the complex projective space over \,$\frakJ$\,, which we denote by \,$[\frakJ] \cong \CP^{26}$\,. For \,$X \in \frakJ \setminus \{0\}$\,, we denote
by \,$[X] := (\R^{\C}) X$\, the projective line through \,$X$\,; for a subset \,$M \subset \frakJ$\,, we put \,$[M] := \Menge{[X]}{X \in M \setminus \{0\}}$\,. 
Following Atsuyama (\cite{Atsuyama:projspaces2-1997}), we consider the submanifold
$$ \wt{\EIII} := 
\left\{ \left. X = \left( \begin{matrix} \xi_1 & x_3 & \overline{x_2} \\ \overline{x_3} & \xi_2 & x_1 \\ x_2 & \overline{x_1} & \xi_3 \end{matrix} \right)
\in \frakJ \;\right|\; \begin{matrix} 
\xi_1,\xi_2,\xi_3 \in \R^{\C}, \; x_1,x_2,x_3 \in \OO^{\C} \\
\xi_2\,\xi_3 = |x_1|^2,\;\;\;\xi_3\,\xi_1 = |x_2|^2,\;\;\; \xi_1\,\xi_2 = |x_3|^2, \\ 
x_2\,x_3 = \xi_1\,\overline{x_1},\;\;\; x_3\,x_1 = \xi_2\,\overline{x_2},\;\;\; x_1\,x_2 = \xi_3\,\overline{x_3}
\end{matrix} \right\} $$
of \,$\frakJ$\,. Then Atsuyama has shown (\cite{Atsuyama:projspaces2-1997}, Lemma~3.1) that \,$[\wt{\EIII}] \subset [\frakJ]$\, is a model of the
exceptional symmetric space \,$\EIII$\,. In the sequel, we denote by \,$\EIII$\, this model.

We will also use the fact that the exceptional Lie group \,$\Esix$\,, which is the transvection group of \,$\EIII$\,, can be realized as a subgroup of the 
group \,$\Aut(\frakJ)$\, of complex-linear automorphisms of \,$(\frakJ,\circ)$\,. More specifically, consider the
inner product \,$\g{\cdot}{\cdot}$\, and the operation \,$A \,\Delta\, B$\, defined on \,$\frakJ$\, in \cite{Atsuyama:projspaces2-1997}, \S 1. 
Then Atsuyama showed in \cite{Atsuyama:projspaces2-1997}, Lemma~1.5(2) that 
$$ \Esix = \Menge{f \in \Aut(\frakJ)}{\;\forall X,Y\in\frakJ: f(X\,\Delta\,Y) = (fX)\,\Delta\,(fY),\; \g{fX}{fY} = \g{X}{Y}\;} $$
is a model of the exceptional Lie group \,$\Esix$\,. This model acts transitively on the model of \,$\EIII$\, described above.

We now define several involutive isometries on \,$\EIII$\, (see also \cite{Yokota:invol1-1990} Section~3, where the involutive automorphisms on the
exceptional Lie group \,$\Esix$\, are classified):

\begin{itemize}
\item
The conjugation map \,$\lambda_0: \OO^{\C} \to \OO^{\C}$\, induced by the real form \,$\OO$\, of \,$\OO^{\C}$\,
(i.e.~the orthogonal involution \,$\lambda_0: \OO^{\C}\to\OO^{\C}$\, characterized by \,$\Fix(\lambda_0) = \OO$\,) induces via the map
$$ \wt{\EIII} \to \wt{\EIII},\; \left( \begin{smallmatrix} \xi_1 & x_3 & \overline{x_2} \\ \overline{x_3} & \xi_2 & x_1 \\ x_2 & \overline{x_1} & \xi_3 \end{smallmatrix} \right)
\mapsto \left( \begin{smallmatrix} \lambda_0(\xi_1) & \lambda_0(x_3) & \overline{\lambda_0(x_2)} \\ \overline{\lambda_0(x_3)} & \lambda_0(\xi_2) & \lambda_0(x_1) \\ \lambda_0(x_2) & \overline{\lambda_0(x_1)} & \lambda_0(\xi_3) \end{smallmatrix} \right) $$
an isometric involution \,$\lambda: \EIII \to \EIII$\,. 
\item The orthogonal involution \,$\gamma_0 : \OO^{\C} \to \OO^{\C}$\, characterized by \,$\Fix(\gamma_0) = \HH^{\C}$\, 
induces via the map
$$ \wt{\EIII} \to \wt{\EIII},\; \left( \begin{smallmatrix} \xi_1 & x_3 & \overline{x_2} \\ \overline{x_3} & \xi_2 & x_1 \\ x_2 & \overline{x_1} & \xi_3 \end{smallmatrix} \right)
\mapsto \left( \begin{smallmatrix} \xi_1 & \gamma_0(x_3) & \overline{\gamma_0(x_2)} \\ \overline{\gamma_0(x_3)} & \xi_2 & \gamma_0(x_1) \\ \gamma_0(x_2) & \overline{\gamma_0(x_1)} & \xi_3 \end{smallmatrix} \right) $$
another isometric involution \,$\gamma: \EIII \to \EIII$\,.
\item The linear map 
$$ \wt{\EIII} \to \wt{\EIII},\; \left( \begin{smallmatrix} \xi_1 & x_3 & \overline{x_2} \\ \overline{x_3} & \xi_2 & x_1 \\ x_2 & \overline{x_1} & \xi_3 \end{smallmatrix} \right)
\mapsto \left( \begin{smallmatrix} \xi_1 & -x_3 & -\overline{x_2} \\ -\overline{x_3} & \xi_2 & x_1 \\ -x_2 & \overline{x_1} & \xi_3 \end{smallmatrix} \right) $$
induces yet another isometric involution \,$\sigma: \EIII \to \EIII$\,. 
\,$\sigma$\, is the geodesic symmetry of the symmetric space \,$\EIII$\, at the point \,$p_0 := \left[ \begin{smallmatrix} 1 & 0 & 0 \\ 0 & 0 & 0 \\ 0 & 0 & 0 \end{smallmatrix} \right]
\in \EIII$\,. 
\end{itemize}

With help of these involutions we can describe the reflective submanifolds of \,$\EIII$\, explicitly.

\paragraph{The type \,$\boldsymbol{(\PP,\vi=\tfrac\pi4,\OP^2)}$\,.}
The fixed point set of \,$\lambda$\, equals \,$[\EIII \cap \frakJ(3,\OO)] \cong \OP^2$\,, a totally geodesic submanifold of \,$\EIII$\, of type
\,$(\PP,\vi=\tfrac\pi4,\OP^2)$\,. Notice that this is a real form of the Hermitian symmetric space \,$\EIII$\,. 

\paragraph{The types \,$\boldsymbol{(G_2\C^6)}$\, and \,$\boldsymbol{(\PP\times\PP^1,(\C,5),\C)}$\,.}
The fixed point set of the involutive isometry \,$\gamma: \EIII \to \EIII$\, has two connected components:
\begin{gather*}
F^\gamma_1 := [\wt{\EIII} \cap \frakJ(3,\HH^{\C})] \;,\\
\qmq{and} F^\gamma_2 := \left\{\left.\left[ \begin{smallmatrix} 0 & a_3\,e & -a_2\,e \\ -a_3\,e & 0 & a_1\,e \\ a_2\,e & -a_1\,e & 0 \end{smallmatrix} \right] \;\right|\;
a_k \in \HH^{\C},\; a_1\,\overline{a_2} = a_2\,\overline{a_3} = a_3\,\overline{a_1} = 0 \right\} \; . 
\end{gather*}

It turns out that the totally geodesic submanifolds \,$F^\gamma_1$\, and \,$F^\gamma_2$\, of \,$\EIII$\, are of type \,$(G_2\C^6)$\, and 
\,$(\PP\times \PP^1,(\C,5),\C)$\,, respectively. To show that they are isomorphic to \,$G_2(\C^6)$\, resp.~to \,$\CP^1 \times \CP^5$\,,
we will now explicitly construct isometries \,$f_1: G_2(\C^6) \to F^\gamma_1$\, and \,$f_2: \CP^1 \times \CP^5 \to F^\gamma_2$\, which are compatible
with the group actions on the symmetric spaces involved.

For this purpose, we note that \,$\Esix$\, contains a subgroup which is isomorphic to \,$(\Sp(1)\times\SU(6))/\Z_2$\,, and which is the fixed point group of
the Lie group automorphism \,$\Esix \to \Esix,\; g \mapsto \gamma\cdot g \cdot \gamma^{-1}$\,. This subgroup has been described explicitly by
Yokota (\cite{Yokota:invol1-1990}, Section~3.5) in the following way:

To associate to a given \,$(b,B) \in \Sp(1) \times \SU(6)$\, a member of \,$\Esix \subset \Aut(\frakJ)$\,, we need to describe an action
of \,$(b,B)$\, on \,$\frakJ$\,. For this purpose we note that \,$\frakJ$\, is \,$(\R^{\C})$-linear isomorphic to \,$\frakJ(3,\HH^{\C}) \oplus (\HH^{\C})^3$\, by the map
$$ \vi_1: \frakJ(3,\HH^{\C}) \oplus (\HH^{\C})^3 \to \frakJ,\; (X,x) \mapsto X + \left( \begin{smallmatrix} 0 & x_3\,e & x_2\,e \\ -x_3\,e & 0 & x_1\, e \\ x_2\,e & -x_1\,e & 0 \end{smallmatrix} \right) \; . $$
Furthermore, \,$M(3,\HH^{\C}) \supset \frakJ(3,\HH^{\C})$\, is \,$(\R^{\C})$-linear isomorphic to \,$M(6,\C^{\C})_J := \Menge{X \in M(6,\C^{\C})}{JX=\overline{X}J}$\,, and
\,$(\HH^{\C})^3$\, is \,$(\R^{\C})$-linear isomorphic to \,$M(2\times 6,\C^{\C})_J := \Menge{X \in M(2\times 6,\C^{\C})}{J'X=\overline{X}J}$\,, where we put
\,$J' := \left( \begin{smallmatrix} 0 & -1 \\ 1 & 0 \end{smallmatrix} \right) \in M(2,\R)$\, and 
\,$J := \diag(J',J',J') \in M(6,\R)$\,. 
These isomorphisms are exhibited by the maps
$$ \vi_2: M(3,\HH^{\C}) \to M(6,\C^{\C})_J \qmq{resp.} \vi_2': (\HH^{\C})^3 \to M(2\times 6,\C^{\C})_J $$
which transform any given matrix \,$X \in M(3,\HH^{\C})$\, resp.~any given row vector \,$x \in (\HH^{\C})^3$\,
into a matrix \,$\vi_2(X) \in M(6,\C^{\C})$\, resp.~\,$\vi_2'(x) \in M(2\times 6,\C^{\C})$\, by mapping
every entry \,$a+bj \in \HH^{\C}$\, (\,$a,b\in\C^{\C}$\,) of \,$X$\, into a \,$(2\times 2)$-block component
\,$\left( \begin{smallmatrix} a & b \\ -\overline{b} & \overline{a} \end{smallmatrix} \right)$\, of \,$\vi_2(X)$\, resp.~\,$\vi_2'(x)$\,. 
We put \,$\frakJ(6,\C^{\C})_J := \vi_2(\frakJ(3,\HH^{\C})) \subset M(6,\C^{\C})_J$\,. 
In this way we obtain an isomorphism
between \,$\frakJ$\, and \,$\bbV := \frakJ(6,\C^{\C})_J \oplus M(2\times 6,\C^{\C})_J$\,:
$$ \vi := (\vi_2 \oplus \vi_2') \circ \vi_1^{-1} : \frakJ \to \bbV \; , $$
which we will use to describe the action of \,$\Sp(1)\times \SU(6)$\, on \,$\frakJ$\,. 

To do so, we consider for \,$\K\in\{\C,\C^{\C}\}$\, besides 
\,$\SU(6,\K) = \Menge{A\in M(6,\K)}{A^*A = \id,\; \det(A)=1}$\, also
\,$\SU^*(6,\K) := \Menge{A \in M(6,\K)}{JA = \overline{A}J,\; \det(A)=1}$\,. 
Then we have the isomorphism of Lie groups 
$$ \Phi: \SU(6,\C^{\C}) \to \SU^*(6,\C^{\C}),\; A \mapsto \eps\,A - \overline{\eps}\,J\,\overline{A}\,J \;, $$
where we put \,$\eps := \tfrac12\,(1+iI) \in \C^{\C}$\,. 

We now consider the action \,$ F_0: (\Sp(1) \times \SU^*(6,\C^{\C})) \,\times\, \bbV \to \bbV $\, given by
$$ F_0(b,B)(X+x) = BXB^* + (\vi_2'\,b\,(\vi_2')^{-1})xB^{-1} $$
for all \,$(b,B) \in \Sp(1)\times\SU^*(6,\C^{\C})$\, and \,$X+x \in \bbV$\,. \,$F_0$\, induces an action \,$F: (\Sp(1) \times \SU(6)) \,\times\, \frakJ \to \frakJ$\,
which is characterized by the fact that the following diagram commutes:
\begin{equation*}
\begin{minipage}{5cm}
\begin{xy}
\xymatrix{
(\Sp(1)\times\SU^*(6,\C^{\C})) \,\times\, \bbV  \ar[r]^{\hspace{2cm} F_0} & \bbV \\
(\Sp(1)\times\SU(6)) \,\times\, \frakJ \ar[r]_{\hspace{1.5cm} F} \ar[u]^{(\id_{\Sp(1)}\times \Phi)\,\times\,\vi} & \frakJ \ar[u]_{\vi} \;.
}
\end{xy}
\end{minipage}
\end{equation*}
It has been shown by Yokota (\cite{Yokota:invol1-1990}, Theorem~3.5.11 and its proof) that \,$F(b,B)\in\Esix$\, holds for all \,$(b,B) \in \Sp(1)\times\SU(6)$\,. 
In this way we obtain a homomorphism of Lie groups \,$F: \Sp(1)\times\SU(6)\to\Esix$\, with \,$\ker(F) = \{\pm(\id,\id)\}$\,. 

We now denote for \,$U \in G_2(\C^6)$\, by \,$P_U \in M(6,\C)$\, the orthogonal projection onto \,$U$\,. Then we have
\,$Q_U := \eps\,P_U - \overline{\eps}\,J\,\overline{P_U}\,J \in \frakJ(6,\C^{\C})_J$\,, and therefore the map
$$ f_1: G_2(\C^6) \to [\frakJ],\; U \mapsto [\vi^{-1}(Q_U + 0_{M(2\times 6,\C^{\C})})] $$
is well-defined.
It turns out that \,$f_1$\, is an isometric embedding and equivariant in the sense that for every \,$B\in\SU(6),\;U\in G_2(\C^6)$\, we have
$$ F(\id,B)f_1(U) = f_1(B\,U) \; . $$
As a consequence of this property and the fact that \,$f_1(\C e_1 \oplus \C e_2) = p_0 \in \EIII$\, holds, 
\,$f_1$\, maps into \,$\EIII$\,, and hence it maps into \,$\EIII \cap [\frakJ(3,\HH^{\C})] = F_1^\gamma$\,. Because
both \,$F_1^\gamma$\, and \,$G_2(\C^6)$\, are compact and connected, and are of the same (real) dimension \,$16$\,, it follows that
the isometric embedding \,$f_1$\, in fact maps \,$G_2(\C^6)$\,
onto the totally geodesic submanifold \,$F_1^\gamma$\, of \,$\EIII$\,. 

To similarly construct a map \,$f_2: \CP^1 \times \CP^5 \to F^\gamma_2$\,, we identify \,$\C^2$\, with \,$\HH$\,. In this way, 
we can regard \,$\CP^1$\, as the space \,$\Menge{\ell\,\C}{\ell\in\Sph(\HH)}$\,. We also identify \,$\C^6$\, with \,$\HH^3$\,. 
Using these identifications, we can interpret for any \,$\ell\in \C^2 \cong \HH$\, and \,$v \in \C^6 \cong \HH^3$\, the expression \,$\ell\eps v$\, as a member
of \,$(\HH^{\C})^3$\,; via this expression we define the map
$$ f_2: \CP^1 \times \CP^5 \to [\frakJ],\; (\ell\,\C,[v]) \mapsto [\vi^{-1}(0_{\frakJ(6,\C^{\C})_J} + \ell\eps v^*)] \; , $$
which turns out to be a well-defined isometric embedding, which is equivariant in the following sense:
For all \,$(b,B) \in \Sp(1) \times \SU(6)$\,, \,$(\ell\,\C,[v]) \in \CP^1 \times \CP^5$\,, we have
$$ F(b,B)\,f_2(\ell\,\C,[v]) = f_2(b\ell\,\C,[Bv]) \; . $$
Because of this property, and the fact that \,$f_2(1\,\C,[e_1]) = \left[ \begin{smallmatrix} 0 & 0 & 0 \\ 0 & 0 & \eps\,e \\ 0 & -\eps\,e & 0  \end{smallmatrix} \right]
\in F_2^\gamma \subset \EIII$\,,
\,$f_2$\, maps into \,$\EIII$\,, and hence it maps into \,$\EIII \cap [\vi(0_{\frakJ(6,\C^{\C})_J}\oplus (\HH^{\C})^3)] = F^\gamma_2$\,. Because
both \,$F_2^\gamma$\, and \,$\CP^1\times \CP^5$\, are compact and connected, and they are of the same (real) dimension \,$12$\,, it follows that the isometric embedding
\,$f_2$\, in fact maps \,$\CP^1 \times \CP^5$\, onto the totally geodesic submanifold \,$F^\gamma_2$\, of \,$\EIII$\,. 

\paragraph{The type \,$\boldsymbol{(G_2\HH^4)}$\,.}
Notice that the involutive isometries \,$\lambda$\, and \,$\gamma$\, commute with each other, and therefore 
\,$\lambda\circ \gamma$\, is another involutive isometry of \,$\EIII$\,. 
The fixed point set of the latter involution equals
$$
F^{\lambda\gamma} := \left\{\left. p := \left[ \begin{smallmatrix} 
r_1 & p_3 + q_3\,e\,I & \overline{p_2} - q_2\,e\,I \\ \overline{p_3} - q_3 \, e \, I & r_2 & p_1 + q_1 \, e\, I \\ p_2 + q_2\,e\,I & \overline{p_1} -q_1\,e\,I & r_3
\end{smallmatrix} \right] \;\right|\;
\begin{matrix} 
r_k \in \R,\; p_k,q_k \in \HH \\
p \in \EIII
\end{matrix} \right\} \;. 
$$

It turns out that the totally geodesic submanifold \,$F^{\lambda\gamma}$\, of \,$\EIII$\, corresponds to the type \,$(G_2\HH^4)$\,. We will show that
\,$F^{\lambda\gamma}$\, is isometric to \,$G_2(\HH^4)/\Z_2$\,. 

\,$\Esix$\, contains a subgroup isomorphic to \,$\Sp(4)/\Z_2$\,, which is the fixed point group of the Lie group
automorphism \,$\Esix\to\Esix,\; g \mapsto (\lambda\gamma)g(\lambda\gamma)^{-1}$\,. Also this subgroup has been described explicitly by Yokota
(\cite{Yokota:invol1-1990}, Section~3.4). We will use his construction, which we now describe, to obtain an action of \,$\Sp(4)$\, on \,$\frakJ$\,.

We continue to use the space \,$\bbV$\, and the linear isomorphism \,$\vi: \frakJ \to \bbV$\, from the previous construction, 
put \,$\frakJ(4,\HH^{\C})_0 := \Menge{X\in\frakJ(4,\HH^{\C})}{\tr(X)=0}$\,, and consider the isomorphism of linear spaces
\,$\psi: \frakJ \to \frakJ(4,\HH^{\C})_0$\, given in the following way: For \,$A\in\frakJ$\,, say \,$\vi(A) = X+x\in\bbV$\,,
we put
$$\psi(A) = \left( \begin{smallmatrix} \tfrac12\,\tr(X) & I\,x \\ I\,x^* & X - \tfrac12\,\tr(X) \cdot \id_{(\HH^{\C})^3} \end{smallmatrix} \right) \;, $$
where the right-hand expression is to be read as a block matrix with respect to the decomposition \,$(\HH^{\C})^4 = \HH^{\C} \oplus (\HH^{\C})^3$\,. 

Notice that \,$\Sp(4)$\, acts on \,$\frakJ(4,\HH^{\C})_0$\, in the canonical way, i.e.~by the action
$$ F_0: \Sp(4) \times \frakJ(4,\HH^{\C})_0 \to \frakJ(4,\HH^{\C})_0,\; (B,X) \mapsto BXB^* \; . $$
Via the linear isomorphism \,$\psi$\,, \,$F_0$\, induces an action \,$F: \Sp(4) \times \frakJ \to \frakJ$\,, characterized by the fact that the diagram
\begin{equation*}
\begin{minipage}{5cm}
\begin{xy}
\xymatrix{
\Sp(4) \times \frakJ(4,\HH^{\C})_0 \ar[r]^{F_0} & \frakJ(4, \HH^{\C})_0 \\
\Sp(4) \times \frakJ \ar[r]_{F} \ar[u]^{\id_{\Sp(4)} \times \psi} & \frakJ \ar[u]_{\psi} 
}
\end{xy}
\end{minipage}
\end{equation*}
commutes. It has been shown by Yokota (\cite{Yokota:invol1-1990}, the proof of Theorem~3.4.2) that for any \,$B \in \Sp(4)$\,, \,$F(B) \in \Esix$\, holds,
and \,$F(B)$\, commutes with \,$\lambda\gamma \in \Esix$\,. 
In this way, we obtain a homomorphism of Lie groups \,$F: \Sp(4)\to \Esix$\, with \,$\ker(F) = \{\pm \id\}$\,. 

We now consider the map
$$ f: G_2(\HH^4) \to [\frakJ], \; U \mapsto [\psi^{-1}(Z_U)] \;, $$
where for any \,$U\in G_2(\HH^4)$\, we denote by \,$Z_U \in \frakJ(4,\HH^{\C})_0$\, the linear map characterized by \,$Z_U|U = \tfrac12\,\id_U$\,,
\,$Z_U|U^\perp = -\tfrac12\,\id_{U^\perp}$\,. It is easy to see that \,$f$\, is a well-defined isometric two-fold covering map onto its image with fibers
\,$\{U,U^\perp\}$\, for \,$U \in G_2(\HH^4)$\,, and that \,$f$\, is equivariant, i.e.~that for any \,$U \in G_2(\HH^4)$\, and \,$B\in\Sp(4)$\, we have
$$ F(B)f(U) = f(BU) \; . $$
Because of the latter property and the fact that \,$f(\HH e_1 \oplus \HH e_2) = p_0 \in \EIII$\, holds, 
\,$f$\, maps into \,$\EIII$\,. Moreover, we have \,$p_0 \in F^{\lambda\gamma}$\, and for every \,$B\in\Sp(4)$\,, \,$F(B)\in\Esix$\, commutes with \,$\lambda\gamma$\,,
and therefore \,$f$\, maps into the totally geodesic submanifold \,$F^{\lambda\gamma}$\, of \,$\EIII$\,. 
Because both \,$F^{\lambda\gamma}$\, and \,$G_2(\HH^4)$\, are compact and connected, and they are of the same dimension \,$16$\,, it follows that the isometric immersion
\,$f$\, in fact maps \,$G_2(\HH^4)$\, onto \,$F^{\lambda\gamma}$\,. Because \,$f$\, is a two-fold covering map
with fibers \,$\{U,U^\perp\}$\,, we conclude that \,$F^{\lambda\gamma}$\, is isometric to \,$G_2(\HH^4)/\Z_2$\,. 

\paragraph{The types \,$\boldsymbol{(Q)}$\, and \,$\boldsymbol{(\mathrm{DIII})}$\,.}
The connected components \,$\neq \{p_0\}$\, of the fixed point set of the geodesic symmetry at \,$p_0$\, are the polars of the symmetric space.
(Also see \cite{Chen/Nagano:totges2-1978}, \S 2, especially Theorem~2.8, where the polars are denoted by \,$M_+$\,.)
In the case of \,$\EIII$\,, the polars have also been investigated by Atsuyama in \cite{Atsuyama:projspaces2-1997}, \S 3. 

It is easily seen that the fixed point set of \,$\sigma$\, consists of two connected components besides \,$\{p_0\}$\,, namely
\begin{gather*}
F^\sigma_1 := \left\{ \left. \left[ \begin{smallmatrix} 0 & 0 & 0 \\ 0 & \xi_2 & x_1 \\ 0 & \overline{x_1} & \xi_3 \end{smallmatrix} \right] \in [\frakJ] \;\right|\;
\xi_2\,\xi_3 = |x_1|^2 \right\} \\
\qmq{and} F^\sigma_2 := \left\{ \left. \left[ \begin{smallmatrix} 0 & x_3 & \overline{x_2} \\ \overline{x_3} & 0 & 0 \\ x_2 & 0 & 0 \end{smallmatrix} \right] \in [\frakJ] \;\right|\;
x_2\,x_3 = 0,\; x_2\,\overline{x_2} = 0,\; x_3\,\overline{x_3} = 0 \right\}
\end{gather*}
It turns out that the totally geodesic submanifolds \,$F^\sigma_1$\, and \,$F^\sigma_2$\, are of type \,$(Q)$\, and \,$(\mathrm{DIII})$\,, respectively.

The complex-8-dimensional submanifold \,$F^\sigma_1$\, of the complex projective space \,$[\frakJ]$\, is defined by a single non-degenerate quadratic equation, which is
adapted to the Fubini-Study metric of \,$[\frakJ]$\,. Hence \,$F^\sigma_1$\, is isometric to the complex quadric \,$Q^8$\,.

Furthermore, it has been shown by Atsuyama that the reflective submanifold \,$F^\sigma_2$\, is isometric to \,$\SO(10)/\Ug(5)$\,, see \cite{Atsuyama:projspaces2-1997}, 
the remark after Lemma~3.2 and \cite{Atsuyama:EIIIproj-1985}, the Remark~(2) after Proposition~5.4.

\subsection[Totally geodesic submanifolds in \,$\Sp(2)$\,]{Totally geodesic submanifolds in \,$\boldsymbol{\Sp(2)}$\,}
\label{SSe:EIII:B2}

Our next objective is the classification of the Lie triple systems in the Lie group \,$\Sp(2)$\,, regarded as a Riemannian symmetric space.
We will use this result also in our classification of Lie triple systems of \,$\SO(10)/\Ug(5)$\, in Section~\ref{SSe:EIII:DIII} below.

We will base our classification on the fact that \,$\Sp(2)$\, is a maximal totally geodesic submanifold of \,$G_2(\HH^4)$\, (of type \,$(\Sp_2)$\,
according to the classification in \cite{Klein:2007-tgG2}, Theorem~5.3). 
Because the Lie triple systems of \,$G_2(\HH^4)$\, have been classified in \cite{Klein:2007-tgG2}, we can therefore obtain a classification
of the Lie triple systems by determining which of the Lie triple systems of \,$G_2(\HH^4)$\, are contained in a Lie triple system of type
\,$(\Sp_2)$\,. 

To do so, we will work in the setting of \cite{Klein:2007-tgG2} in the present section. We consider the space
\,$G_2(\HH^4) = \Sp(4)/(\Sp(2)\times\Sp(2))$\,. We let
 \,$\lieg = \liek \oplus \liem$\, be the canonical decomposition associated with this space, i.e.~we have \,$\lieg = \liesp(4)$\,, \,$\liek = \liesp(2)\oplus\liesp(2)
\subset \lieg$\,, and \,$\liem$\, is isomorphic to the tangent space of \,$G_2(\HH^4)$\,. 
We will use the notations of Section~5 of \cite{Klein:2007-tgG2} in the sequel,
especially we use the types of Lie triple systems defined in Theorem~5.3 of \cite{Klein:2007-tgG2} for \,$G_2(\HH^4)$\,, i.e.~for \,$n=2$\,. 
We let \,$\liem_1 \subset \liem$\, be a Lie triple system of \,$\liem$\, of type \,$(\Sp_2)$\,. 

\begin{Theorem}
\label{EIII:SO5:cla}
Exactly the following types of Lie triple systems of \,$\liem$\,,
as defined in Theorem~5.3 of \cite{Klein:2007-tgG2}, have representatives which are contained in \,$\liem_1$\,:
\begin{itemize}
\item \,$(\mathrm{Geo},\vi=t)$\, with \,$t \in [0,\tfrac\pi4]$\,
\item \,$(\Sph,\vi=\arctan(\tfrac13),\ell)$\, with \,$\ell \in \{2,3\}$\,
\item \,$(\PP,\vi=\tfrac\pi4,\tau)$\, with \,$\tau\in \{\;(\R,1)\;,\;(\C,1)\;,\;(\Sph^3)\;,\;(\HH,1)\;\}$\,
\item \,$(\PP\times\PP,\tau_1,\tau_2)$\, with \,$\tau_1,\tau_2 \in \{\;(\R,1)\;,\; (\C,1)\;,\; (\Sph^3)\;\}$\,
\item \,$(\Sph^1\times\Sph^5,\ell)$\, with \,$2\leq \ell\leq 3$\,
\item \,$(Q_3)$\,
\end{itemize}
The maximal Lie triple systems of \,$\liem_1$\, are those which are of the types: \,$(\Sph,\vi=\arctan(\tfrac13),3)$\,, \,$(\PP,\vi=\tfrac\pi4,(\HH,1))$\,, \,$(\PP\times\PP, (\Sph^3), (\Sph^3))$\, and \,$(Q_3)$\,. 
\end{Theorem}

\begin{Remark}
\label{R:EIII:B2:CN}
The maximal totally geodesic submanifolds of \,$\Sp(2)$\, of types \,$(\Sph,\vi=\arctan(\tfrac13),3)$\, and
\,$(Q_3)$\, are missing from \cite{Chen/Nagano:totges2-1978}, Table~VIII. Their isometry type is that of a 3-sphere of radius
\,$\tfrac\pi2\,\sqrt{10}$\,
resp.~of a complex quadric \,$Q^3$\,. The totally geodesic submanifolds of the former type are once again in a ``skew'' position
in the sense that their geodesic diameter is strictly larger than the geodesic diameter \,$\pi$\, of \,$\Sp(2)$\,. 
\end{Remark}

\beweis
It is easily seen that the prototypes for the types listed, as they are given in \cite{Klein:2007-tgG2}, Theorem~5.3, are contained in Lie triple systems
of type \,$(\Sp_2)$\,. 
Therefore, we only need to show that no other types of Lie triple systems of \,$G_2(\HH^4)$\, have representatives which are
contained in \,$\liem_1$\,. 

For this purpose, we let \,$\liem'$\, be a Lie triple system of \,$\liem$\, which is contained in \,$\liem_1$\,. We are to show
that \,$\liem'$\, is of one of the types listed in Theorem~\ref{EIII:SO5:cla}.

If \,$\liem'$\, is of rank \,$2$\,, then for \,$\liem'$\, to be contained in \,$\liem_1$\,, it is necessary that all the roots of \,$\liem'$\, have at
most the multiplicity of the corresponding root in \,$\liem_1$\,. Because the Dynkin diagram of \,$\liem_1$\, is 
$\xymatrix@=.4cm{ \mathop{\bullet}^2 \ar@{=>}[r] & \mathop{\bullet}^2 }$, we see by this argument that \,$\liem'$\, cannot be of one of the types
\,$(\mathrm{G}_2,\tau)$\,, \,$(\PP\times\PP,\tau_1,\tau_2)$\, where either of the \HP-types%
\footnote{See \cite{Klein:2007-tgG2}, Definition~5.1}
\,$\tau_1$\, and \,$\tau_2$\, have dimension \,$\geq 2$\, or width \,$4$\,, or \,$(\Sph^1\times\Sph^5,\ell)$\, where \,$\ell \geq 4$\,.
This already shows that among the types of Lie triple systems of rank \,$2$\, of \,$G_2(\HH^4)$\,, only those which are listed in Theorem~\ref{EIII:SO5:cla}
remain.

If \,$\liem'$\, is of rank \,$1$\,, we note that if \,$\liem'$\, is of type  \,$(\PP,\vi=0,(\C,1))$\, or of type \,$(\PP,\vi=\arctan(\tfrac12),\tau)$\,
with \,$\tau\neq(\R,1)$\,, it cannot be contained in \,$\liem_1$\, because the roots \,$2\lambda_k$\, are not present in \,$\liem_1$\,. 
Because the types \,$(\PP,\vi=0,(\R,1))$\, and \,$(\PP,\vi=\arctan(\tfrac12),(\R,1))$\, are identical to \,$(\mathrm{Geo},\vi=t)$\, with \,$t=0$\,
resp.~with \,$t=\arctan(\tfrac12)$\,, this argument again leaves only the types of rank \,$1$\, which have been listed in the theorem.

For the statements on the maximality: \,$(\PP\times\PP,(\Sph^3),(\Sph^3))$\, and \,$(Q_3)$\, are Lie triple systems of rank \,$2$\,, and therefore can be
contained only in other Lie triple systems of this rank. Because they have the same dimension \,$6$\, and are clearly not isomorphic, neither of
them can be contained in the other, and also for reason of dimension, neither can be contained in \,$(\Sph^1\times\Sph^5,\ell)$\,
with \,$\ell\leq 3$\,. Therefore these two types
are maximal in \,$\liem_1$\,. From a consideration of the root systems it can also be seen that \,$(\Sph,\vi=\arctan(\tfrac13),3)$\, and
\,$(\PP,\vi=\tfrac\pi4,(\HH,1))$\, are maximal.
On the other hand, \,$(\mathrm{Geo},\vi=t)$\,, \,$(\PP,\vi=\tfrac\pi4,(\K,1))$\, with \,$\K\in\{\R,\C\}$\,, \,$(\PP\times\PP,\tau_1,\tau_2)$\, with
\,$\tau_1,\tau_2 \in \{\;(\R,1)\;,\; (\C,1)\;,\; (\Sph^3)\;\}$\, and \,$(\Sph^1\times\Sph^5,\ell)$\, with \,$\ell \leq 3$\, are all contained in
\,$(\PP\times\PP,(\Sph^3),(\Sph^3))$\,, whereas \,$(\Sph,\vi=\arctan(\tfrac13),2)$\, is contained in \,$(\Sph,\vi=\arctan(\tfrac13),3)$\,. Therefore
these types cannot be maximal.
\beweisende

We can obtain the global isometry types of the totally geodesic submanifolds corresponding to the Lie triple systems of \,$\Sp(2)$\, as listed in
Theorem~\ref{EIII:SO5:cla} from the totally geodesic embeddings into \,$G_2(\HH^n)$\, described in \cite{Klein:2007-tgG2}, Section~6. When applying the
information from that paper, one needs to take into account, however, that in the \,$\Sp(2)$\, as totally geodesic submanifold of \,$G_2(\HH^n)$\, 
(with the Riemannian metric considered in that paper) the shortest restricted root has length \,$\sqrt{2}$\,, whereas here we want to view \,$\Sp(2)$\,
with the metric so that the shortest restricted root has length \,$1$\,. Therefore the curvatures of the projective spaces have to be multiplied
with \,$\tfrac12$\,, and the radii of the spheres have to be multiplied with \,$\sqrt{2}$\,, to translate from the situation in \cite{Klein:2007-tgG2}
to the present situation.
In this way, we obtain the following information on the totally geodesic submanifolds
of \,$\Sp(2)_{\mathrm{srr}=1}$\,, where we again use the notations introduced in Section~\ref{SSe:EIII:tgsub}.

\begin{center}
\begin{longtable}{|c|c|c|}
\hline
type of Lie triple system & isometry type & properties \\
\hline
\endhead
\hline
\endfoot
$(\mathrm{Geo},\vi=t)$ & \,$\R$\, or \,$\Sph^1$ & \\
$(\Sph,\vi=\arctan(\tfrac13),\ell)$ & $\Sph^\ell_{r=\sqrt{5}}$\, & \,$\ell=3$\,: maximal \\
$(\PP,\vi=\tfrac\pi4,(\K,1))$ & \,$\KP^1_{\vkap=1/4}$ & \,$\K=\HH$\,: polar, maximal \\
$(\PP,\vi=\tfrac\pi4,(\Sph^3))$ & \,$\Sph^3_{r=1}$\, & \\
$(\PP\times\PP,(\K_1,1),(\K_2,1))$ & \,$\K_1\mathrm{P}^1_{\vkap=1/2} \times \K_2\mathrm{P}^1_{\vkap=1/2}$\, & \\
$(\PP\times\PP,(\K,1),(\Sph^3))$ & \,$\K\mathrm{P}^1_{\vkap=1/2} \times \Sph^3_{r=1/\sqrt{2}}$\, & \,$\K=\R$\,: meridian to \,$(Q_3)$\, \\
$(\PP\times\PP,(\Sph^3),(\Sph^3))$ & \,$\Sph^3_{r=1/\sqrt{2}} \times \Sph^3_{r=1/\sqrt{2}}$\, & meridian to \,$(\PP,\vi=\tfrac\pi4,(\HH,1))$\,, maximal \\
$(\Sph^1\times\Sph^5,\ell)$ & \,$(\Sph^1_{r=1} \times \Sph^\ell_{r=1})/\Z_2$\, & \\
$(Q_3)$ & \,$Q^3_{\mathrm{srr}=1}$\, & polar, maximal \\
\end{longtable}
\end{center}

\subsection[Totally geodesic submanifolds in \,$\SO(10)/\Ug(5)$\,]{Totally geodesic submanifolds in \,$\boldsymbol{\SO(10)/\Ug(5)}$\,}
\label{SSe:EIII:DIII}

We now want to classify the Lie triple systems of \,$\SO(10)/\Ug(5)$\,. Note that this symmetric space occurs as a maximal totally geodesic
submanifold of \,$\EIII$\,. We will use the classification of Lie triple systems of \,$\EIII$\, from Section~\ref{SSe:EIII:lts} to obtain the classification
for \,$\SO(10)/\Ug(5)$\, in an analogous way as we used the classification in \,$G_2(\HH^4)$\, to obtain the classification for \,$\Sp(2)$\,
in the previous section.

Thus we now return to the situation studied in Section~\ref{SSe:EIII:lts}. We consider the Riemannian symmetric space \,$\EIII$\,, and let \,$\lieg=\liek\oplus\liem$\,
be the canonical decomposition of \,$\lieg=\liee_6$\, associated with this space, i.e.~we have \,$\liek = \R \oplus \lieo(10)$\, and \,$\liem$\, is isomorphic
to the tangent space of \,$\EIII$\,. We will use the names for the types of Lie triple systems of \,$\liem$\, as introduced in Theorem~\ref{EIII:EIII:cla}
and Remark~\ref{R:EIII:EIII:extratypes}.

Further, we let \,$\liem_1$\, be a Lie triple system of \,$\liem$\, of type \,$(\mathrm{DIII})$\,, i.e.~the totally geodesic submanifold of \,$\EIII$\, corresponding
to \,$\liem_1$\, is isometric to \,$\SO(10)/\Ug(5)$\,.

\begin{Theorem}
\label{EIII:DIII:cla}
Exactly the following types of Lie triple systems of \,$\EIII$\, have representatives which are contained in \,$\liem_1$\,:
\begin{itemize}
\item \,$(\mathrm{Geo},\vi=t)$\, with \,$t \in [0,\tfrac\pi4]$\,
\item \,$(\PP,\vi=0,(\K,4))$\, with \,$\K\in\{\R,\C\}$\,
\item \,$(\PP,\vi=\tfrac\pi4,(\Sph^k))$\, with \,$k\in\{5,6\}$\,
\item \,$(\PP\times\PP^1,(\K_1,3),\K_2)$\, with \,$\K_1,\K_2 \in \{\R,\C\}$\,
\item \,$(Q,(\mathrm{G1},6))$\,
\item The types \,$(Q,\tau)$\,, where \,$\tau$\, is one of the types of Lie triple systems in the complex quadric as defined in 
\cite{Klein:2007-claQ}, Theorem~4.1, for \,$m=6$\,, i.e.~\,$\tau$\, is one of the following:
\,$(\mathrm{G1},k)$\, with \,$k \leq 5$\,, 
\,$(\mathrm{G2},k_1,k_2)$\, with \,$k_1+k_2 \leq 6$\,, \,$(\mathrm{G3})$\,, \,$(\mathrm{P1},k)$\, with \,$k \leq 6$\,,
\,$(\mathrm{P2})$\,, \,$(\mathrm{A})$\,, \,$(\mathrm{I1},k)$\, with \,$k \leq 3$\,,
and \,$(\mathrm{I2},k)$\, with \,$k \leq 3$\,.
\item \,$(G_2\C^6,(\mathrm{G}_2,(\C,3)))$\,
\item The types \,$(G_2\C^6,\tau)$\,, where \,$\tau$\, is one of the following:
\,$(\PP,\vi=\arctan(\tfrac12),(\K,k))$\, with \,$(\K,k) \in \{(\R,1),(\C,1),(\R,2)\}$\,,
\,$(\mathrm{G}_2,(\K,k))$\, with \,$\K\in\{\R,\C\}$\, and \,$k \leq 3$\,, and
\,$(\PP\times\PP,(\K_1,k_1),(\K_2,k_2))$\, with \,$\K_1,\K_2 \in \{\R,\C\}$\, and \,$k_1+k_2 \leq 3$\,.
\item \,$(G_2\HH^4,(\Sp_2))$\,
\item The types \,$(G_2\HH^4,\tau)$\,, where \,$\tau$\, is one of the following:
\,$(\Sph,\vi=\arctan(\tfrac13),3)$\,, 
\,$(\PP,\vi=\tfrac\pi4,(\K,1))$\, with \,$\K\in\{\R,\C,\HH\}$\,, and
\,$(\Sph^1\times\Sph^5,3)$\,.
\end{itemize}
The maximal Lie triple systems of \,$\liem_1$\, are those of the types: \,$(\PP,\vi=0,(\C,4))$\,, \,$(\PP\times\PP^1,(\C,3),\C)$\,,
\,$(Q,(\mathrm{G1},6))$\,, \,$(G_2\C^6,(\C,3))$\, and \,$(G_2\HH^4,(\Sp_2))$\,. 
\end{Theorem}

\begin{Remark}
\label{R:EIII:DIII:CN}
Chen and Nagano correctly list the local isometry type of all the maximal totally geodesic submanifolds of \,$\SO(10)/\Ug(5)$\, in their
Table~VIII of \cite{Chen/Nagano:totges2-1978}. However the isometry types of the types \,$(Q,(\mathrm{G1},6))$\, 
resp.~\,$(G_2\HH^4,(\Sp_2)))$\, are \,$Q^6$\, resp.~\,$\SO(5)$\, (where Chen/Nagano claim \,$G_2(\R^8) \cong Q^6/\Z_2$\, 
and \,$\Sp(2)\cong \Spin(5)$\, respectively). Moreover, it should be mentioned that also \,$\SO(10)/\Ug(5)$\, contains ``skew'' totally geodesic submanifolds,
namely the totally geodesic submanifolds of the types \,$(Q,(A))$\, and \,$(G_2\HH^4,(\Sph,\vi=\arctan(\tfrac13),3))$\,, which are isometric to
a 2-sphere resp.~a 3-sphere of radius \,$\sqrt5$\,, so that their geodesic diameter \,$\sqrt{5}\,\pi$\, is strictly larger than the geodesic diameter
\,$\pi$\, of \,$\SO(10)/\Ug(5)$\,. They are not maximal in \,$\SO(10)/\Ug(5)$\,; their presence can not be inferred from Table~VIII of \cite{Chen/Nagano:totges2-1978}
because of the missing entries for the spaces \,$G_2^+(\R^5)$\, and \,$\Sp(2)$\,. 
\end{Remark}

{
\footnotesize
\emph{Proof of Theorem~\ref{EIII:DIII:cla}.}
For the maximal ones among the types listed, corresponding totally geodesic embeddings into \,$\SO(10)/\Ug(5)$\, are described 
below, so we know that these types, and therefore also all the other types listed, have representatives contained in \,$\liem_1$\,.
Therefore, we only need to show that no other types of Lie triple systems of \,$\EIII$\, have representatives which are
contained in \,$\liem_1$\,. 

For this purpose, we let \,$\liem'$\, be a Lie triple system of \,$\liem$\, which is contained in \,$\liem_1$\,. We are to show
that \,$\liem'$\, is of one of the types listed in Theorem~\ref{EIII:DIII:cla}.

If \,$\liem'$\, is of rank \,$2$\,, all the roots of \,$\liem'$\,
have at most the multiplicity of the corresponding root in \,$\liem_1$\,. Because the Dynkin diagram of \,$\liem_1$\, is 
$\xymatrix@=.4cm{ \mathop{\bullet}^4 \ar@{<=>}[r] & \mathop{\bullet\hspace{-.29cm}\bigcirc}^{4[1]} }$, we see by this argument that \,$\liem'$\,
cannot be one of the types \,$(\PP\times\PP^1,(\K_1,\ell),\K_2)$\, with \,$\ell \geq 4$\,, \,$(Q)$\,, \,$(G_2\C^6)$\, and \,$(G_2\HH^4)$\,. Moreover, 
we note that the intersection of \,$\liem_1$\, with a Lie triple system of type \,$(Q)$\, is of type \,$(Q,(\mathrm{G1},6))$\, (corresponding to \,$Q^6
\subset Q^8 \subset \EIII$\,), with a Lie triple system of type \,$(G_2\C^6)$\, is of type \,$(G_2\C^6,(G_2,(\C,3)))$\, (corresponding to
\,$G_2(\C^5) \subset G_2(\C^6) \subset \EIII$\,), and with a Lie triple system of type \,$G_2(\HH^4)$\, is of type \,$(G_2\HH^4,(\Sp_2))$\,
(corresponding to \,$\SO(5) \subset G_2(\HH^4)/\Z_2 \subset \EIII$\,). Therefrom everything about the rank 2 Lie triple systems follows.

For the spaces of rank \,$1$\, a similar consideration of the multiplicities of the roots shows that \,$\liem'$\, is of one
of the types listed in the theorem.
\strut\hfill $\Box$

}

Because \,$\SO(10)/\Ug(5)$\, is a totally geodesic submanifold of \,$\EIII$\,, the isometry types of the totally geodesic submanifolds in \,$\SO(10)/\Ug(5)$\,
corresponding to the various types of Lie triple systems are the same as the isometry types of the totally geodesic submanifolds in \,$\EIII$\,
of those types, which were described in Section~\ref{SSe:EIII:tgsub}. In particular, the isometry types of the maximal totally geodesic
submanifolds of \,$\SO(10)/\Ug(5)$\,, and some of their properties, are as follows:

\begin{center}
\begin{tabular}{|c|c|c|}
\hline
type & isometry type & properties \\
\hline
$(\PP,\vi=0,(\C,4))$ & \,$\CP^4_{\vkap=1}$\, & polar \\
\,$(\PP\times\PP^1,(\C,3),\C)$\, & $\CP^3_{\vkap=1} \times \CP^1_{\vkap=1}$ & meridian for \,$(G_2\C^6,(\mathrm{G}_2,(\C,3)))$\, \\
\,$(Q,(\mathrm{G1},6))$\, & $Q^6_{\mathrm{srr}=\sqrt{2}}$ & meridian for \,$(\PP,\vi=0,(\C,4))$\, \\
\,$(G_2\C^6,(\mathrm{G}_2,(\C,3)))$\, & $G_2(\C^5)_{\mathrm{srr}=1}$ & polar \\
\,$(G_2\HH^4,(\Sp_2))$\, & $\SO(5)_{\mathrm{srr}=1}$ & reflective \\
\hline
\end{tabular}
\end{center}
To elucidate the position of the maximal totally geodesic submanifolds, we describe totally geodesic embeddings for these types:

\paragraph{The types \,$\boldsymbol{(\PP,\vi=0,(\C,4))}$\, and \,$\boldsymbol{(G_2\C^6,(\mathrm{G}_2,(\C,3)))}$\,.}
The totally geodesic submanifolds of these types are the polars in \,$\SO(10)/\Ug(5)$\,, and can therefore be obtained as \,$\Ug(5)$-orbits through
points of \,$\SO(10)/\Ug(5)$\, which are antipodal to the origin point \,$p_0 := \Ug(5) \in \SO(10)/\Ug(5)$\, in this space.

For an explicit construction, we consider both \,$\Ug(5)$\, and \,$\SO(10)$\, acting on \,$\C^5$\,; in the latter case the action is only \,$\R$-linear
on \,$\C^5 \cong \R^{10}$\,. We fix a real form \,$V$\, of \,$\C^5$\, (i.e.~a 5-dimensional real linear subspace \,$V \subset \C^5$\, so that \,$i\cdot V$\,
is orthogonal to \,$V$\, with respect to the real inner product on \,$\C^5$\,). Then we can describe \,$\Ug(5)$\, as a subgroup of \,$\SO(10)$\, by
$$ \Ug(5) = \left\{ \, g \in \SO(10) \, \left| \, g = \left( \begin{smallmatrix} A & -B \\ B & A \end{smallmatrix} \right),\; A,B \in M(5\times 5,\R) \right. \right\} \;, $$
where the matrix expression is to be read as a block matrix with respect to the splitting \,$\C^5 = V \operp i\cdot V$\,. In the same way, we can describe
the involutive automorphism describing the symmetric structure of \,$\SO(10)/\Ug(5)$\,:
$$ \sigma : \SO(10) \to \SO(10),\; \left( \begin{smallmatrix} A & C \\ B & D \end{smallmatrix} \right) \;\mapsto\; 
\left( \begin{smallmatrix} D & -B \\ -C & A \end{smallmatrix} \right) \; . $$
Via the linearization of \,$\sigma$\,, we obtain the space \,$\liem$\, in the splitting \,$\lieo(10) = \lieu(5) \oplus \liem$\, induced by the symmetric
structure of \,$\SO(10)/\Ug(5)$\,:
$$ \liem = \left\{ \left. \left( \begin{smallmatrix} A & B \\ B & -A \end{smallmatrix} \right) \right| A,B \in \lieo(5) \right\} \; . $$

We now fix a \,$2k$-dimensional real linear subspace \,$W \subset V$\, (with \,$k\in\{1,2\}$\,) and a ``partial complex structure with respect to \,$W$\,'',
i.e.~a skew-adjoint transformation \,$J: V \to V$\, with \,$J^3 = -J$\, and \,$J(V)=W$\,. Then we have \,$X := \left( \begin{smallmatrix} J & 0 \\ 0 & -J
\end{smallmatrix} \right) \in \liem$\,, and therefore \,$\gamma: \R \to \SO(10)/\Ug(5),\;t \mapsto \exp(tX)\cdot p_0$\, is a geodesic of \,$\SO(10)/\Ug(5)$\,.
For \,$t\in\R$\, and \,$w\in W$\, we have
$$ \exp(tX)w = \cos(t)\,w + \sin(t)\,Jw \qmq{and} \exp(tX)iw = \cos(t)\,iw - \sin(t)\,iJw $$
as well as \,$\exp(tX)w' = w'$\, for any \,$w' \in (W \oplus iW)^{\perp}$\,. We have \,$\gamma(t)=p_0$\, if and only if \,$\exp(tX)\in\Ug(5)$\,;
from the above description it follows that this is the case if and only if \,$\sin(t)=0$\, holds, i.e.~if we have \,$t\in\pi\,\Z$\,. Hence the geodesic \,$\gamma$\,
is periodic with period \,$\pi$\,, and therefore \,$p_1 := \gamma(\tfrac\pi2)$\, is an antipodal point of \,$p_0$\, in \,$\SO(10)/\Ug(5)$\,. By general results
(see \cite{Chen/Nagano:totges2-1978}, Lemma~2.1), it is known that the polar \,$M := \Ug(5)\cdot p_1$\, is a totally geodesic submanifold of \,$\SO(10)/\Ug(5)$\,. 

To determine the isometry type of the totally geodesic submanifold \,$M$\,, we calculate the isotropy group of the action of \,$\Ug(5)$\,
at \,$p_1$\,: We have \,$p_1 = S\cdot \Ug(5)$\, with \,$ S := \exp(\tfrac\pi2\,X) \in \SO(10)$\,; from the explicit description of \,$X$\, we 
obtain the explicit description
\begin{equation}
\label{eq:EIII:DIII:polar-S}
S|W = J|W\,, \quad S|iW = -J|iW\,,\quad S|(W\oplus iW)^\perp = \id_{(W \oplus iW)^\perp} 
\end{equation}
of \,$S$\,. Therefore we have for \,$g \in \Ug(5)$\,:
$$ g\cdot p_1 = p_1 \;\Longleftrightarrow\; g\cdot S \cdot \Ug(5) = S \cdot \Ug(5) \;\Longleftrightarrow\; S^{-1} \, g \, S \in \Ug(5)
\;\Longleftrightarrow\; g(W\oplus iW) = W\oplus iW \;, $$
where the last equivalence follows from Equations~\eqref{eq:EIII:DIII:polar-S}. Therefore the isotropy group of the action of \,$\Ug(5)$\, at \,$p_1$\,
is isomorphic to \,$\Ug(W\oplus iW) \times \Ug((W\oplus iW)^\perp) \cong \Ug(2k) \times \Ug(5-2k)$\,. It follows that the totally geodesic submanifold \,$M$\, of 
\,$\SO(10)/\Ug(5)$\, is isometric to \,$\Ug(5)/ (\Ug(2k)\times\Ug(5-2k))$\,.

In the case \,$k=1$\,, \,$M$\, is thus isometric to \,$\Ug(5)/(\Ug(2) \times \Ug(3))\cong G_2(\C^5)$\,; this totally geodesic submanifold turns out to be of
type \,$(G_2\C^6,(\C,3))$\,. 

In the case \,$k=2$\,, \,$M$\, is isometric to \,$\Ug(5) / (\Ug(4) \times \Ug(1)) \cong \CP^4$\,; this totally geodesic submanifold is of type
\,$(\PP,\vi=0,(\C,4))$\,. 

\paragraph{The types \,$\boldsymbol{(\PP\times\PP^1,(\C,3),\C)}$\, and \,$\boldsymbol{(Q,(\mathrm{G1},6))}$\,.}
These types are the meridians corresponding to the polars of type \,$(G_2\C^6,(\C,3))$\, and \,$(\PP,\vi=0,(\C,4))$\,, respectively. 
This means that they are totally geodesic submanifolds which intersect the corresponding polar orthogonally and transversally in one point.

However, in the present situation there is an easier way to describe the totally geodesic submanifolds of these types. Note that there are canonical embeddings
\,$\SO(4) \times \SO(6) \subset \SO(10)$\, and \,$\SO(8) \subset \SO(10)$\, which are compatible with the symmetric structure of \,$\SO(10)/\Ug(5)$\,.
In this way, we get totally geodesic embeddings of \,$(\SO(4)/\Ug(2)) \times (\SO(6)/\Ug(3)) \cong \CP^1 \times \CP^3$\, and of \,$\SO(8)/\Ug(4) \cong Q^6$\,
into \,$\SO(10)/\Ug(5)$\,; they are of type \,${(\PP\times\PP^1,(\C,3),\C)}$\, and \,${(Q,(\mathrm{G1},6))}$\,, respectively.

\paragraph{The type \,$\boldsymbol{(G_2\HH^4,(\Sp_2))}$\,.} Consider the map
$$ \Phi: \SO(5) \to \SO(10), \; B \mapsto \left( \begin{smallmatrix} B & 0 \\ 0 & B^{-1} \end{smallmatrix} \right) \; . $$
For \,$B\in\SO(5)$\, we have \,$\Phi(B) \in \Ug(5) \Longleftrightarrow B = \id$\,, and therefore \,$\Phi$\, induces an embedding
\,$\underline{\Phi}: \SO(5) \to \SO(10)/\Ug(5)$\,. Its linearization maps \,$X \in \lieo(5)$\, onto \,$\left( \begin{smallmatrix} X & 0 \\ 0 & -X \end{smallmatrix}
\right) \in \liem$\,, and therefore \,$\underline{\Phi}$\, is totally geodesic. It turns out to be of type \,$(G_2\HH^4,(\Sp_2))$\,. 

\section[The symmetric spaces \,$E_6/F_4$\,, \,$\SU(6)/\Sp(3)$\,, \,$\SU(3)$\, and \,$\SU(3)/\SO(3)$\,]{The symmetric spaces \,$\boldsymbol{E_6/F_4}$\,, \,$\boldsymbol{\SU(6)/\Sp(3)}$\,, \,$\boldsymbol{\SU(3)}$\, and \,$\boldsymbol{\SU(3)/\SO(3)}$\,}
\label{Se:EIV}

\subsection[The geometry of \,$E_6/F_4$\,]{The geometry of \,$\boldsymbol{E_6/F_4}$\,}
\label{SSe:EIV:geometry}

In this section we will study the Riemannian symmetric space \,$\EIV := \Esix/\Ffour$\,, which has the Satake diagram
\begin{center}
\begin{minipage}{5cm}
\medskip
\xymatrix@=.4cm{
{\displaystyle \mathop{\circ}^{1}} \ar@{-}[r]<-.8ex> 
        & \displaystyle \mathop{\bullet}^{3} \ar@{-}[r]<-.8ex>
        & \displaystyle \mathop{\bullet}^{4} \ar@{-}[r]<-.8ex> \ar@{-}[d]
        & \displaystyle \mathop{\bullet}^{5} \ar@{-}[r]<-.8ex>
        & {\displaystyle \mathop{\circ}^{6}} \\
&& {\displaystyle \mathop{\bullet}_{2}} &&
}
\medskip
\end{minipage}
\end{center}
\,$\EIV$\, does not have an invariant Hermitian structure.

We consider the Lie algebra \,$\lieg:=\liee_6$\, of the transvection group \,$\Esix$\, of \,$\EIV$\,
and the splitting \,$\lieg = \liek \oplus \liem$\, induced by the symmetric structure
of \,$\EIV$\,. Herein \,$\liek = \lief_4$\, is the Lie algebra of the isotropy group of \,$\EIV$\,, and \,$\liem$\, is isomorphic to the tangent space of \,$\EIV$\,
in the origin.
The \,$E_6$-invariant Riemannian metric on \,$\EIV$\, induces an \,$\Ad(F_4)$-invariant Riemannian
metric on \,$\liem$\,. As was explained in Section~\ref{Se:generallts}, this metric is only unique up to a factor; we choose the factor in such a way
that the restricted roots of \,$\EIV$\, (see below) have the length \,$1$\,.

\paragraph{The root space decomposition.}
Let \,$\liet$\, be a Cartan subalgebra of \,$\lieg$\, which is maximally non-compact, i.e.~it is such that the flat subspace
\,$\liea := \liet \cap \liem$\, of \,$\liem$\, is of the maximal dimension \,$2$\,, and hence a Cartan subalgebra of \,$\liem$\,. 
Then we consider the root system \,$\Delta^\lieg \subset \liet^*$\, of \,$\lieg$\, with respect to \,$\liet$\,,
as well as the restricted root system \,$\Delta \subset \liea^*$\, of the symmetric space \,$\EIV$\, with respect to \,$\liea$\,. 
\,$\EIV$\, has the restricted Dynkin diagram  $\xymatrix@=.4cm{ \mathop{\bullet}^8 \ar@{-}[r] & \mathop{\bullet}^{8} }$,
in other words: its restricted root system \,$\Delta$\, is of type \,$A_2$\,, i.e.~we have
\,$\Delta = \{\pm\lambda_1, \pm\lambda_2, \pm\lambda_3\}$\,, where
\,$(\lambda_1,\lambda_2)$\, is a system of simple roots of \,$\Delta$\,, these two roots are at an angle of \,$\tfrac{2}{3}\,\pi$\, and have the same length,
and we have \,$\lambda_3 = \lambda_1+\lambda_2$\,. All roots in \,$\Delta$\, have the multiplicity \,$8$\,, and \,$\Delta$\, has the following
graphical representation:
\vspace{-0.5cm}
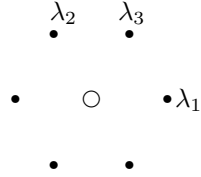
\begin{center}
\begin{minipage}{5cm}
\begin{center}
\strut \\[.3cm]
\setlength{\unitlength}{1cm}
\begin{picture}(2,2)
\put(1,1){\circle{0.2}}
\put(2,1){\circle*{0.1}}
\put(1.5,1.866){\circle*{0.1}}
\put(0.5,1.866){\circle*{0.1}}
\put(1.5,0.134){\circle*{0.1}}
\put(0.5,0.134){\circle*{0.1}}
\put(0,1){\circle*{0.1}}

\put(2.1,0.90){{\small $\lambda_1$}}
\put(1.35,2.05){{\small $\lambda_3$}}
\put(0.45,2.05){{\small $\lambda_2$}}
\end{picture}
\strut \\[1cm]
\end{center}
\end{minipage}
\end{center}
\vspace{.2cm}

To be able to apply the results from \cite{Klein:2007-Satake} and the corresponding computer package for the calculation of the curvature tensor of \,$\EIV$\,,
we again need to describe the relationship between the restricted roots of the symmetric space \,$\EIV$\, and the (non-restricted) roots of the Lie algebra \,$\liee_6$\,.
For this purpose, we again denote the positive roots of \,$\liee_6$\, by \,$\alpha_1,\dotsc,\alpha_{36}$\, in the way described in Section~\ref{SSe:EIII:geometry}.
To find out which restricted root of \,$\EIV$\, corresponds to each root of \,$\liee_6$\,,
we tabulate the orbits of the action of \,$\sigma$\, on the root system \,$\Delta^\lieg$\,, and the restricted root of \,$\EIV$\, corresponding to each orbit
(compare Section~4 of \cite{Klein:2007-Satake}):

\medskip
{\footnotesize
\begin{center}
\begin{tabular}{|c||c|c|c|c|}
\hline
orbit & $\{\alpha_1,-\alpha_{30}\}$ & $\{\alpha_7,-\alpha_{27}\}$ & $\{\alpha_{12},-\alpha_{22}\}$ & $\{\alpha_{17},-\alpha_{18}\}$ \\
\hline
corresp.~restr.~root & $\lambda_1$ & $\lambda_1$ & $\lambda_1$ & $\lambda_1$ \\
\hline
\end{tabular}
\end{center}

\medskip

\begin{center}
\begin{tabular}{|c||c|c|c|c|}
\hline
orbit & $\{\alpha_6,-\alpha_{31}\}$ & $\{\alpha_{11},-\alpha_{28}\}$ & $\{\alpha_{16},-\alpha_{24}\}$ & $\{\alpha_{20},-\alpha_{21}\}$ \\
\hline
corresp.~restr.~root & $\lambda_2$ & $\lambda_2$ & $\lambda_2$ & $\lambda_2$ \\
\hline
\end{tabular}
\end{center}

\medskip

\begin{center}
\begin{tabular}{|c||c|c|c|c|}
\hline
orbit & $\{\alpha_{23},-\alpha_{36}\}$ & $\{\alpha_{26},-\alpha_{35}\}$ & $\{\alpha_{29},-\alpha_{34}\}$ & $\{\alpha_{32},-\alpha_{33}\}$ \\
\hline
corresp.~restr.~root & $\lambda_3$ & $\lambda_3$ & $\lambda_3$ & $\lambda_3$ \\
\hline
\end{tabular}
\end{center}

}

\medskip

Moreover, we have \,$\sigma(\alpha_k)=\alpha_k$\, for \,$k \in \{2,3,4,5,8,9,10,13,14,15,19,25\}$\,.

Using the notations of \cite{Klein:2007-Satake}, Proposition~5.2(a) we now put for \,$c_1,\dotsc,c_4 \in \C$\, and \,$t\in\R$\,, and where \,$A$\, denotes one
of the letters \,$K$\, and \,$M$\,:
\begin{align*}
A_{\lambda_1}(c_1,c_2,c_3,c_4) & := A_{\alpha_{1}}(c_1) + A_{\alpha_{7}}(c_2) + A_{\alpha_{12}}(c_3) + A_{\alpha_{17}}(c_4) \; , \\
A_{\lambda_2}(c_1,c_2,c_3,c_4) & := A_{\alpha_{6}}(c_1) + A_{\alpha_{11}}(c_2) + A_{\alpha_{16}}(c_3) + A_{\alpha_{20}}(c_4) \; , \\
A_{\lambda_3}(c_1,c_2,c_3,c_4) & := A_{\alpha_{23}}(c_1) + A_{\alpha_{26}}(c_2) + A_{\alpha_{29}}(c_3) + A_{\alpha_{32}}(c_4) \; . 
\end{align*}
Then we have \,$\liem_{\lambda_k} = M_{\lambda_k}(\C,\C,\C,\C)$\, for \,$k\in\{1,2,3\}$\,.

\paragraph{The action of the isotropy group.} We now look at the isotropy action of \,$\EIV$\,. Regarding it, we again use the notations introduced at the end
of Section~\ref{Se:generallts}, in particular we have the continuous function \,$\vi: \liem \setminus \{0\} \to [0,\tfrac\pi3]$\, parametrizing the
orbits of the isotropy action.
For the elements of the closure \,$\overline{\liec}$\, of the positive Weyl chamber \,$\liec := \Menge{v\in\liea}{\lambda_1(v)\geq 0,
\lambda_2(v)\geq 0}$\, we again explicitly describe the relation to their isotropy angle: \,$(\tfrac{\lambda_1^\sharp+\lambda_3^\sharp}{\sqrt{3}},
\lambda_2^\sharp)$\, is an orthonormal basis of \,$\liea$\, 
so that with \,$v_t := \cos(t)\,\tfrac{\lambda_1^\sharp+\lambda_3^\sharp}{\sqrt{3}} + \sin(t)\lambda_2^\sharp$\, we have
\begin{equation}
\label{eq:EIV:isotropy:liec}
\overline{\liec} = \Menge{s\cdot v_t}{t\in[0,\tfrac\pi3],s\in\R_{\geq 0}}\;, 
\end{equation}
and because the Weyl chamber \,$\liec$\, is bordered
by the two vectors \,$v_0 = \tfrac{\lambda_1^\sharp+\lambda_3^\sharp}{\sqrt{3}}$\, with \,$\vi(v_0)=0$\, 
and \,$v_{\pi/3}=\tfrac{\lambda_1^\sharp+\lambda_3^\sharp}{\sqrt{3}}$\, with \,$\vi(v_{\pi/3})=\tfrac\pi3$\,, 
we have
\begin{equation}
\label{eq:EIV:isotropy:vivt}
\vi(s\cdot v_t)=t \qmq{for all} t\in[0,\tfrac\pi3]\,, s\in\R_+ \; . 
\end{equation}

The isotropy action of \,$K=F_4$\, on \,$\liem$\, corresponds to the irreducible 26-dimensional representation of \,$F_4$\, (see \cite{Adams:1996},
Lemma~14.4(i), p.~95). It can be described as an action of \,$F_4$\, on the 26-dimensional
space \,$\frakJ(3,\OO)_0 := \Menge{X\in M(3\times 3,\OO)}{X^*=X, \tr(X)=0}$\, of trace-free, Hermitian \,$(3\times 3)$-matrices over the division algebra
of octonions \,$\OO$\,, for the details see \cite{Adams:1996}, Chapter~16. Under the identification of \,$\liem$\, with \,$\frakJ(3,\OO)_0$\, induced thereby,
the Cartan subalgebra \,$\liea$\, corresponds to the space of trace-free diagonal matrices, and the three root spaces \,$\liem_{\lambda_k}$\, (\,$k=1,2,3$\,) correspond
to the subspaces \,$\frakJ_k := \Menge{X=(x_{ij})\in \frakJ(3,\OO)_0}{x_{11}=x_{22}=x_{33}=x_{k\ell}=x_{km}=0}$\, of \,$\frakJ(3,\OO)_0$\,, where \,$\ell$\, and \,$m$\,
are the two members of \,$\{1,2,3\}\setminus\{k\}$\,. 

The subgroup \,$K_0$\, of \,$F_4$\, with Lie algebra \,$\liek^\liea := \Menge{X\in\liek}{[X,\liea]=0}$\, consists of all those \,$g\in F_4$\, 
which leave the \,$\frakJ_k$\, invariant, and 
is therefore isomorphic to \,$\Spin(8)$\, (see \cite{Adams:1996}, Theorem~16.7(iii)). \,$K_0$\, acts 
on the three spaces \,$\frakJ_k$\, as the three irreducible 8-dimensional representations of \,$\Spin(8)$\,: the vector representation, and the two spin representations;
these representations are ``intertwined'' by the triality automorphism of \,$\Spin(8)$\,.

\begin{Prop}
\label{P:EIV:isotropy}
We regard \,$\R^8$\, as the real linear space underlying \,$\C^4$\, and for \,$k\in\{1,2,3\}$\, we consider the linear isometry
$$ \vi: \R^8 \to \liem_{\lambda_k},\; (c_1,c_2,c_3,c_4) \mapsto M_{\lambda_k}(c_1,c_2,c_3,c_4) \; . $$
Then there exists an isomorphism of Lie groups \,$\Phi: \Spin(8) \to K_0$\, so that the following diagram commutes:
\begin{equation*}
\begin{minipage}{5cm}
\begin{xy}
\xymatrix{
\Spin(8)\times \R^8 \ar[r]^{\Phi\times\vi} \ar[d] & K_0 \times \liem_{\lambda_1} \ar[d]^{\Ad} \\
\R^8 \ar[r]_{\vi} & \liem_{\lambda_1} \;,
}
\end{xy}
\end{minipage}
\end{equation*}
where the left vertical arrow represents the canonical action of \,$\Spin(8)$\, on \,$\R^8$\,.

If we fix \,$v_1 \in \liem_{\lambda_1} \setminus \{0\}$\,, then the Lie subgroup \,$\Spin' := \Menge{B\in \Spin(8)}{B(\vi^{-1}v_1) = \vi^{-1}v_1}$\, of \,$\Spin(8)$\,
is isomorphic to \,$\Spin(7)$\,, and the subgroup \,$K_0' := \Phi(\Spin')$\,, which is isomorphic to \,$\Spin(7)$\,, acts transitively on \,$\liem_{\lambda_2}$\,.

If we now also fix \,$v_2 \in \liem_{\lambda_2} \setminus \{0\}$\,, then the Lie subgroup \,$\Spin'' := \Menge{B\in \Spin'}{B(\vi^{-1}v_2) = \vi^{-1}v_2}$\, of 
\,$\Spin'$\, is isomorphic to the exceptional Lie group \,$G_2$\,, and hence \,$\Phi(\Spin'')$\, is also isomorphic to \,$G_2$\,.

These statements are also true for an arbitrary permutation of the indices \,$1$\,, \,$2$\,, \,$3$\, of the root spaces \,$\liem_{\lambda_k}$\,. 
\end{Prop}

\beweis
Most statements follow from the preceding discussion of the isotropy action. For the transitivity statements, see \cite{Adams:1996}, Lemma~14.13, p.~100.
\beweisende

\subsection[Lie triple systems in \,$E_6/F_4$\,.]{Lie triple systems in \,$\boldsymbol{E_6/F_4}$\,.}
\label{SSe:EIV:lts}

\begin{Theorem}
\label{T:EIV:cla} 
The linear subspaces \,$\liem' \subset \liem$\, given in the following are Lie triple systems, 
and every Lie triple system \,$\{0\} \neq \liem' \subsetneq \liem$\, is congruent under the isotropy action to one of them. 

\begin{itemize}
\item \,$\boldsymbol{(\mathrm{Geo},\vi=t)}$\, with \,$t\in[0,\tfrac\pi3]$\,. \\
\,$\liem' = \R(\cos(t)\,\tfrac{\lambda_1^\sharp+\lambda_3^\sharp}{\sqrt{3}} + \sin(t)\,\lambda_2^\sharp)$\, 
(compare Equation~\eqref{eq:EIV:isotropy:liec}).
\item \,$\boldsymbol{(\Sph,\vi=\tfrac\pi6,\ell)}$\, with \,$\ell \leq 9$\,. \\
\,$\liem'$\, is an \,$\ell$-dimensional linear subspace of \,$\R\lambda_1^\sharp \oplus \liem_{\lambda_1}$\,. 
\item \,$\boldsymbol{(\PP,\vi=\tfrac\pi6,(\K,\ell))}$\, with \,$\K\in\{\R,\C,\HH\}$\, and \,$\ell\in\{2,3\}$\,, or with \,$(\K,\ell)=(\OO,2)$\,. \\
We define the following vectors:
\begin{center}
\begin{tabular}{cc}
$v_0 := M_{\lambda_1}(1,0,0,0)+ M_{\lambda_2}(1,0,0,0)$ & $v_1 := M_{\lambda_1}(i,0,0,0) + M_{\lambda_2}(-i,0,0,0)$ \\
$v_0^C := M_{\lambda_1}(0,0,0,i)+M_{\lambda_2}(0,0,0,-i)$ & $v_1^C := M_{\lambda_1}(0,0,0,1) + M_{\lambda_2}(0,0,0,1)$ \\
$v_0^H := M_{\lambda_1}(0,0,i,0)+M_{\lambda_2}(0,0,-i,0)$ & $v_1^H := M_{\lambda_1}(0,0,1,0) + M_{\lambda_2}(0,0,1,0)$ \\
$v_0^{CH} := M_{\lambda_1}(0,1,0,0)+M_{\lambda_2}(0,1,0,0)$ & $v_1^{CH} := M_{\lambda_1}(0,-i,0,0) + M_{\lambda_2}(0,i,0,0)$ \\
$v_0^O := M_{\lambda_1}(i,0,0,0)+M_{\lambda_2}(i,0,0,0)$ & \\
$v_0^{CO} := M_{\lambda_1}(0,0,0,1)+M_{\lambda_2}(0,0,0,-1)$ & \\
$v_0^{HO} := M_{\lambda_1}(0,0,1,0)+M_{\lambda_2}(0,0,-1,0)$ & \\
$v_0^{CHO} := M_{\lambda_1}(0,-i,0,0)+M_{\lambda_2}(0,-i,0,0)$ & \\
$H := \lambda_3^\sharp$ & $w_4 := M_{\lambda_3}(0,0,0,1)$ \\
$w_1 := M_{\lambda_3}(1,0,0,0)$ & $w_5 := M_{\lambda_3}(-i,0,0,0)$ \\
$w_2 := M_{\lambda_3}(0,1,0,0)$ & $w_6 := M_{\lambda_3}(0,i,0,0)$ \\
$w_3 := M_{\lambda_3}(0,0,i,0)$ & $w_7 := M_{\lambda_3}(0,0,1,0)$ \\
\end{tabular}
\end{center}
Then \,$\liem'$\, is spanned by the following vectors, in dependence of \,$(\K,\ell)$\,: \\
For \,$(\K,\ell)=(\R,2)$\,: \,$H, v_0$\, \\
For \,$(\K,\ell)=(\R,3)$\,: \,$H, v_0, v_1$\, \\
For \,$(\K,\ell)=(\C,2)$\,: \,$H, v_0, v_0^C, w_1$\, \\
For \,$(\K,\ell)=(\C,3)$\,: \,$H, v_0, v_0^C, v_1, v_1^C, w_1$\, \\
For \,$(\K,\ell)=(\HH,2)$\,: \,$H, v_0, v_0^C, v_0^H, v_0^{CH}, w_1, w_2, w_3$\, \\
For \,$(\K,\ell)=(\HH,3)$\,: \,$H, v_0, v_0^C, v_0^H, v_0^{CH}, v_1, v_1^C, v_1^H, v_1^{CH}, w_1, w_2, w_3$\, \\
For \,$(\K,\ell)=(\OO,2)$\,: \,$H, v_0, v_0^C, v_0^H, v_0^{CH}, v_0^O, v_0^{CO}, v_0^{HO}, v_0^{CHO}, w_1, w_2, w_3, w_4, w_5, w_6, w_7$\, 
\item \,$\boldsymbol{(\mathrm{AI})}$\,. \\
\,$\liem' = \liea \oplus M_{\lambda_1}(\R,0,0,0) \oplus M_{\lambda_2}(\R,0,0,0) \oplus M_{\lambda_3}(0,0,0,i\R)$\,.
\item \,$\boldsymbol{(\mathrm{A_2})}$\,. \\
\,$\liem' = \liea \oplus M_{\lambda_1}(\C,0,0,0) \oplus M_{\lambda_2}(\C,0,0,0) \oplus M_{\lambda_3}(0,0,0,\C)$\,.
\item \,$\boldsymbol{(\mathrm{AII})}$\,. \\
\,$\liem' = \liea \oplus M_{\lambda_1}(\C,\C,0,0) \oplus M_{\lambda_2}(\C,\C,0,0) \oplus M_{\lambda_3}(0,0,\C,\C)$\,.
\item \,$\boldsymbol{(\Sph \times \Sph^1,\ell)}$\, with \,$\ell \leq 9$\,. \\
\,$\liem' = \liea \oplus \liem_{\lambda_1}'$\, with an \,$(\ell-1)$-dimensional linear subspace \,$\liem_{\lambda_1}' \subset \liem_{\lambda_1}$\,. 
\end{itemize}
We call the full name \,$(\mathrm{Geo},\vi=t)$\, etc.~given in the above table the \emph{type} of the Lie triple systems which are congruent under the
adjoint action to the space given in that entry. Then every Lie triple system of \,$\liem$\, is of exactly one type.

The Lie triple systems \,$\liem'$\, of the various types have the properties given in the following table. 
The column ``isometry type'' again gives the isometry type of the totally geodesic submanifolds corresponding to the Lie triple systems of the respective type
in abbreviated form, for the details see Section~\ref{SSe:EIV:tgsub}.
\begin{center}
\begin{tabular}{|c|c|c|c|c|c|}
\hline
type of \,$\liem'$ & $\dim(\liem')$ & $\rk(\liem')$ & \,$\liem'$\, maximal & isometry type \\
\hline
\,$(\mathrm{Geo},\vi=t)$\, & $1$ & $1$ & no & \,$\R$\, or \,$\Sph^1$\, \\
\,$(\Sph,\vi=\tfrac\pi6,\ell)$\, & $\ell$ & $1$ & no & \,$\Sph^\ell$\, \\
\,$(\PP,\vi=\tfrac\pi6,(\K,\ell))$\, & $\dim_{\R}\K\cdot \ell$ & $1$ & for \,$(\K,\ell) \in \{(\HH,3),(\OO,2)\}$\,  & $\KP^\ell$ \\
\hline
\,$(\mathrm{AI})$\, & $5$ & $2$ & no & $(\SU(3)/\SO(3))/\Z_3$ \\
\,$(\mathrm{A_2})$\, & $8$ & $2$ & no & $\SU(3)/\Z_3$ \\
\,$(\mathrm{AII})$\, & $14$ & $2$ & yes & \,$(\SU(6)/\Sp(3))/\Z_3$ \\
\,$(\Sph\times\Sph^1,\ell)$\, & \,$\ell+1$\, & $2$ & for \,$\ell=9$\, & \,$(\Sph^\ell\times\Sph^1)/\Z_4$\, \\
\hline
\end{tabular}
\end{center}
\end{Theorem}

\begin{Remark}
\label{R:EIV:EIV:CN}
For the symmetric space \,$\EIV$\,, Chen and Nagano correctly list the \emph{local} isometry types of the maximal totally geodesic submanifolds. However,
the global isometry types of the totally geodesic submanifolds of type \,$(\mathrm{AII})$\, resp.~\,$(\Sph\times\Sph^1,9)$\, is
\,$(\SU(6)/\Sp(3))/\Z_3$\, resp.~\,$(\Sph^1\times\Sph^9)/\Z_4$\, (and not \,$\SU(6)/\Sp(3)$\, resp.~\,$\Sph^1\times\Sph^9$\,, as Chen and Nagano claim). 
\end{Remark}

\emph{Proof of Theorem~\ref{T:EIV:cla}.}
We first mention that it is easily checked using the \textsf{Maple} implementation that the spaces defined in the theorem, and therefore
also the linear subspaces \,$\liem' \subset \liem$\, which are congruent to one of them, are Lie triple systems. It is also
easily seen that the information in the table concerning the dimension and the rank of the Lie triple systems is correct.
The information on the isometry type of the corresponding totally geodesic submanifolds will
be discussed in Section~\ref{SSe:EIV:tgsub}.

We next show that the information on the maximality of the Lie triple systems given in the table is correct.
For this purpose, we presume that the list of Lie triple systems given in the theorem is in fact complete; this will
be proved in the remainder of the present section.

Proof that the Lie triple systems which are claimed to be maximal in the table indeed are: This is clear for the type
\,$(\PP,\vi=\tfrac\pi6,(\OO,2))$\,, because it has the maximal dimension among all the Lie triple systems of \,$\EIV$\,. 
It is also clear for the type \,$(\mathrm{AII})$\, because it has rank \,$2$\, and maximal dimension among all
the Lie triple systems of \,$\EIV$\, of that rank. 
For the type \,$(\Sph\times\Sph^1,9)$\,: For reason of dimension and rank, a Lie triple system \,$\liem'$\, of this type
could only be contained in a Lie triple system of type \,$(\mathrm{AII})$\,; however \,$\liem'$\, 
has a root of multiplicity \,$8$\,, whereas all the roots of Lie triple systems
of type \,$(\mathrm{AII})$\, have multiplicity \,$4$\,, so such an inclusion is impossible.
For the type \,$(\PP,\vi=\tfrac\pi6,(\HH,3))$\,: For reason of dimension, a Lie triple system \,$\liem'$\, of this type
could again only be contained in a Lie triple system of type \,$(\mathrm{AII})$\,; however this is impossible because
\,$\liem'$\, requires the multiplicity \,$8$\, for the ``collapsing'' roots \,$\lambda_1$\, and \,$\lambda_2$\,. 

That no Lie triple systems are maximal besides those mentioned above follows from the following table:

\begin{center}
\begin{tabular}{|c|c|}
\hline
Every Lie triple system of type ... & is contained in a Lie triple system of type ... \\
\hline
$(\mathrm{Geo},\vi=t)$ & $(\Sph\times\Sph^1,1)$ \\
$(\Sph,\vi=\tfrac\pi6,\ell)$ & $(\Sph\times\Sph^1,\ell)$ \\
$(\PP,\vi=\tfrac\pi6,(\K,2))$ with \,$\K\in\{\R,\C,\HH\}$\, & $(\PP,\vi=\tfrac\pi6,(\OO,2))$ \\
$(\PP,\vi=\tfrac\pi6,(\K,3))$ with \,$\K\in\{\R,\C\}$\, & $(\PP,\vi=\tfrac\pi6,(\HH,3))$ \\
\hline
$(\mathrm{AI})$ & $(\mathrm{A}_2)$ \\
$(\mathrm{A}_2)$ & $(\mathrm{AII})$ \\
$(\Sph\times\Sph^1,\ell)$ with \,$\ell \leq 8$\, & $(\Sph\times\Sph^1,9)$ \\
\hline
\end{tabular}
\end{center}

We now turn to the proof that the list of Lie triple systems of \,$\EIV$\, given in Theorem~\ref{T:EIV:cla} is indeed
complete. For this purpose, we let an arbitrary Lie triple system \,$\liem'$\, of \,$\liem$\,, \,$\{0\} \neq \liem' \subsetneq \liem$\,,
be given. 
Because the symmetric space \,$\EIV$\, is of rank \,$2$\,, the rank of \,$\liem'$\, is either \,$1$\, or \,$2$\,. We will handle
these two cases separately in the sequel.

We first suppose that \,$\liem'$\, is a Lie triple system of rank \,$2$\,. Let us fix a Cartan subalgebra \,$\liea$\, of \,$\liem'$\,;
because of \,$\rk(\liem') = \rk(\liem)$\,, \,$\liea$\, is then also a Cartan subalgebra of \,$\liem$\,. In relation to this situation, we use
the notations introduced in Sections~\ref{Se:generallts} and \ref{SSe:EIV:geometry}. In particular,
we consider the positive root system
\,$\Delta_+ := \{\lambda_1,\lambda_2,\lambda_3\}$\, of the root system \,$\Delta := \Delta(\liem,\liea)$\, of \,$\liem$\,,
and also the root system \,$\Delta' := \Delta(\liem',\liea)$\, of \,$\liem'$\,. By Proposition~\ref{P:cla:subroots:subroots-neu}(b), \,$\Delta'$\,
is a root subsystem of \,$\Delta$\,, and therefore \,$\Delta_+' := \Delta' \cap \Delta_+$\, is a positive system of roots for \,$\Delta'$\,.
Moreover, in the root space decompositions of \,$\liem$\, and \,$\liem'$\,
\begin{equation}
\label{eq:EIV:rk2:decomp}
\liem = \liea \;\oplus\; \bigoplus_{\lambda\in\Delta_+} \liem_\lambda \qmq{and}
\liem' = \liea \;\oplus\; \bigoplus_{\lambda\in\Delta_+'} \liem_\lambda'
\end{equation}
the root space \,$\liem_\lambda'$\, of \,$\liem'$\, with respect to \,$\lambda \in \Delta_+'$\, is related to the corresponding root space \,$\liem_\lambda$\, of \,$\liem$\,
by \,$\liem_\lambda' = \liem_\lambda \cap \liem'$\,. 

Because the subset \,$\Delta'$\, of \,$\Delta$\, is invariant under its own Weyl transformation group, we have (up to Weyl transformation) only the
following possibilities for \,$\Delta_+'$\,, which we will treat individually in the sequel:
$$ \Delta_+' = \Delta_+, \quad \Delta_+' = \{\lambda_1\} \qmq{and} \Delta_+' = \varnothing \; . $$

\textbf{The case \,$\boldsymbol{\Delta_+' = \Delta_+}$\,.}
In this case, the restricted Dynkin diagram with multiplicities of \,$\liem'$\, is 
$\xymatrix@=.4cm{ \mathop{\bullet}^{n_{\lambda_1}'} \ar@{-}[r] & \mathop{\bullet}^{n_{\lambda_2}'} }$,
and the classification of the Riemannian symmetric spaces (see, for example, \cite{Loos:1969-2}, p.~119, 146) shows that
\,$n' := n_{\lambda_1}' = n_{\lambda_2}' = n_{\lambda_3}' \in \{1,2,4,8\}$\, holds. 

If \,$n' = 1$\, holds, we may suppose without loss of generality by 
Proposition~\ref{P:EIV:isotropy} that \,$\liem_{\lambda_k}$\, is spanned by \,$v_k := M_{\lambda_k}(1,0,0,0)$\, for \,$k\in\{1,2\}$\,. 
Then we have \,$\liem' \ni R(\lambda_1^\sharp,v_1)v_2 = \tfrac{\sqrt{2}}{8}\,M_{\lambda_3}(0,0,0,i)$\,, and therefore \,$\liem_{\lambda_3}'$\, is 
spanned by \,$v_3 := M_{\lambda_3}(0,0,0,i)$\,. Thus \,$\liem' = \liea \oplus \bigoplus_{k=1}^3 \liem_{\lambda_k}'$\, is of type \,$(\mathrm{AI})$\,.

If \,$n' = 2$\, holds, we may suppose without loss of generality \,$\liem_{\lambda_1}' = M_{\lambda_1}(\C,0,0,0)$\, and \,$v_2 \in \liem_{\lambda_2}'$\,. 
We then obtain \,$v_3 \in \liem_{\lambda_3}'$\, as before, also from the equality \,$\liem' \ni R(\lambda_1^\sharp,M_{\lambda_1}(i,0,0,0))v_2 = 
-\tfrac{\sqrt{2}}{8}\,M_{\lambda_3}(0,0,0,1)$\, the fact \,$M_{\lambda_3}(0,0,0,1) \in \liem_{\lambda_3}'$\, and then from the equality
\,$\liem' \ni R(\lambda_1^\sharp,v_1)M_{\lambda_3}(0,0,0,1) = -\tfrac{\sqrt{2}}{8}\,M_{\lambda_2}(i,0,0,0)$\, the fact
\,$M_{\lambda_2}(i,0,0,0) \in\liem_{\lambda_2}'$\,. Thus we have besides \,$\liem_{\lambda_1}' = M_{\lambda_1}(\C,0,0,0)$\, also 
\,$\liem_{\lambda_2}' = M_{\lambda_2}(\C,0,0,0)$\, and \,$\liem_{\lambda_3}' = M_{\lambda_3}(0,0,0,\C)$\,, and therefore
\,$\liem' = \liea \oplus \bigoplus_{k=1}^3 \liem_{\lambda_k}'$\, is of type \,$(\mathrm{A_2})$\,.

If \,$n' = 4$\, holds, we may suppose without loss of generality \,$\liem_{\lambda_1}' = M_{\lambda_1}(\C,\C,0,0)$\, and \,$v_2 \in \liem_{\lambda_2}'$\,. 
Then as above we obtain \,$M_{\lambda_2}(\C,0,0,0) \subset \liem_{\lambda_2}'$\, and \,$M_{\lambda_3}'(0,0,0,\C) \subset \liem_{\lambda_3}'$\,. 
Moreover for \,$c\in\C$\, we have \,$\liem' \ni R(\lambda_1^\sharp,M_{\lambda_1}(0,c,0,0))v_2 = \tfrac{\sqrt{2}}{8}\,M_{\lambda_3}(0,0,\overline{c}\,i,0)$\,,
hence \,$\liem_{\lambda_3}' = M_{\lambda_3}(0,0,\C,\C)$\,, and \,$\liem' \ni R(\lambda_1^\sharp,v_1)M_{\lambda_3}(0,0,c,0)
= \tfrac{\sqrt{2}}{8}\,M_{\lambda_2}(0,\overline{c}\,i,0,0)$\,, hence \,$\liem_{\lambda_2}' = M_{\lambda_2}(\C,\C,0,0)$\,. 
This shows that \,$\liem' = \liea \oplus \bigoplus_{k=1}^3 \liem_{\lambda_k}'$\, is of type \,$(\mathrm{AII})$\,.

Finally, if \,$n' = 8$\, holds, we have \,$\liem_{\lambda_k}' = \liem_{\lambda_k}$\, for \,$k\in\{1,2,3\}$\, and therefore
\,$\liem' = \liea \oplus \bigoplus_{k=1}^3 \liem_{\lambda_k}' = \liem$\,.

\textbf{The case \,$\boldsymbol{\Delta_+' = \{\lambda_1\}}$\,.}
In this case we have \,$\liem' = \liea \oplus \liem_{\lambda_1}'$\, with a linear subspace \,$\liem_{\lambda_1}' \subset \liem_{\lambda_1}$\,,
and therefore \,$\liem'$\, is of type \,$(\Sph\times\Sph^1,\ell)$\, with \,$\ell := 1 + n_{\lambda_1}' \leq 9$\,. 

\textbf{The case \,$\boldsymbol{\Delta_+' = \varnothing}$\,.}
In this case we have \,$\liem' = \liea$\,, and therefore \,$\liem'$\, is of type \,$(\Sph\times\Sph^1,1)$\,.

\bigskip

We now turn our attention to the case where \,$\liem'$\, is a Lie triple system of rank \,$1$\,. Via the application of the isotropy action of \,$\EIV$\,, we may
suppose without loss of generality 
that \,$\liem'$\, contains a unit vector \,$H$\, from the closure \,$\overline{\liec}$\, 
of the positive Weyl chamber \,$\liec$\, of \,$\liem$\, (with respect to \,$\liea$\, and our
choice of positive roots). Then we have by Equations~\eqref{eq:EIV:isotropy:liec} and \eqref{eq:EIV:isotropy:vivt}
with \,$\vi_0 := \vi(H) \in [0,\tfrac{\pi}{3}]$\, 
\begin{equation}
\label{eq:EIV:rk1:H}
H = \cos(\vi_0)\,\tfrac{\lambda_1^\sharp + \lambda_3^\sharp}{\sqrt{3}} + \sin(\vi_0)\,\lambda_2^\sharp \; .
\end{equation}

Because of \,$\rk(\liem')=1$\,, \,$\liea' := \R\,H$\, is a Cartan subalgebra of \,$\liem'$\,, and we have \,$\liea' = \liea \cap \liem'$\,. 
It follows from Proposition~\ref{P:cla:subroots:subroots-neu}(a) that the root systems \,$\Delta'$\,
and \,$\Delta$\, of \,$\liem'$\, resp.~\,$\liem$\, with respect to \,$\liea'$\, resp.~to \,$\liea$\, are related by
\begin{equation}
\label{eq:EIV:rk1:Delta'Delta}
\Delta' \;\subset\; \Mengegr{\lambda(H)\,\alpha_0}{\lambda\in\Delta,\,\lambda(H)\neq 0}
\end{equation}
with the linear form \,$\alpha_0: \liea' \to \R,\; tH\mapsto t$\,; moreover for \,$\liem'$\, we have the root space decomposition
\begin{equation}
\label{eq:EIV:rk1:m'decomp}
\liem' = \liea' \oplus \bigoplus_{\alpha\in\Delta_+'} \liem_\alpha' 
\end{equation}
where for any root \,$\alpha\in\Delta'$\,, the corresponding root space \,$\liem_\alpha'$\, is given by
\begin{equation}
\label{eq:EIV:rk1:malpha'}
\liem_\alpha' = \left( \bigoplus_{\substack{\lambda \in \Delta \\ \lambda(H) = \alpha(H)}} \liem_\lambda \right) \;\cap\; \liem' \; . 
\end{equation}

If \,$\Delta' = \varnothing$\, holds, then we have \,$\liem' = \R H$\,, and therefore \,$\liem'$\, is then of type \,$(\mathrm{Geo},\vi=\vi_0)$\,. Otherwise
it follows from Proposition~\ref{P:cla:subroots:Comp} that one of the following two conditions holds: Either \,$H$\, is proportional to a root vector \,$\lambda^\sharp$\,
with \,$\lambda\in\Delta$\,, or there exist two \,$\lambda,\mu\in\Delta$\, (\,$\lambda\neq\mu$\,) so that \,$H$\, is orthogonal to \,$\lambda^\sharp-\mu^\sharp$\,.
Evaluating all possible values for \,$\lambda$\, and \,$\mu$\,, we see that \,$\vi_0\in\{0,\tfrac\pi6,\tfrac\pi3\}$\, holds.

In the sequel we consider the three possible values for \,$\vi_0$\, individually.

\textbf{The case \,$\boldsymbol{\vi_0=0}$\,.} In this case we have \,$H = \tfrac1{\sqrt{3}}\,(\lambda_1^\sharp+\lambda_3^\sharp) = \tfrac1{\sqrt{3}}\,(2\lambda_1^\sharp
+\lambda_2^\sharp)$\, 
by Equation~\eqref{eq:EIII:rk1:H} and therefore
$$ \lambda_1(H) = \tfrac12\,\sqrt{3},\quad \lambda_2(H) = 0,\quad \lambda_3(H) = \tfrac12\,\sqrt{3} \; . $$
Thus we have \,$\Delta' = \{\pm \alpha\}$\, with \,$\alpha := \lambda_1|\liea' = \lambda_3|\liea'$\, by Equation~\eqref{eq:EIV:rk1:Delta'Delta},
\,$\liem' = \R\,H \oplus \liem_{\alpha}'$\, by Equation~\eqref{eq:EIV:rk1:m'decomp}, and
\,$\liem_{\alpha}' \subset \liem_{\lambda_1} \oplus \liem_{\lambda_3}$\, by Equation~\eqref{eq:EIII:rk1:malpha'}.

Assume that \,$\liem_{\alpha}' \neq \{0\}$\, holds.
We have \,$\alpha^\sharp = \tfrac12\,\sqrt{3}\,H = \tfrac12\,\lambda_1^\sharp + \tfrac12\,\lambda_3^\sharp$\, and therefore by Proposition~\ref{P:cla:skew}, for any 
\,$v \in \liem_{\alpha}'$\,, say \,$v = M_{\lambda_1}(a_1,\dotsc,a_4) + M_{\lambda_3}(b_1,\dotsc,b_4)$\, with \,$a_1,\dotsc,a_4,b_1,\dotsc,b_4 \in \C$\,,
we have \,$\|a\| = \|b\|$\,. Therefore we can suppose without loss of generality via Proposition~\ref{P:EIV:isotropy} that \,$v_0 := M_{\lambda_1}(1,0,0,0)
+ M_{\lambda_3}(1,0,0,0) \in \liem_{\alpha}'$\, holds. Then we have \,$\liem' \ni R(H,v_0)v_0 = \tfrac34\,H + \tfrac{\sqrt{6}}{8}\,M_{\lambda_2}(0,0,0,i)$\,.
However, this is a contradiction to the fact that because of \,$\lambda_2(H)=0$\,, no element of \,$\liem'$\, can have a non-zero \,$\liem_{\lambda_2}$-component.
So we in fact have \,$\liem_{\alpha}' = \{0\}$\,, hence \,$\liem' = \R\,H$\,. This shows that (for \,$\dim(\liem') \geq 2$\,) the case \,$\vi_0=0$\,
cannot in fact occur.

\textbf{The case \,$\boldsymbol{\vi_0=\tfrac\pi6}$\,.} In this case we have 
\,$H = \tfrac{\sqrt{3}}2\cdot \tfrac1{\sqrt{3}}\,(\lambda_1^\sharp+\lambda_3^\sharp) + \tfrac12\,\lambda_2^\sharp = \lambda_3^\sharp$\,
by Equation~\eqref{eq:EIII:rk1:H} and therefore
$$ \lambda_1(H) = \tfrac12,\quad \lambda_2(H) = \tfrac12,\quad \lambda_3(H) = 1 \; . $$
Thus we have \,$\Delta' \subset \{\pm \alpha\}$\, with \,$\alpha := \lambda_1|\liea' = \lambda_2|\liea'$\, by Equation~\eqref{eq:EIV:rk1:Delta'Delta},
\,$\liem' = \R\,H \oplus \liem_{\alpha}' \oplus \liem_{2\alpha}'$\, by Equation~\eqref{eq:EIV:rk1:m'decomp}, and
\,$\liem_{\alpha}' \subset \liem_{\lambda_1} \oplus \liem_{\lambda_2}$\, and \,$\liem_{2\alpha}' \subset \liem_{\lambda_3}$\, by Equation~\eqref{eq:EIII:rk1:malpha'}.

If \,$\alpha\not\in\Delta'$\, holds, we thus have \,$\liem' = \R\,\lambda_1^\sharp \oplus \liem_{2\alpha}' \subset \R\,\lambda_1^\sharp \oplus \liem_{\lambda_3}$\,,
and therefore \,$\liem'$\, then is of type \,$(\Sph,\vi=\tfrac\pi6,\ell)$\, with \,$\ell := 1+n_{2\alpha}'$\,.

So we now suppose \,$\alpha\in\Delta'$\,. By the classification of the Riemannian symmetric spaces of rank \,$1$\, we then have \,$n_{2\alpha}' \in \{0,1,3,7\}$\,,
and the totally geodesic submanifold corresponding to \,$\liem'$\, is isometric either to \,$\RP^k$\,, to \,$\CP^k$\,, to \,$\HP^k$\, or to 
the Cayley projective plane \,$\OP^{k=2}$\,,
depending on whether \,$n_{2\alpha}'$\, equals \,$0$\,, \,$1$\,, \,$3$\, or \,$7$\,, respectively; here we have \,$k = n_{\alpha}' / (n_{2\alpha}'+1)$\,. 

It should also be noted that we have \,$\alpha^\sharp = \tfrac12\,H = \tfrac12\,\lambda_1^\sharp + \tfrac12\,\lambda_2^\sharp$\,, and therefore we
have for any \,$c_1,\dotsc,c_4,d_1,\dotsc,d_4 \in \C$\, by Proposition~\ref{P:cla:skew}
\begin{equation}
\label{eq:EIV:rk1:pi6:alpha-ab}
M_{\lambda_1}(c_1,\dotsc,c_4) + M_{\lambda_2}(d_1,\dotsc,d_4) \in \liem_\alpha' \;\Longrightarrow\; \|c\| = \|d\| \; .
\end{equation}

In the sequel, we consider the four possible values for \,$n_{2\alpha}'$\, individually. In our calculations we will use the vectors \,$v_0,v_0^C,\dotsc$\,
as they are defined in the entry for the types \,$(\PP,\vi=\tfrac\pi6,(\K,\ell))$\, in Theorem~\ref{T:EIV:cla}. 

Let us first suppose \,$n_{2\alpha}'=0$\,, i.e.~\,$\Delta' = \{\pm \alpha\}$\,. By Proposition~\ref{P:EIV:isotropy} and because of \eqref{eq:EIV:rk1:pi6:alpha-ab}
we may suppose without loss of generality that \,$v_0 \in \liem_\alpha'$\, holds. If \,$n_\alpha'=1$\, holds, we then have \,$\liem_\alpha' = \R v_0$\,
and therefore \,$\liem' = \R H \oplus \liem_\alpha'$\, is of type \,$(\PP,\vi=\tfrac\pi6,(\R,2))$\,. Otherwise we choose \,$v\in\liem_\alpha'$\, to
be orthogonal to \,$v_0$\,, say \,$v = M_{\lambda_1}(c_1,\dotsc,c_4) + M_{\lambda_2}(d_1,\dotsc,d_4)$\,. Then we have
$$ \liem' \ni R(H,v_0)v = \tfrac14\RE(c_1)\,\lambda_1^\sharp + \tfrac14\RE(d_1)\,\lambda_2^\sharp 
+ \tfrac{\sqrt{2}}{16}\,M_{\lambda_3}\bigr(\,i(\overline{c_4}-\overline{d_4})\,,\, i(-\overline{c_3}+\overline{d_3}) \,,\, -i(\overline{c_2}+\overline{d_2})
\,,\, i(\overline{d_1}-c_1)\,\bigr) \; . $$
Because the \,$\liea$-component of this vector is proportional to \,$H$\,, we have \,$\RE(c_1) = \RE(d_1)$\,; this equation together with our requirement
that \,$v$\, be orthogonal to \,$v_0$\, shows \,$\RE(c_1) = \RE(d_1) = 0$\, and hence
\,$ c_1,d_1 \in i\R$\,. 
Moreover, because of \,$2\alpha\not\in\Delta'$\,, the \,$\liem_{\lambda_3}$-component of the above vector vanishes, and thus we have
$$ c_1 = -d_1\in i\R, \quad c_2 = -d_2, \quad c_3 = d_3 \qmq{and} c_4 = d_4 \; , $$
hence 
\begin{equation}
\label{eq:EIV:rk1:pi6:0:v}
v = M_{\lambda_1}(it,c_2,c_3,c_4) + M_{\lambda_2}(-it,-c_2,c_3,c_4)
\end{equation}
with \,$t\in\R$\,. By application of another isotropy transformation, we can now arrange that \,$v$\, is proportional to \,$v_1$\,. Thus we have
\,$v_0,v_1\in\liem_\alpha'$\,, and therefore in the case \,$n_\alpha'=2$\,, \,$\liem' = \R H \oplus \liem_\alpha'$\, is of type \,$(\PP,\vi=\tfrac\pi6,(\R,3))$\,.
We now show that the case \,$n_\alpha' \geq 3$\, does not occur. For this purpose, we again let \,$v\in \liem_\alpha'$\, be given, but now suppose
that \,$v$\, is orthogonal to both \,$v_0$\, and \,$v_1$\,. Then \,$v$\, again has the form of Equation~\eqref{eq:EIV:rk1:pi6:0:v}, however
the requirement that \,$v$\, be orthogonal to \,$v_1$\, implies \,$t=0$\,. Moreover, we have
$$ \liem' \ni R(H,v)v_1 = \tfrac{\sqrt{2}}{16}\,M_{\lambda_3}(\overline{c_4}\,i,-\overline{c_3}\,i,-\overline{c_2}\,i,0) \; . $$
Because of \,$2\alpha\not\in\Delta'$\,, the \,$\liem_{\lambda_3}$-component of this vector vanishes, and thus we have \,$c_2=c_3=c_4=0$\,, hence \,$v=0$\,.
This shows that \,$n_{\alpha}' \geq 3$\, is impossible.

Next we suppose \,$n_{2\alpha}'=1$\,. Then the Lie triple system \,$\liem'$\, corresponds to a complex projective space \,$\CP^\ell$\,, which is a Hermitian symmetric
space. Let \,$\liem'' \subset \liem'$\, be the tangent space of a real form of this space, 
then \,$\liem''$\, will also be a Lie triple system of \,$\liem$\,, it will be of rank \,$1$\, and
correspond to the isotropy angle \,$\vi=\tfrac\pi6$\,, and it will have only the root \,$\alpha$\,, not \,$2\alpha$\,. As a consequence of the preceding
classification of the Lie triple systems with these properties, \,$\liem''$\, is of type \,$(\PP,\vi=\tfrac\pi6,(\R,\ell))$\, with \,$\ell\in\{2,3\}$\,. 
Without loss of generality,
we may therefore suppose that \,$\liem''$\, is the prototype Lie triple system of the type \,$(\PP,\vi=\tfrac\pi6,(\R,\ell))$\, as given in Theorem~\ref{T:EIV:cla}. 
Thus we have \,$v_0\in\liem_\alpha'$\, and in the case \,$\ell=3$\, also \,$v_1\in\liem_\alpha'$\,. Further we may suppose without loss of generality
\,$\liem_{2\alpha}' = \R\,w_1$\,. Then we have
$$  R(w_1,v_k)H = \tfrac{\sqrt{2}}{16}\,v_k^C $$
for \,$k\in\{0,1\}$\,, and therefore \,$v_k \in \liem_\alpha'$\, implies also \,$v_k^C \in \liem_\alpha'$\,. This shows that \,$\liem' = \R H \oplus \liem_\alpha'
\oplus \liem_{2\alpha}'$\, is of type \,$(\PP,\vi=\tfrac\pi6,(\C,\ell))$\,.

Now we suppose \,$n_{2\alpha}'=3$\,. Then \,$\liem'$\, corresponds to a quaternionic projective space \,$\HP^\ell$\,, and therefore an analogous
argument as in the treatment of the case \,$n_{2\alpha}'=1$\, shows that \,$\liem'$\, contains as a complex form a Lie triple system \,$\liem''$\,
of type \,$(\PP,\vi=\tfrac\pi6,(\C,\ell))$\, with \,$\ell\in\{2,3\}$\,. Without loss of generality, we may suppose that \,$\liem''$\, is the
prototype Lie triple system of that type as given in Theorem~\ref{T:EIV:cla}, and therefore \,$w_1 \in \liem_{2\alpha}'$\, and 
\,$v_k,v_k^C\in\liem_\alpha'$\, holds, where \,$k=0$\, for \,$\ell=2$\, and \,$k=0,1$\, for \,$\ell=3$\,. Further we may suppose without loss of generality
that also \,$w_2 \in \liem_{2\alpha}'$\, holds. We have for \,$k\in\{0,1\}$\,
$$ R(w_2,v_k)H = -\tfrac{\sqrt{2}}{16}\,v_k^H \qmq{and} R(w_2,v_k^C)H = -\tfrac{\sqrt{2}}{16}\,v_k^{CH} \; . $$ 
Therefore \,$v_k,v_k^C \in \liem_\alpha'$\, implies also \,$v_k^H,v_k^{CH}\in\liem_\alpha'$\, Moreover, we have
$$ R(v_0,v_0^{CH})H = \tfrac{\sqrt{2}}{4}\,w_3 \;, $$
and therefore \,$w_3 \in \liem_{2\alpha}'$\,. Thus \,$\liem' = \R H \oplus \liem_\alpha' \oplus \liem_{2\alpha}'$\, is of type
\,$(\PP,\vi=\tfrac\pi6,(\HH,\ell))$\,. 

Finally we suppose \,$n_{2\alpha}'=7$\,. Then \,$n_\alpha'=8$\, is the only possibility by the classification of the Riemannian symmetric spaces of rank \,$1$\,,
and \,$\liem'$\, corresponds to the Cayley projective plane \,$\OP^2$\,. \,$\OP^2$\, contains a \,$\HP^2$\, as totally geodesic submanifold, and
thus by an analogous argument as before, we see that \,$\liem'$\, contains a Lie triple system \,$\liem''$\, of type \,$(\PP,\vi=\tfrac\pi6,(\HH,2))$\,;
without loss of generality we may suppose that \,$\liem''$\, is the prototype Lie triple system of that type. Thus we have \,$v_0,v_0^C,v_0^H,v_0^{CH}\in\liem_\alpha'$\,
and \,$w_1,w_2,w_3 \in \liem_{2\alpha}'$\,. Without loss of generality we may further suppose \,$w_4 \in \liem_{2\alpha}'$\,. We have
\begin{gather*}
R(w_4,v_0)H = \tfrac{\sqrt{2}}{16}\,v_0^O \;, \quad R(w_4,v_0^C)H = \tfrac{\sqrt{2}}{16}\,v_0^{CO} \;, \\
R(w_4,v_0^H)H = \tfrac{\sqrt{2}}{16}\,v_0^{HO} \qmq{and} \quad R(w_4,v_0^{CH})H = \tfrac{\sqrt{2}}{16}\,v_0^{CHO} \;, 
\end{gather*} 
and therefore \,$\liem_{\alpha}'$\, is spanned by \,$v_0,v_0^C,v_0^H,v_0^{CH},v_0^O,v_0^{CO},v_0^{HO},v_0^{CHO}$\,. Moreover, we have
$$ R(v_0,v_0^{CO})H = \tfrac{\sqrt{2}}{4}\, w_5\;,\quad  R(v_0,v_0^{HO})H = \tfrac{\sqrt{2}}{4}\, w_6 \qmq{and}  R(v_0,v_0^{CHO})H = \tfrac{\sqrt{2}}{4}\, w_7 $$
and therefore \,$\liem_{2\alpha}'$\, is spanned by \,$w_1,\dotsc,w_7$\,. Therefore \,$\liem' = \R H \oplus \liem_\alpha' \oplus \liem_{2\alpha}'$\, is
of type \,$(\PP,\vi=\tfrac\pi6,(\OO,2))$\,.

\textbf{The case \,$\boldsymbol{\vi_0=\tfrac\pi3}$\,.} By an analogous argument as in the case \,$\vi_0=0$\,, one shows that this case cannot occur.

This completes the classification of the Lie triple systems in the Riemannian symmetric space \,$\EIV$\,. 
\strut\hfill$\Box$

\subsection[Totally geodesic submanifolds in \,$E_6/F_4$\,]{Totally geodesic submanifolds in \,$\boldsymbol{E_6/F_4}$\,}
\label{SSe:EIV:tgsub}

We are interested in determining the global isometry types of the totally geodesic submanifolds of \,$\EIV$\, corresponding to the various
types of Lie triple systems as they were classified in Theorem~\ref{T:EIV:cla}. In the case of \,$\EIV$\, all the maximal totally geodesic
submanifolds are reflective, so we can derive this information from the classification of reflective submanifolds due to \textsc{Leung},
see \cite{Leung:reflective-1979}. 

Using the information from that paper, we obtain the results of the following table. In it, we again use the notations introduced at
the beginning of Section~\ref{SSe:EIII:tgsub}.

\begin{center}
\begin{tabular}{|c|c|c|}
\hline
type of Lie triple system & isometry type & properties\footnotemark \\
\hline
$(\mathrm{Geo},\vi=t)$ & \,$\R$\, or \,$\Sph^1$ & \\
$(\Sph,\vi=\tfrac\pi6,\ell)$ & $\Sph^\ell_{r=1}$ & \,$\ell=9$\,: Helgason sphere \\
$(\PP,\vi=\tfrac\pi6,(\K,\ell))$ & $\KP^\ell_{\vkap=1/4}$ & \,$(\K,\ell)=(\OO,2)$\,: polar, maximal \\
& & \,$(\K,\ell)=(\HH,3)$\,: reflective, maximal \\
$(\mathrm{AI})$ & $((\SU(3)/\SO(3))/\Z_3)_{\mathrm{srr}=1}$ & \\
$(\mathrm{A}_2)$ & $(\SU(3)/\Z_3)_{\mathrm{srr}=1}$ & \\
$(\mathrm{AII})$ & $((\SU(6)/\Sp(3))/\Z_3)_{\mathrm{srr}=1}$ & reflective, maximal \\
$(\Sph\times\Sph^1,\ell)$ & $(\Sph^\ell_{r=1}\times\Sph^1_{r=\sqrt{3}})/\Z_4$ & \,$\ell=9$\,: meridian for \,$(\PP,\vi=\tfrac\pi6,(\OO,2))$\,, maximal \\
\hline
\end{tabular}
\end{center}
\footnotetext{The polars and meridians are also reflective, without this fact being noted explicitly in the table.}

\subsection[Totally geodesic submanifolds in \,$\SU(6)/\Sp(3)$\,]{Totally geodesic submanifolds in \,$\boldsymbol{\SU(6)/\Sp(3)}$\,}
\label{SSe:EIV:AII}

Similarly as we derived the classification of the Lie triple systems resp.~the totally geodesic submanifolds in \,$\SO(10)/\Ug(5)$\, from that
classification in \,$\EIII$\, in Section~\ref{SSe:EIII:DIII}, we now derive the classification for \,$\SU(6)/\Sp(3)$\, from the classification
in \,$\EIV$\,, using the fact that \,$\SU(6)/\Sp(3)$\, is the local isometry type of a maximal totally geodesic submanifold of \,$\EIV$\,. 

Thus we remain in the situation studied in Section~\ref{SSe:EIV:lts}. We consider the Riemannian symmetric space \,$\EIV$\,, and let \,$\lieg=\liek\oplus\liem$\,
be the canonical decomposition of \,$\lieg=\liee_6$\, associated with this space, i.e.~we have \,$\liek = \lief_4$\, and \,$\liem$\, is isomorphic
to the tangent space of \,$\EIV$\,. We will use the names for the types of Lie triple systems of \,$\liem$\, as introduced in Theorem~\ref{T:EIV:cla}.

Further, we let \,$\liem_1$\, be a Lie triple system of \,$\liem$\, of type \,$(\mathrm{AII})$\,, i.e.~\,$\liem_1$\, corresponds to a totally geodesic submanifold
which is locally isometric to \,$\SU(6)/\Sp(3)$\,.

\begin{Theorem}
\label{EIV:AII:cla}
Exactly the following types of Lie triple systems of \,$\EIV$\, have representatives which are contained in \,$\liem_1$\,:
\begin{itemize}
\item \,$(\mathrm{Geo},\vi=t)$\, with \,$t \in [0,\tfrac\pi3]$\,
\item \,$(\Sph,\vi=\tfrac\pi6,\ell)$\, with \,$\ell \leq 5$\,
\item \,$(\PP,\vi=\tfrac\pi6,(\K,2))$\, with \,$\K\in\{\R,\C,\HH\}$\,
\item \,$(\PP,\vi=\tfrac\pi6,(\K,3))$\, with \,$\K\in\{\R,\C\}$\,
\item \,$(\mathrm{AI})$\,
\item \,$(\mathrm{A}_2)$\,
\item \,$(\Sph\times\Sph^1,\ell)$\, with \,$\ell\leq 5$\,
\end{itemize}
The maximal Lie triple systems of \,$\liem_1$\, are those of the types: \,$(\PP,\vi=\tfrac\pi6,(\HH,2))$\,, 
\,$(\PP,\vi=\tfrac\pi6,(\C,3))$\,, \,$(\mathrm{A}_2)$\, and \,$(\Sph\times\Sph^1,5)$\,. 
\end{Theorem}

\beweis
Similar to the proofs of Theorems~\ref{EIII:SO5:cla} and \ref{EIII:DIII:cla}.
\beweisende

\begin{Remark}
\label{R:EIV:AII:CN}
Chen/Nagano incorrectly state in \cite{Chen/Nagano:totges2-1978} that the Lie triple systems of type \,$(\mathrm{AI})$\,
(corresponding to \,$\SU(3)/\SO(3)$\,) were maximal in \,$\SU(6)/\Sp(3)$\,, rather these Lie triple systems are contained in Lie triple systems
of type \,$(\mathrm{A}_2)$\, (corresponding to \,$\SU(3)$\,).
\end{Remark}

Also for \,$\SU(6)/\Sp(3)$\,, the maximal totally geodesic submanifolds are all reflective. Using the information from \cite{Leung:reflective-1979},
we obtain the following information on the global isometry type of the totally geodesic submanifolds of \,$\SU(6)/\Sp(3)$\, 
corresponding to the various types of Lie triple systems:

\begin{center}
\begin{tabular}{|c|c|c|}
\hline
type of Lie triple system & isometry type & properties\footnotemark \\
\hline
$(\mathrm{Geo},\vi=t)$ & \,$\R$\, or \,$\Sph^1$ & \\
$(\Sph,\vi=\tfrac\pi6,\ell)$ & $\Sph^\ell_{r=1}$ & \,$\ell=5$\,: Helgason sphere \\
$(\PP,\vi=\tfrac\pi6,(\K,\ell))$ & $\KP^\ell_{\vkap=1/4}$ & \,$(\K,\ell)=(\HH,2)$\,: polar, maximal \\
& & \,$(\K,\ell)=(\C,3)$\,: reflective, maximal \\
$(\mathrm{AI})$ & $(\SU(3)/\SO(3))_{\mathrm{srr}=1}$ & \\
$(\mathrm{A}_2)$ & $\SU(3)_{\mathrm{srr}=1}$ & reflective \\
$(\Sph\times\Sph^1,\ell)$ & $(\Sph^\ell_{r=1}\times\Sph^1_{r=\sqrt{3}})/\Z_2$ & \,$\ell=5$\,: meridian for \,$(\PP,\vi=\tfrac\pi6,(\HH,2))$\,, maximal \\
\hline
\end{tabular}
\end{center}
\footnotetext{The polars and meridians are also reflective, without this fact being noted explicitly in the table.}

\subsection[Totally geodesic submanifolds in \,$\SU(3)$\,]{Totally geodesic submanifolds in \,$\boldsymbol{\SU(3)}$\,}
\label{SSe:EIV:A2}

Using the same strategy as before, we next classify the totally geodesic submanifolds of \,$\SU(3)$\,, regarded as a Riemannian symmetric space. 
We again let \,$\lieg=\liek\oplus\liem$\,
be the splitting corresponding to \,$\EIV$\,, and let \,$\liem_1$\, now be a Lie triple system of \,$\liem$\, of type \,$(\mathrm{A}_2)$\,;
then the totally geodesic submanifold of \,$\EIV$\, corresponding to \,$\liem_1$\, is locally isometric to \,$\SU(3)$\,.

\begin{Theorem}
\label{EIV:A2:cla}
Exactly the following types of Lie triple systems of \,$\EIV$\, have representatives which are contained in \,$\liem_1$\,:
\begin{itemize}
\item \,$(\mathrm{Geo},\vi=t)$\, with \,$t \in [0,\tfrac\pi3]$\,
\item \,$(\Sph,\vi=\tfrac\pi6,\ell)$\, with \,$\ell \leq 3$\,
\item \,$(\PP,\vi=\tfrac\pi6,(\K,2))$\, with \,$\K\in\{\R,\C\}$\,
\item \,$(\PP,\vi=\tfrac\pi6,(\R,3))$\, 
\item \,$(\mathrm{AI})$\,
\item \,$(\Sph\times\Sph^1,\ell)$\, with \,$\ell\leq 3$\,
\end{itemize}
The maximal Lie triple systems of \,$\liem_1$\, are those of the types: \,$(\PP,\vi=\tfrac\pi6,(\C,2))$\,, 
\,$(\PP,\vi=\tfrac\pi6,(\R,3))$\,, \,$(\mathrm{AI})$\, and \,$(\Sph\times\Sph^1,3)$\,. 
\end{Theorem}

\beweis
Similar to the proofs of Theorems~\ref{EIII:SO5:cla} and \ref{EIII:DIII:cla}.
\beweisende

\begin{Remark}
\label{R:EIV:A2:CN}
Chen/Nagano incorrectly state in \cite{Chen/Nagano:totges2-1978} that \,$\SU(3)$\, contains totally geodesic submanifolds
isometric to \,$\SU(2)\times\SU(2)$\, and \,$\SU(3)/(\SU(2)\times\SU(2))$\,. This is impossible, because \,$\SU(2)\times\SU(2)$\,
has the same rank as \,$\SU(3)$\,, but whereas the former group has two orthogonal roots, the latter has not.
\end{Remark}

Once again, also for the Riemannian symmetric space \,$\SU(3)$\,, all the maximal totally geodesic submanifolds are reflective. 
Using the classification of the reflective submanifolds by Leung (in \cite{Leung:reflective-1974}, Theorem~3.3 for the group manifolds,
see also \cite{Leung:reflective-errata-1975}),
we obtain the following information on the global isometry type of the totally geodesic submanifolds of \,$\SU(3)$\, 
corresponding to the various types of Lie triple systems:

\begin{center}
\begin{tabular}{|c|c|c|}
\hline
type of Lie triple system & isometry type & properties\footnotemark \\
\hline
$(\mathrm{Geo},\vi=t)$ & \,$\R$\, or \,$\Sph^1$ & \\
$(\Sph,\vi=\tfrac\pi6,\ell)$ & $\Sph^\ell_{r=1}$ & \,$\ell=3$\,: Helgason sphere \\
$(\PP,\vi=\tfrac\pi6,(\K,\ell))$ & $\KP^\ell_{\vkap=1/4}$ & \,$(\K,\ell)=(\C,2)$\,: polar, maximal \\
& & \,$(\K,\ell)=(\R,3)$\,: reflective, maximal \\
$(\mathrm{AI})$ & $(\SU(3)/\SO(3))_{\mathrm{srr}=1}$ & reflective, maximal \\
$(\Sph\times\Sph^1,\ell)$ & $(\Sph^\ell_{r=1}\times\Sph^1_{r=\sqrt{3}})/\Z_2$ & \,$\ell=3$\,: meridian for \,$(\PP,\vi=\tfrac\pi6,(\C,2))$\,, maximal \\
\hline
\end{tabular}
\end{center}
\footnotetext{The polars and meridians are also reflective, without this fact being noted explicitly in the table.}

\subsection[Totally geodesic submanifolds in \,$\SU(3)/\SO(3)$\,]{Totally geodesic submanifolds in \,$\boldsymbol{\SU(3)/\SO(3)}$\,}
\label{SSe:EIV:AI}

The totally geodesic submanifolds of \,$\SU(3)/\SO(3)$\, have already been classified in \cite{Klein:2007-Satake}, 
Section~6. Because the totally geodesic submanifolds of \,$\EIV$\, of type \,$(\mathrm{AI})$\, are locally isometric to \,$\SU(3)/\SO(3)$\,,
the Lie triple systems of \,$\SU(3)/\SO(3)$\, also occur as Lie triple systems of \,$\EIV$\,. In the following table, we list the correspondence
between the types 
of Lie triple systems of \,$\SU(3)/\SO(3)$\, as defined in \cite{Klein:2007-Satake}, Proposition~6.1,
and types of Lie triple systems of \,$\EIV$\, as defined in Theorem~\ref{T:EIV:cla} of the present paper. We also give the isometry type of
the corresponding totally geodesic submanifolds, as it has been determined in \cite{Klein:2007-Satake}, Section~6; for the application of this
information it should be noted that there the metric of \,$\SU(3)/\SO(3)$\, has been normalized in such a way that the roots have length \,$\sqrt{2}$\,,
whereas we now want to normalize the metric in such a way that the roots have length \,$1$\,. 

\begin{center}
\begin{tabular}{|c|c|c|c|}
\hline
type (\cite{Klein:2007-Satake}, Prop.~6.1) & type (Thm.~\ref{T:EIV:cla}) & isometry type & properties \\
\hline
(G) & $(\mathrm{Geo},\vi=t)$ & \,$\R$\, or \,$\Sph^1$\, & \\
(T) & $(\Sph\times\Sph^1,1)$ & $(\Sph^1_{r=1}\times\Sph^1_{r=\sqrt{3}})/\Z_2$ & \\
(S) & $(\Sph,\vi=\tfrac\pi6,2)$ & $\Sph^2_{r=1}$ & Helgason sphere \\
(M) & $(\PP,\vi=\tfrac\pi6,(\R,2))$ & $\RP^2_{\vkap=1/4}$ & polar, maximal  \\
(P) & $(\Sph\times\Sph^1,2)$ & $(\Sph^2_{r=1}\times\Sph^1_{r=\sqrt{3}})/\Z_2$ & meridian, maximal \\
\hline
\end{tabular}
\end{center}

\section[The symmetric spaces \,$G_2$\, and \,$G_2/\SO(4)$\,]{The symmetric spaces \,$\boldsymbol{G_2}$\, and \,$\boldsymbol{G_2/\SO(4)}$\,}
\label{Se:G2}

\subsection[The geometry of the Lie group \,$G_2$\,, regarded as a symmetric space]{The geometry of the Lie group \,$\boldsymbol{G_2}$\,, regarded as a symmetric space}
\label{SSe:G2:geometry}

In this section we will study the exceptional compact Lie group \,$G_2$\,, regarded as a Riemannian symmetric space. In particular we need to obtain
its curvature tensor. The usual way to do so would be to regard \,$G_2$\, as the quotient space \,$(G_2 \times G_2)/\Delta(G_2)$\,, where
\,$\Delta(G_2) := \Menge{(g,g)}{g \in G_2}$\, is the diagonal of the product \,$G_2\times G_2$\,, and to apply the results of 
\cite{Klein:2007-Satake} to this space. 

However, we can reduce the effort involved in the calculations by noting that in that model of the symmetric space \,$G_2$\,, the space \,$\liem$\,
which corresponds to the tangent space in the origin, is given by \,$\liem = \Menge{(X,-X)}{X\in\lieg_2} \subset \lieg_2 \oplus \lieg_2$\,, and that
for elements \,$(X,-X),(Y,-Y),(Z,-Z) \in \liem$\,, the curvature tensor is given by
$$ -\bigr[ \, [(X,-X)\,,\,(Y,-Y)]\,,\,(Z,-Z)\,\bigr] = \, -\bigr(\,[[X,Y],Z]\,,\,-[[X,Y],Z]\,\bigr) \; . $$
Under the canonical isomorphism \,$\liem \to \lieg_2,\; (X,-X)\mapsto X$\,, the curvature tensor of these elements of \,$\liem$\, therefore
corresponds to \,$-[[X,Y],Z] \in \lieg_2$\,, hence the Lie triple systems 
in \,$\liem \subset \lieg_2 \oplus \lieg_2$\, correspond to the Lie triple systems in \,$\lieg_2$\, (i.e.~to the linear subspaces of \,$\lieg_2$\, which 
are invariant under the Lie triple bracket \,$[[\,\cdot\,,\,\cdot\,],\,\cdot\,]$\, of \,$\lieg_2$\,).
Moreover, the isotropy action of \,$\Delta(G_2)$\, on \,$\liem$\, corresponds to the adjoint action of \,$G_2$\, on \,$\lieg_2$\,. 
For this reason, we can carry out the classification of Lie triple systems by calculation in \,$\lieg_2$\, itself (instead of in \,$\liem \subset 
\lieg_2 \oplus \lieg_2$\,). In doing so, we will only need the description of the root system and the Lie bracket of \,$\lieg_2$\,, which we obtain 
by application of the results of Sections~2, 3 of \cite{Klein:2007-Satake}.

In the sequel, we will consider also an \,$\Ad(G_2)$-invariant inner product on  \,$\lieg_2$\,. Such an inner product is unique up to a positive constant,
which we choose so that the shortest roots of \,$\lieg_2$\, (see below) have the length \,$1$\,. 

We now fix a Cartan subalgebra \,$\liea \subset \lieg_2$\, and a choice of positive roots in the root system \,$\Delta$\, of \,$\lieg_2$\, with respect to \,$\liea$\,. 
The Dynkin diagram of \,$\lieg_2$\, is $\xymatrix@=.4cm{ \bullet \ar@3{<-}[r] & \bullet }$, and therefore the simple roots of \,$\lieg_2$\,,
which we denote by \,$\lambda_1$\, and \,$\lambda_2$\,, have an angle of \,$\tfrac{5\pi}{6}$\, to each other, where \,$\lambda_2$\, is the longer root by a factor 
of \,$\sqrt{3}$\,. The other positive roots of \,$G_2$\, are
$$ \lambda_3 := \lambda_1+\lambda_2,\quad \lambda_4 := 2\lambda_1+\lambda_2,\quad \lambda_5 := 3\lambda_1 + \lambda_2 \qmq{and} \lambda_6 := 3\lambda_1 + 2\lambda_2\; . $$
In this way we obtain the following root diagram for \,$G_2$\,: 

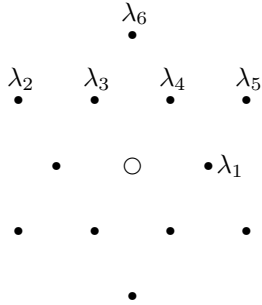
\begin{center}
\begin{minipage}{5cm}
\begin{center}
\strut \\[.3cm]
\setlength{\unitlength}{1cm}
\begin{picture}(3.5,3.5)
\put(2,2){\circle{0.2}}
\put(3,2){\circle*{0.1}}
\put(2.5,2.866){\circle*{0.1}}
\put(1.5,2.866){\circle*{0.1}}
\put(2.5,1.134){\circle*{0.1}}
\put(1.5,1.134){\circle*{0.1}}
\put(1,2){\circle*{0.1}}
\put(2,3.732){\circle*{0.1}}
\put(2,0.268){\circle*{0.1}}
\put(0.5,2.866){\circle*{0.1}}
\put(3.5,2.866){\circle*{0.1}}
\put(0.5,1.134){\circle*{0.1}}
\put(3.5,1.134){\circle*{0.1}}

\put(3.1,1.90){{\small $\lambda_1$}}
\put(0.35,3.05){{\small $\lambda_2$}}
\put(1.35,3.05){{\small $\lambda_3$}}
\put(2.35,3.05){{\small $\lambda_4$}}
\put(3.35,3.05){{\small $\lambda_5$}}
\put(1.85,3.92){{\small $\lambda_6$}}
\end{picture}
\strut \\[1cm]
\end{center}
\end{minipage}
\end{center}
\vspace{.2cm}

In the sequel, we will use the notation \,$V_{\lambda_k}(c)$\, for \,$k\in\{1,\dotsc,6\}$\, and \,$c\in\C$\, to denote an element of the root space of \,$\lieg_2$\,
corresponding to the root \,$\lambda_k$\, as defined in \cite{Klein:2007-Satake}, Proposition~3.3(d). Then the root space corresponding to \,$\lambda_k$\,
equals \,$V_{\lambda_k}(\C)$\,. 

We will also use the isotropy angle function \,$\vi$\, defined at the end of Section~\ref{Se:generallts} for \,$\lieg_2$\,; remember that
in the present situation, the isotropy action of the symmetric space \,$G_2$\, is given simply by the adjoint action of \,$G_2$\, on \,$\lieg_2$\,. 
We have
\,$\vi_{\textrm{max}}=\tfrac\pi6$\, and thus we obtain an isotropy angle function \,$\vi: \lieg_2\setminus\{0\} \to [0,\tfrac\pi6]$\,. 
For the elements of the closure \,$\overline{\liec}$\, of the positive Weyl chamber \,$\liec := \Menge{v\in\liea}{\lambda_1(v)\geq 0,
\lambda_2(v)\geq 0}$\, we once again explicitly describe the relation to their isotropy angle: \,$(\lambda_4^\sharp,\tfrac{1}{\sqrt{3}}\,\lambda_2^\sharp)$\,
is an orthonormal basis of \,$\liea$\, so that with \,$v_t := \cos(t)\lambda_4^\sharp + \sin(t)\,\tfrac{1}{\sqrt{3}}\,\lambda_2^\sharp$\, we have
\begin{equation}
\label{eq:G2:isotropy:liec}
\overline{\liec} = \Menge{s\cdot v_t}{t\in[0,\tfrac\pi6],s\in\R_{\geq 0}}\;, 
\end{equation}
and because the Weyl chamber \,$\liec$\, is bordered
by the two vectors \,$v_0 = \lambda_4^\sharp$\, with \,$\vi(v_0)=0$\, and \,$v_{\pi/6}=\tfrac{1}{\sqrt{3}}\,\lambda_6^\sharp$\, with \,$\vi(v_{\pi/6})=\tfrac\pi6$\,, 
we have
\begin{equation}
\label{eq:G2:isotropy:vivt}
\vi(s\cdot v_t)=t \qmq{for all} t\in[0,\tfrac\pi6]\,, s\in\R_+ \; . 
\end{equation}

Further we note the following simple fact on the adjoint action of \,$G_2$\,:

\begin{Prop}
\label{P:G2:geometry:isotropy}
Let \,$\lambda$\, be a short root, and \,$\lambda'$\, be a long root of \,$G_2$\,. Then the adjoint action of the maximal torus \,$T := \exp(\liea)$\, on \,$\lieg_2$\,
leaves \,$\liea$\, pointwise fixed, and acts ``jointly transitively'' on the unit spheres in the root spaces \,$V_\lambda(\C)$\, and \,$V_{\lambda'}(\C)$\,
in the sense that for any given \,$c_1,c_2,c_1',c_2'\in\C$\, with \,$|c_1|=|c_2|$\, and \,$|c_1'|=|c_2'|$\, there exists \,$g\in T$\, with
\,$\Ad(g)V_\lambda(c_1) = V_\lambda(c_2)$\, and \,$\Ad(g)V_{\lambda'}(c_1') = V_{\lambda'}(c_2')$\,. 
\end{Prop}

\subsection[Lie triple systems in \,$G_2$\,]{Lie triple systems in \,$\boldsymbol{G_2}$\,}
\label{SSe:G2:lts}

We continue to use the notations of the preceding section.

\begin{Theorem}
\label{T:G2:cla}
The linear subspaces \,$\liem' \subset \lieg_2$\, given in the following are Lie triple systems, 
and every Lie triple system \,$\{0\} \neq \liem' \subsetneq \lieg_2$\, is congruent under the adjoint action to one of them. 

\begin{itemize}
\item \,$\boldsymbol{(\mathrm{Geo},\vi=t)}$\, with \,$t\in[0,\tfrac\pi6]$\,. \\
\,$\liem' = \R(\cos(t)\,\lambda_4^\sharp + \sin(t)\,\tfrac{1}{\sqrt{3}}\,\lambda_2^\sharp)$\, (see Equation~\eqref{eq:G2:isotropy:vivt}).

\item \,$\boldsymbol{(\Sph,\vi=0,\ell})$\, with \,$\ell \in \{2,3\}$\,. \\
\,$\liem'$\, is an \,$\ell$-dimensional linear subspace of \,$\R\,\lambda_1^\sharp \oplus \liem_{\lambda_1}$\,. 

\item \,$\boldsymbol{(\Sph,\vi=\arctan(\tfrac{1}{3\sqrt{3}}),\ell)}$\, with \,$\ell \in \{2,3\}$\,. \\
\,$\liem'$\, is an \,$\ell$-dimensional subspace of 
\,$\mathrm{span}\{\,9\lambda_1^\sharp + 5\lambda_2^\sharp,\, V_{\lambda_1}(1) + V_{\lambda_2}(\tfrac13\,\sqrt{5}),\, V_{\lambda_1}(i) + V_{\lambda_2}(\tfrac13\,\sqrt{5}\,i)\}$\,. 

\item \,$\boldsymbol{(\Sph,\vi=\tfrac\pi6,\ell})$\, with \,$\ell \in \{2,3\}$\,. \\
\,$\liem'$\, is an \,$\ell$-dimensional linear subspace of \,$\R\,\lambda_6^\sharp \oplus \liem_{\lambda_6}$\,. 

\item \,$\boldsymbol{(\PP,\vi=\tfrac\pi6,(\R,\ell))}$\, with \,$\ell\in \{2,3\}$\,. \\
\,$\liem'$\, is an \,$\ell$-dimensional subspace of \,$\mathrm{span}\{\,\lambda_6^\sharp,\, V_{\lambda_2}(1)+V_{\lambda_4}(\sqrt{3}),\, 
V_{\lambda_2}(i)-V_{\lambda_4}(\sqrt{3}\,i)\,\}$\,. 

\item \,$\boldsymbol{(\PP,\vi=\tfrac\pi6,(\C,2))}$\,. \\
\,$\liem' = \mathrm{span}\{\,\lambda_6^\sharp,\, V_{\lambda_2}(1)+V_{\lambda_4}(\sqrt{3}),\, V_{\lambda_3}(\sqrt{3}\,i)+V_{\lambda_5}(i),\,V_{\lambda_6}(1)\,\}$\,. 

\item \,$\boldsymbol{(\Sph\times\Sph,\ell,\ell')}$\, with \,$\ell,\ell' \leq 3$\,. \\
\,$\liem' = \liea \oplus \liem_{\lambda_1}' \oplus \liem_{\lambda_6}'$\,, where \,$\liem_{\lambda_1}' \subset V_{\lambda_1}(\C)$\, and \,$\liem_{\lambda_6}' \subset V_{\lambda_6}(\C)$\,
are linear subspaces of dimension \,$\ell-1$\, resp.~\,$\ell'-1$\,. 

\item \,$\boldsymbol{(\mathrm{AI})}$\,. \\
\,$\liem' = \liea \oplus V_{\lambda_2}(\R) \oplus V_{\lambda_5}(\R) \oplus V_{\lambda_6}(i\R)$\,.

\item \,$\boldsymbol{(\mathrm{A}_2)}$\,. \\
\,$\liem' = \liea \oplus V_{\lambda_2}(\C) \oplus V_{\lambda_5}(\C) \oplus V_{\lambda_6}(\C)$\,.

\item \,$\boldsymbol{(\mathrm{G})}$\,. \\
\,$\liem' = \liea \oplus V_{\lambda_1}(\R) \oplus V_{\lambda_2}(\R) \oplus V_{\lambda_3}(i\R) \oplus V_{\lambda_4}(\R) \oplus V_{\lambda_5}(i\R) \oplus V_{\lambda_6}(\R)$\,.
\end{itemize}
We call the full name \,$(\mathrm{Geo},\vi=t)$\, etc.~given in the above table the \emph{type} of the Lie triple systems which are congruent under the
adjoint action to the space given in that entry.%
\footnote{Notice that in this case, the types \,$(\Sph\times\Sph,\ell,\ell')$\, and \,$(\Sph\times\Sph,\ell',\ell)$\, with
\,$\ell\neq\ell'$\, are not equivalent, because the two irreducible components of the Lie triple systems of this type
correspond to spheres of different radius, see Section~\ref{SSe:G2:tgsub}.}
Then no Lie triple system is of more than one type.

The Lie triple systems \,$\liem'$\, of the various types have the properties given in the following table. 
The column ``isometry type'' again gives the isometry type of the totally geodesic submanifolds corresponding to the Lie triple systems of the respective type
in abbreviated form, for the details see Section~\ref{SSe:G2:tgsub}.
\begin{center}
\begin{longtable}{|c|c|c|c|c|c|}
\hline
type of \,$\liem'$ & $\dim(\liem')$ & $\rk(\liem')$ & $\liem'$\, Lie subalgebra & \,$\liem'$\, maximal & isometry type \\
\hline
\endhead
\hline
\endfoot
\,$(\mathrm{Geo},\vi=t)$\, & $1$ & $1$ & yes & no & \,$\R$\, or \,$\Sph^1$\, \\
\,$(\Sph,\vi=0,\ell)$\, & $\ell$ & $1$ & for \,$\ell=3$\, & no & \,$\Sph^2_{r=1}$\, \\
\,$(\Sph,\vi=\arctan(\tfrac{1}{3\sqrt{3}}),\ell)$ & $\ell$ & $1$ & for \,$\ell=3$\, & for \,$\ell=3$\, & $\Sph^\ell_{r=\tfrac23\,\sqrt{21}}$ \\
\,$(\Sph,\vi=\tfrac\pi6,\ell)$ & $\ell$ & $1$ & for \,$\ell=3$\, & no & \,$\Sph^\ell_{r=1/\sqrt{3}}$\, \\
\,$(\PP,\vi=\tfrac\pi6,(\R,\ell))$ & $\ell$ & $1$ & no & no & \,$\RP^\ell$\, \\
\,$(\PP,\vi=\tfrac\pi6,(\C,2))$ & $4$ & $1$ & no & no & \,$\CP^2$\, \\
\hline
\,$(\Sph\times\Sph,\ell,\ell')$\, & \,$\ell+\ell'$\, & $2$ & for \,$\ell,\ell'\in\{1,3\}$\, & for \,$\ell=\ell'=3$\, & \,$(\Sph^\ell_{r=1} \times \Sph^{\ell'}_{r=1/\sqrt{3}})/\Z_2$\, \\
\,$(\mathrm{AI})$\, & $5$ & $2$ & no & no & \,$\SU(3)/\SO(3)$\, \\
\,$(\mathrm{A}_2)$\, & $8$ & $2$ & yes & yes & \,$\SU(3)$\, \\
\,$(\mathrm{G})$\, & $8$ & $2$ & no & yes & \,$G_2/\SO(4)$\, \\
\end{longtable}
\end{center}
\end{Theorem}

\begin{Remark}
\label{R:G2:CN}
The maximal totally geodesic submanifolds of \,$G_2$\, of type \,$(\Sph,\vi=\arctan(\tfrac{1}{3\sqrt{3}}),3)$\,, which are isometric
to a 3-sphere of radius \,$\tfrac23\,\sqrt{21}$\,, are missing from the classification by \textsc{Chen} and \textsc{Nagano} 
in Table~VIII of \cite{Chen/Nagano:totges2-1978}. Once again, these submanifolds are in a ``skew'' position in the ambient manifold \,$G_2$\, 
in the sense that their geodesic diameter \,$\tfrac23\,\sqrt{21}\,\pi$\, is strictly larger than the geodesic diameter \,$\tfrac{2}{3}\,\sqrt{3}\,\pi$\,
of \,$G_2$\,.
\end{Remark}

\emph{Proof of Theorem~\ref{T:G2:cla}.}
Once again, it is easily checked that the spaces defined in the theorem are Lie triple systems, and thus the spaces which are conjugate to one of them under
the adjoint action also are. It is also easily seen that the information in the table on the dimension and the rank of the Lie triple systems and regarding
the question which of them are Lie subalgebras of \,$\lieg_2$\, is correct. The information on the isometry type of the totally geodesic submanifolds
corresponding to the various types of Lie triple systems will be proved in Section~\ref{SSe:G2:tgsub}.

We next show that the information on the maximality of the Lie triple systems is correct. 
For this purpose, we presume that the list of Lie triple systems given in the theorem is in fact complete; this will
be proved in the remainder of the present section.

That the Lie triple systems which are claimed to be maximal in the table indeed are: This is clear for the types \,$(\mathrm{A}_2)$\, and \,$(\mathrm{G})$\,
because there are no Lie triple systems of greater dimension. For the type \,$(\Sph,3,3)$\,, we note that if it were not maximal, it could only be included
in a Lie triple system of type \,$(\mathrm{A}_2)$\, or \,$(\mathrm{G})$\, for dimension reasons. However, the Lie triple systems of type \,$(\Sph,3,3)$\,
have two orthogonal roots of multiplicity \,$2$\,, whereas the systems of type \,$(\mathrm{A}_2)$\, do not have a pair of orthogonal roots, and in the
systems of type \,$(\mathrm{G})$\,, all roots have multiplicity \,$1$\,. So such an inclusion is not in fact possible, and hence the Lie triple systems of
type \,$(\Sph,3,3)$\, are maximal. 
For the type \,$(\PP,\vi=\arctan(\tfrac{1}{3\sqrt{3}}),3)$\,: Let \,$\wh{\liem}' \subset \lieg_2$\, be a Lie triple system
with \,$\liem' \subsetneq \wh{\liem}'$\,. If \,$\wh{\liem}'$\, were of rank \,$1$\,, then it would need to have the same isotropy angle \,$\vi=\arctan(\tfrac{1}{3\sqrt{3}})$\,
and a strictly greater dimension than \,$\liem'$\,, but no such Lie triple system exists. So \,$\wh{\liem}'$\, is of rank \,$2$\,.
It now follows from the description of the type \,$(\PP,\vi=\arctan(\tfrac{1}{3\sqrt{3}}),3)$\, that \,$\wh{\liem}'$\, has two roots
at an angle of \,$\tfrac{5\pi}{6}$\, to each other and these two roots both have multiplicity \,$2$\,. Therefore \,$\wh{\liem}' = \lieg_2$\, holds,
and hence \,$\liem'$\, is maximal.

That no Lie triple systems are maximal besides those mentioned above follows from the following table:

\begin{center}
\begin{longtable}{|c|c|}
\hline
Every Lie triple system of type ... & is contained in a Lie triple system of type ... \\
\hline
\endhead
\hline
\endfoot
$(\mathrm{Geo},\vi=t)$ & $(\Sph\times\Sph,1,1)$ \\
$(\Sph,\vi=0,\ell)$ & $(\Sph\times\Sph,\ell,1)$ \\
$(\Sph,\vi=\arctan(\tfrac{1}{3\sqrt{3}}),2)$ & $(\mathrm{G})$ \\
$(\Sph,\vi=\tfrac\pi6,\ell)$ & $(\Sph\times\Sph,1,\ell)$ \\
$(\PP,\vi=\tfrac\pi6,(\R,\ell))$ & $(\Sph\times\Sph,\ell,\ell)$ \\
$(\PP,\vi=\tfrac\pi6,(\C,2))$ & $(\mathrm{G})$ \\
$(\Sph\times\Sph,\ell,\ell')$\, with \,$(\ell,\ell') \neq (3,3)$\, & \,$(\Sph\times\Sph,3,3)$\, \\
$(\mathrm{AI})$ & $(\mathrm{A}_2)$ \\
\end{longtable}
\end{center}

We now turn to the proof that the list of Lie triple systems of \,$\lieg_2$\, given in Theorem~\ref{T:G2:cla} is indeed
complete. For this purpose, we let an arbitrary Lie triple system \,$\liem'$\, of \,$\lieg_2$\,, \,$\{0\} \neq \liem' \subsetneq \lieg_2$\,,
be given. Because the Lie algebra \,$\lieg_2$\, is of rank \,$2$\,, the rank of \,$\liem'$\, is either \,$1$\, or \,$2$\,. We will handle
these two cases separately in the sequel.

We first suppose that \,$\liem'$\, is a Lie triple system of rank \,$2$\,. Let us fix a Cartan subalgebra \,$\liea$\, of \,$\liem'$\,;
because of \,$\rk(\liem') = \rk(\lieg_2)$\,, \,$\liea$\, is then also a Cartan subalgebra of \,$\lieg_2$\,. 
In relation to this situation, we again use
the notations introduced in Sections~\ref{Se:generallts} and \ref{SSe:G2:geometry}. In particular,
we consider the positive root system
\,$\Delta_+ := \{\lambda_1,\dotsc,\lambda_6\}$\, of the root system \,$\Delta := \Delta(\lieg_2,\liea)$\, of \,$\lieg_2$\,,
and also the root system \,$\Delta' := \Delta(\liem',\liea)$\, of \,$\liem'$\,. By Proposition~\ref{P:cla:subroots:subroots-neu}(b), \,$\Delta'$\,
is a root subsystem of \,$\Delta$\,, and therefore \,$\Delta_+' := \Delta' \cap \Delta_+$\, is a positive system of roots for \,$\Delta'$\,.
Moreover, in the root space decompositions of \,$\lieg_2$\, and \,$\liem'$\,
\begin{equation}
\label{eq:G:rk2:decomp}
\lieg_2 = \liea \;\oplus\; \bigoplus_{k=1}^6 V_{\lambda_k}(\C) \qmq{and}
\liem' = \liea \;\oplus\; \bigoplus_{\lambda\in\Delta_+'} \liem_\lambda'
\end{equation}
the root space \,$\liem_{\lambda}'$\, of \,$\liem'$\, with respect to \,$\lambda \in \Delta_+'$\, 
is related to the corresponding root space \,$V_{\lambda}(\C)$\, of \,$\lieg_2$\,
by \,$\liem_\lambda' = V_{\lambda}(\C) \cap \liem'$\,. 

Because the subset \,$\Delta'$\, of \,$\Delta$\, is invariant under its own Weyl group, we have (up to Weyl transformation) the following possibilities
for \,$\Delta_+'$\,, which we will treat individually in the sequel:
\begin{gather*}
\Delta_+' = \Delta_+,\quad \Delta_+' = \{\lambda_2,\lambda_5,\lambda_6\},\quad \Delta_+' = \{\lambda_1,\lambda_3,\lambda_4\}, \\
\Delta_+' = \{\lambda_1,\lambda_6\}, \quad \Delta_+' = \{\lambda_6\}, \quad \Delta_+' = \{\lambda_1\} \qmq{and} \Delta_+' = \varnothing \; . 
\end{gather*}

\textbf{The case \,$\boldsymbol{\Delta_+' = \Delta}$\,.}
In this case the Dynkin diagram with multiplicities of \,$\liem'$\, is $\xymatrix@=.4cm{ \bullet^{n_1} \ar@3{<-}[r] & \bullet^{n_2} }$ with
\,$n_1,n_2 \in \{1,2\}$\,. From the classification of the irreducible Riemannian symmetric spaces (see for example \cite{Loos:1969-2}, p.~119, 146),
we see that \,$n_1=n_2 =: n \in \{1,2\}$\, holds. If \,$n=2$\, holds, we have \,$\liem' = \lieg_2$\,. If \,$n=1$\, holds, we may by virtue of 
Proposition~\ref{P:G2:geometry:isotropy} suppose without loss of generality that \,$\liem_{\lambda_1}' = V_{\lambda_1}(\R)$\, and \,$\liem_{\lambda_2}' = V_{\lambda_2}(\R)$\, holds.
Then we can calculate the remaining root spaces of \,$\liem'$\, one by one: 
We have \,$R(\lambda_1^\sharp,V_{\lambda_2}(1))V_{\lambda_1}(1) = \tfrac{3\sqrt{3}}{4}\cdot V_{\lambda_3}(i)$\, and 
therefore \,$\liem_{\lambda_3}' = V_{\lambda_3}(i\R)$\,. 
We have \,$R(\lambda_1^\sharp,V_{\lambda_3}(i))V_{\lambda_1}(1) = \tfrac{\sqrt{3}}{4}\cdot V_{\lambda_2}(1) - \tfrac12\cdot V_{\lambda_4}(1)$\, and
therefore \,$\liem_{\lambda_4}' = V_{\lambda_4}(\R)$\,. 
We have \,$R(\lambda_1^\sharp,V_{\lambda_4}(1))V_{\lambda_1}(1) = \tfrac12\cdot V_{\lambda_3}(i) - \tfrac{\sqrt{3}}{4}\cdot V_{\lambda_5}(i)$\, and
therefore \,$\liem_{\lambda_5}' = V_{\lambda_5}(i\R)$\,. 
Finally, we have \,$R(\lambda_1^\sharp,V_{\lambda_5}(i))V_{\lambda_2}(1) = \tfrac{3\sqrt{3}}{4} \cdot V_{\lambda_6}(1)$\, and
therefore \,$\liem_{\lambda_6}' = V_{\lambda_6}(\R)$\,. Thus it follows from Equation~\eqref{eq:G:rk2:decomp} that
$$ \liem' = \liea \oplus V_{\lambda_1}(\R) \oplus V_{\lambda_2}(\R) \oplus V_{\lambda_3}(i\R) \oplus V_{\lambda_4}(\R) \oplus V_{\lambda_5}(i\R) \oplus V_{\lambda_6}(\R) $$
holds, and therefore \,$\liem'$\, is of type \,$(\mathrm{G})$\,.

\textbf{The case \,$\boldsymbol{\Delta_+' = \{\lambda_2,\lambda_5,\lambda_6\}}$\,.} 
In this case, the Dynkin diagram with multiplicities of \,$\liem'$\, is $\xymatrix@=.4cm{ \bullet^{n} \ar@{-}[r] & \bullet^{n} }$ with \,$n \in \{1,2\}$\,.
In the case \,$n=2$\, we have \,$\liem' = \liea \oplus V_{\lambda_2}(\C) \oplus V_{\lambda_5}(\C) \oplus V_{\lambda_6}(\C)$\, 
and therefore \,$\liem'$\, is of type \,$(\mathrm{A}_2)$\,. In the case \,$n=1$\, we may suppose without loss of generality \,$\liem_{\lambda_2}' = V_{\lambda_2}(\R)$\,
and \,$\liem_{\lambda_5}' = V_{\lambda_5}(\R)$\,; then we have \,$\liem' \in R(\lambda_1^\sharp,V_{\lambda_2}(1))V_{\lambda_5}(1) = -\tfrac{3\sqrt{3}}{4}\,V_{\lambda_6}(i)$\,
and hence \,$\liem_{\lambda_6}' = V_{\lambda_6}(i\R)$\,. Therefore we have \,$\liem' = \liea \oplus V_{\lambda_2}(\R) \oplus V_{\lambda_5}(\R) \oplus V_{\lambda_6}(i\R)$\,,
thus \,$\liem'$\, is of type \,$(\mathrm{AI})$\,.

\textbf{The case \,$\boldsymbol{\Delta_+' = \{\lambda_1,\lambda_3,\lambda_4\}}$\,.}
Assume that \,$\liem'$\, is a Lie triple system with this root system. Then there exist \,$c,d\in\C^\times$\, so that \,$V_{\lambda_1}(c), V_{\lambda_3}(d) \in \liem'$\,
holds. We have \,$\liem' \ni R(\lambda_1^\sharp, V_{\lambda_1}(c))V_{\lambda_3}(d) = \tfrac{\sqrt{3}}{2}\,V_{\lambda_2}(\overline{c}\,d\,i) + V_{\lambda_4}(c\,d\,i)$\,
and therefore in particular \,$\lambda_2 \in \Delta_+'$\,, contrary to the hypothesis \,$\Delta_+' = \{\lambda_1,\lambda_3,\lambda_4\}$\,. This calculation shows that 
there do not exist any Lie triple systems \,$\liem'$\, of \,$\lieg_2$\, with \,$\Delta_+' = \{\lambda_1,\lambda_3,\lambda_4\}$\,.

\textbf{The cases \,$\boldsymbol{\Delta_+' \subset \{\lambda_1,\lambda_6\}}$\,.}
In this case we have \,$\liem' = \liea \oplus \liem_{\lambda_1}' \oplus \liem_{\lambda_6}'$\, by Equation~\eqref{eq:G:rk2:decomp}, therefore \,$\liem'$\,
is of type \,$(\Sph\times\Sph,\ell,\ell')$\, with  \,$\ell := 1+\dim(\liem_{\lambda_1}')$\, and \,$\ell' := 1+\dim(\liem_{\lambda_6}')$\,.

This completes the treatment of the case where \,$\liem'$\, is of rank \,$2$\,. 

\bigskip

We now suppose that \,$\liem' \subset \lieg_2$\, is a Lie triple system of rank \,$1$\,. We may suppose without loss of generality that \,$\liem'$\,
contains a unit vector \,$H$\, from the closure 
of the positive Weyl chamber \,$\liec$\,.
By Equations~\eqref{eq:G2:isotropy:liec} and \eqref{eq:G2:isotropy:vivt}, we then have with \,$\vi_0 := \vi(H) \in [0,\tfrac\pi6]$\,
\begin{equation}
\label{eq:G2:rk1:H}
H = \cos(\vi_0)\,\lambda_4^\sharp + \sin(\vi_0)\,\tfrac{1}{\sqrt{3}}\,\lambda_2^\sharp \; .
\end{equation}
Because of \,$\rk(\liem')=1$\,, \,$\liea' := \R\,H$\, is a Cartan subalgebra of \,$\liem'$\,, and we have \,$\liea' = \liea \cap \liem'$\,. 
It follows from Proposition~\ref{P:cla:subroots:subroots-neu}(a) that the root systems \,$\Delta'$\,
and \,$\Delta$\, of \,$\liem'$\, resp.~\,$\lieg_2$\, with respect to \,$\liea'$\, resp.~to \,$\liea$\, are related by
\begin{equation}
\label{eq:G2:rk1:Delta'Delta}
\Delta' \;\subset\; \Mengegr{\lambda(H)\,\alpha_0}{\lambda\in\Delta,\,\lambda(H)\neq 0}
\end{equation}
with the linear form \,$\alpha_0: \liea' \to \R,\; tH\mapsto t$\,; moreover for \,$\liem'$\, we have the root space decomposition
\begin{equation}
\label{eq:G2:rk1:m'decomp}
\liem' = \liea' \oplus \bigoplus_{\alpha\in\Delta_+'} \liem_\alpha' 
\end{equation}
where for any root \,$\alpha\in\Delta'$\,, the corresponding root space \,$\liem_\alpha'$\, is given by
\begin{equation}
\label{eq:G2:rk1:malpha'}
\liem_\alpha' = \left( \bigoplus_{\substack{\lambda \in \Delta \\ \lambda(H) = \alpha(H)}} V_{\lambda}(\C) \right) \;\cap\; \liem' \; . 
\end{equation}

If \,$\Delta' = \varnothing$\, holds, we have \,$\liem' = \R H$\,, and therefore \,$\liem'$\, is then of type \,$(\mathrm{Geo},\vi=\vi_0)$\,. Otherwise
it follows from Proposition~\ref{P:cla:subroots:Comp} that one of the following two conditions holds: Either \,$H$\, is proportional to a root vector \,$\lambda^\sharp$\,
with \,$\lambda\in\Delta$\,, or there exist two \,$\lambda,\mu\in\Delta$\, (\,$\lambda\neq\mu$\,) so that \,$H$\, is orthogonal to \,$\lambda^\sharp-\mu^\sharp$\,.
Evaluating all possible values for \,$\lambda$\, and \,$\mu$\,, we see that \,$\vi_0\in\{0,\arctan(\tfrac{1}{3\sqrt{3}}),\tfrac\pi6\}$\, holds.

In the sequel we consider the three possible values for \,$\vi_0$\, individually.

\textbf{The case \,$\boldsymbol{\vi_0 = 0}$\,.}
In this case we have \,$H = \lambda_4^\sharp$\, by Equation~\eqref{eq:G2:rk1:H} and therefore 
$$ \lambda_1(H) = \tfrac12,\quad \lambda_2(H) = 0,\quad \lambda_3(H) = \tfrac12,\quad 
\lambda_4(H) = 1,\quad \lambda_5(H) = \tfrac32,\quad \lambda_6(H) = \tfrac32 \; . $$
Thus we have \,$\Delta' \subset \{\pm \alpha,\pm 2\alpha,\pm 3\alpha\}$\, with \,$\alpha := \lambda_1|\liea' = \lambda_3|\liea'$\, by Equation~\eqref{eq:G2:rk1:Delta'Delta},
\,$\liem' = \R\,H \oplus \liem_{\alpha}' \oplus \liem_{2\alpha}' \oplus \liem_{3\alpha}'$\, by Equation~\eqref{eq:G2:rk1:m'decomp} and
\,$\liem_\alpha' \subset V_{\lambda_1}(\C) \oplus V_{\lambda_3}(\C)$\,, \,$\liem_{2\alpha}' \subset V_{\lambda_4}(\C)$\, and \,$\liem_{3\alpha}' \subset V_{\lambda_5}(\C)
\oplus V_{\lambda_6}(\C)$\, by Equation~\eqref{eq:G2:rk1:malpha'}.

We now show that actually \,$\alpha,3\alpha\not\in\Delta_+'$\, holds.

Indeed, let \,$v\in\liem_{\alpha}'$\, be given. Then there exist \,$c,d\in\C$\, so that \,$v = V_{\lambda_1}(c) + V_{\lambda_3}(d)$\, holds. We have
\,$|c|=|d|$\, because of Proposition~\ref{P:cla:skew} and the fact that \,$\alpha^\sharp = \tfrac12\,\lambda_1^\sharp + \tfrac12\,\lambda_3^\sharp$\, holds. Next we notice
that because of \,$\lambda_2(H)=0$\,, the \,$V_{\lambda_2}(\C)$-component of every vector in \,$\liem'$\, must vanish. However, 
the \,$V_{\lambda_2}(\C)$-component of \,$R(H,v)v \in \liem'$\, equals \,$\tfrac{\sqrt{3}}{2}\,V_{\lambda_2}(\overline{c}\,d\,i)$\,, and so we conclude
\,$\overline{c}\,d = 0$\,. Because of \,$|c|=|d|$\,, it follows that we have \,$c=d=0$\, and hence \,$v=0$\,. Thus we have \,$\liem_{\alpha}' = \{0\}$\,
and hence \,$\alpha\not\in\Delta'$\,.

A similar calculation also shows \,$3\alpha \not\in\Delta'$\,, and therefore we have \,$\Delta' = \{\pm 2\alpha\}$\, and hence \,$\liem' = \R\,H \oplus
\liem_{2\alpha}'$\, with a linear subspace \,$\{0\} \neq \liem_{2\alpha}' \subset V_{\lambda_4}(\C)$\,. It follows that \,$\liem'$\, is of type
\,$(\Sph,\vi=0,\ell)$\, with \,$\ell := 1+ n_{2\alpha}'$\,. 

\textbf{The case \,$\boldsymbol{\vi_0 = \arctan(\tfrac{1}{3\sqrt{3}})}$\,.}
In this case we have \,$H = \tfrac{\sqrt{21}}{42}\,(9\lambda_4^\sharp + \lambda_2^\sharp)$\, 
by Equation~\eqref{eq:G2:rk1:H} and therefore 
\begin{gather*}
\lambda_1(H) = \tfrac{\sqrt{21}}{14}\cdot 1,\quad
\lambda_2(H) = \tfrac{\sqrt{21}}{14}\cdot 1,\quad
\lambda_3(H) = \tfrac{\sqrt{21}}{14}\cdot 2, \\ 
\lambda_4(H) = \tfrac{\sqrt{21}}{14}\cdot 3,\quad
\lambda_5(H) = \tfrac{\sqrt{21}}{14}\cdot 4,\quad
\lambda_6(H) = \tfrac{\sqrt{21}}{14}\cdot 5 \; . 
\end{gather*}
In the present case we have \,$\lambda^\sharp\not\in\liea'$\, for every \,$\lambda\in\Delta$\,, therefore \,$\liem'$\, can only have composite
roots (see Definition~\ref{D:cla:subroots:Elemcomp}) by Proposition~\ref{P:cla:subroots:Comp}(a). This fact, together with the above values of \,$\lambda(H)$\,
and Equation~\eqref{eq:G2:rk1:Delta'Delta}, shows that we have \,$\Delta' = \{\pm \alpha\}$\, with  \,$\alpha := \lambda_1|\liea' = \lambda_2|\liea'$\,,
Moreover, we have \,$\liem' = \R\,H \oplus \liem_\alpha'$\, by Equation~\eqref{eq:G2:rk1:m'decomp} and
\,$\liem_\alpha' \subset V_{\lambda_1}(\C) \oplus V_{\lambda_2}(\C)$\, by Equation~\eqref{eq:G2:rk1:malpha'}.

Let \,$v\in\liem_\alpha'$\, be given, say \,$v = V_{\lambda_1}(c_1) + V_{\lambda_2}(c_2)$\, with \,$c_1,c_2 \in \C$\,. We have
\,$\alpha^\sharp = \alpha(H)\cdot H = \lambda_1(H)\cdot H = \tfrac{\sqrt{21}}{14}\cdot H = \tfrac{1}{28}(9\lambda_4^\sharp + \lambda_2^\sharp)
= \tfrac{9}{14}\,\lambda_1^\sharp + \tfrac{5}{14}\,\lambda_2^\sharp$\,, and therefore Proposition~\ref{P:cla:skew} shows that we have
\,$|c_2| = \sqrt{\tfrac{5/14}{9/14}}\,|c_1| = \tfrac13\,\sqrt{5}\,|c_1|$\,. 

By Proposition~\ref{P:G2:geometry:isotropy} we may therefore suppose without loss of generality that \,$v_0 := V_{\lambda_1}(1) + V_{\lambda_2}(\tfrac13\,\sqrt{5})
\in \liem_{\alpha}'$\, holds. Then we have
\begin{equation}
\label{eq:G2:rk1:arctan33:Rv0vH}
\liem' \ni R(v_0,v)H = -\tfrac{\sqrt{7}}{14}\,V_{\lambda_3}(i(3c_2-\sqrt{5}c_1)) \; .
\end{equation}
Because of \,$\lambda_3|\liea' = 2\alpha \not\in\Delta'$\,, the vector \eqref{eq:G2:rk1:arctan33:Rv0vH} must vanish, and thus we have \,$3c_2 = \sqrt{5}\,c_1$\,. 
This shows that \,$\liem_\alpha'$\, is a linear subspace of \,$\Menge{V_{\lambda_1}(c) + V_{\lambda_2}(\tfrac{\sqrt{5}}{3}\,c)}{c\in\C}$\,, and therefore
\,$\liem' = \R\,H \oplus \liem_{\alpha}'$\, is of type \,$(\PP,\vi=\arctan(\tfrac{1}{3\sqrt{3}}),\ell)$\, with \,$\ell := 1+n_\alpha'$\,. 

\textbf{The case \,$\boldsymbol{\vi_0 = \tfrac\pi6}$\,.}
In this case we have \,$H = \tfrac{1}{\sqrt{3}}\,\lambda_6^\sharp$\, by Equation~\eqref{eq:G2:rk1:H} and therefore
$$ \lambda_1(H) = 0,\quad \lambda_2(H) = \tfrac{\sqrt{3}}{2},\quad \lambda_3(H) = \tfrac{\sqrt{3}}{2},\quad 
\lambda_4(H) = \tfrac{\sqrt{3}}{2},\quad \lambda_5(H) = \tfrac{\sqrt{3}}{2},\quad \lambda_6(H) = \sqrt{3} \; . $$
Thus we have \,$\Delta' \subset \{\pm \alpha,\pm 2\alpha\}$\, with \,$\alpha := \lambda_2|\liea' = \lambda_3|\liea' = \lambda_4|\liea' = \lambda_5|\liea'$\, 
by Equation~\eqref{eq:G2:rk1:Delta'Delta},
\,$\liem' = \R\,H \oplus \liem_{\alpha}' \oplus \liem_{2\alpha}'$\, by Equation~\eqref{eq:G2:rk1:m'decomp} and
\,$\liem_\alpha' \subset V_{\lambda_2}(\C) \oplus V_{\lambda_3}(\C) \oplus V_{\lambda_4}(\C) \oplus V_{\lambda_5}(\C)$\, and \,$\liem_{2\alpha}' \subset V_{\lambda_6}(\C)$\,
by Equation~\eqref{eq:G2:rk1:malpha'}.

Let us first consider the case \,$\alpha\not\in\Delta'$\, and therefore \,$\Delta' = \{\pm 2\alpha\}$\,. Then we have \,$\liem' = \R\,H \oplus \liem_{2\alpha}'$\,,
and therefore \,$\liem'$\, then is of type \,$(\Sph,\vi=\tfrac\pi6,\ell)$\, with \,$\ell := 1+n_{2\alpha}'$\,. 

So we now suppose \,$\alpha\in\Delta'$\,. Then we have for any \,$v\in \liem_\alpha'$\,, say \,$v = \sum_{k=2}^5 V_{\lambda_k}(c_k)$\, with \,$c_2,\dotsc,c_5 \in \C$\, 
\begin{equation}
\label{eq:G2:rk1:RHvv}
\sqrt{3}\cdot R(H,v)v 
= \bigr( \tfrac32\,|c_3|^2 + 3\,|c_4|^2 + \tfrac92\,|c_5|^2\bigr)\,\lambda_1^\sharp
+ \tfrac32\,\|c\|^2\,\lambda_2^\sharp
+ V_{\lambda_1}\bigr(\,-\tfrac{3\sqrt{3}}{2}\,i\,(\overline{c_2}\,c_3 + \overline{c_4}\,c_5) - 3i\,\overline{c_3}\,c_4 \,\bigr) \; .
\end{equation}
Because this vector is a member of \,$\liem'$\,, its \,$\liea$-component must be proportional to \,$H = \tfrac{1}{\sqrt{3}}\,(3\lambda_1^\sharp + 2\lambda_2^\sharp)$\,,
and therefore we have
$$ 2 \cdot \bigr( \tfrac32\,|c_3|^2 + 3\,|c_4|^2 + \tfrac92\,|c_5|^2\bigr) = 3 \cdot \tfrac32\,\|c\|^2 \; , $$
hence 
\begin{equation}
\label{eq:G2:rk1:lengthrel}
3\,|c_2|^2 + |c_3|^2 = |c_4|^2 + 3\,|c_5|^2 \; .
\end{equation}
As a consequence of this equation we have \,$n_{\alpha}' \leq 4$\,. Moreover, because we have \,$\lambda_1(H)=0$\,, the \,$V_{\lambda_1}(\C)$-component
of the vector \eqref{eq:G2:rk1:RHvv} must vanish, and thus we have 
\begin{equation}
\label{eq:G2:rk1:c2345rel}
\tfrac{\sqrt{3}}{2}\,\bigr(\,\overline{c_2}\,c_3 + \overline{c_4}\,c_5 \,\bigr) = -\,\overline{c_3}\,c_4 \; . 
\end{equation}
It is a consequence of Equations~\eqref{eq:G2:rk1:lengthrel} and \eqref{eq:G2:rk1:c2345rel} that by application of the adjoint action of the subgroup
of \,$G_2$\, with the Lie algebra \,$\liea \oplus V_{\lambda_1}(\C)$\,, we can arrange that \,$v_0 := V_{\lambda_2}(1) + V_{\lambda_4}(\sqrt{3}) \in \liem_\alpha'$\, 
holds. 

Then we have for any \,$v \in \liem_\alpha'$\, as above
\begin{equation}
\label{eq:G2:rk1:RHvv0}
R(H,v)v_0
= 3\,\RE(c_4)\,\lambda_1^\sharp + \bigr(\,\tfrac{\sqrt{3}}{2}\,\RE(c_2) + \tfrac32\,\RE(c_4)\,\bigr)\lambda_2^\sharp
- V_{\lambda_1}\bigr(\,\tfrac34\,i\,c_3 + \tfrac32\,i\,\overline{c_3} + \tfrac{3\sqrt{3}}{4}\,i\,c_5\,\bigr)
- V_{\lambda_6}\bigr(\,\tfrac{3\sqrt{3}}{4}\,i\,c_3 + \tfrac34\,i\,c_5\,\bigr) \; . 
\end{equation}
Because this vector is again a member of \,$\liem'$\,, its \,$\liea$-component must be proportional to \,$H$\,, and thus we have
$$ 2\cdot \bigr(3\,\RE(c_4)) = 3\cdot \bigr(\,\tfrac{\sqrt{3}}{2}\,\RE(c_2) + \tfrac32\,\RE(c_4)\,\bigr) \;, $$
hence
\begin{equation}
\label{eq:G2:rk1:vRe24}
\RE(c_4) = \sqrt{3}\,\RE(c_2) \; .
\end{equation}
Moreover, the \,$V_{\lambda_1}(\C)$-component of \eqref{eq:G2:rk1:RHvv0} vanishes, and thus we have
$$ \tfrac34\,i\,c_3 + \tfrac32\,i\,\overline{c_3} + \tfrac{3\sqrt{3}}{4}\,i\,c_5 = 0 \;, $$
hence
\begin{equation}
\label{eq:G2:rk1:v35}
c_3 + 2\,\overline{c_3} = -\sqrt{3}\,c_5 \; .
\end{equation}

Let us now first consider the case \,$2\alpha\not\in\Delta'$\, and thus \,$\liem' = \R\,H \oplus \liem_\alpha'$\,. If \,$n_\alpha'=1$\, holds, we have
\,$\liem' = \R\,H \oplus \R\,v_0$\,, and therefore \,$\liem'$\, is of type \,$(\PP,\vi=\tfrac\pi6,(\R,2))$\,. Otherwise, we let \,$v\in\liem_{\alpha}'$\, be given
as above, and suppose that \,$v$\, is orthogonal to \,$v_0$\,. The latter condition, together with Equation~\eqref{eq:G2:rk1:vRe24}, shows that
\begin{equation}
\label{eq:G2:rk1:alpha:Re24}
\RE(c_2) = \RE(c_4) = 0
\end{equation}
holds. Because of \,$2\alpha\not\in\Delta'$\,, also the \,$V_{\lambda_6}(\C)$-component of the vector \eqref{eq:G2:rk1:RHvv0} vanishes in this case, 
and therefore we have
$$ \tfrac{3\sqrt{3}}{4}\,i\,c_3 + \tfrac34\,i\,c_5 = 0 \;, $$
hence
\begin{equation}
\label{eq:G2:rk1:alpha:c35}
c_5 = -\sqrt{3}\,c_3 \; . 
\end{equation}
Now we obtain from Equations~\eqref{eq:G2:rk1:v35} and \eqref{eq:G2:rk1:alpha:c35}: \,$c_3 + 2\,\overline{c_3} = -\sqrt{3}\,c_5 = 3\,c_3$\,, hence
\,$c_3 = \overline{c_3}$\, and thus \,$c_3 \in \R$\,. Equation~\eqref{eq:G2:rk1:alpha:c35} then also implies \,$c_5\in\R$\,. 
This, together with \eqref{eq:G2:rk1:alpha:Re24}, shows that there exist \,$t_2,\dotsc,t_5\in\R$\, so that 
$$ v = V_{\lambda_2}(it_2) + V_{\lambda_3}(t_3) + V_{\lambda_4}(it_4) + V_{\lambda_5}(t_5) $$
holds. It now follows from Equation~\eqref{eq:G2:rk1:alpha:c35} that we have
\begin{equation}
\label{eq:G2:rk1:alpha:t5}
t_5 = -\sqrt{3}\,t_3 \;,
\end{equation}
and, using this equation we obtain from Equation~\eqref{eq:G2:rk1:c2345rel}
$$ t_3 \cdot \bigr(-\tfrac{\sqrt{3}}{2}\,t_2 + \tfrac52\,t_4\bigr) = 0 $$
and therefore
$$ \qmq{either} t_3 = 0 \qmq{or} \sqrt{3}\,t_2 = 5\,t_4 \; . $$
If \,$t_3=0$\, holds, then we also have \,$t_5=0$\, by Equation~\eqref{eq:G2:rk1:alpha:t5},
and Equation~\eqref{eq:G2:rk1:lengthrel} shows that we have \,$t_4 = \pm \sqrt{3}\,t_2$\,. If \,$t_3\neq 0$\, holds,
we have \,$t_4 = \tfrac15\,\sqrt{3}\,t_2$\,, and therefore by Equation~\eqref{eq:G2:rk1:lengthrel} \,$t_3 = \pm\,\tfrac35\,t_2$\,,
hence \,$t_5 = \mp \tfrac{3\sqrt{3}}{5}\,t_2$\, by Equation~\eqref{eq:G2:rk1:alpha:t5}. 

This consideration shows that we have \,$\liem' = \R\,H \oplus \liem_\alpha'$\,, where either of the following two equations
holds:
\begin{align}
\label{eq:G2:rk1:alpha:malpha1}
\liem_\alpha' & = \R\,v_0 \oplus \R\,(V_{\lambda_2}(i)+\eps\,V_{\lambda_4}(\sqrt{3}\,i)) \\
\label{eq:G2:rk1:alpha:malpha2}
\text{or}\quad \liem_\alpha' & = \R\,v_0 \oplus \R\,(V_{\lambda_2}(i) + \eps\,V_{\lambda_3}(\tfrac35) + V_{\lambda_4}(\tfrac{\sqrt{3}}{5}i)
- \eps\,V_{\lambda_5}(\tfrac{3\sqrt{3}}{5}))
\end{align}
with \,$\eps\in\{\pm 1\}$\,. 

In either case \,$\liem'$\, is of type \,$(\PP,\vi=\tfrac\pi6,(\R,3))$\,: If \,$\liem_{\alpha}'$\, is given by 
Equation~\eqref{eq:G2:rk1:alpha:malpha1}, this is obvious. On the other hand, if \,$\liem_{\alpha}'$\, is given by 
Equation~\eqref{eq:G2:rk1:alpha:malpha2} (without loss of generality with \,$\eps=1$\,), 
we note that \,$\liem'$\, is contained in the linear space \,$\wh{\liem}'$\,
spanned by the vectors
\begin{align*}
3\lambda_1^\sharp + 2\lambda_2^\sharp, \quad 
& \lambda_1^\sharp + \tfrac43\,\lambda_2^\sharp + V_{\lambda_1}(\sqrt{3}) , \\
V_{\lambda_2}(1) + V_{\lambda_4}(\sqrt{3}) , \quad
& V_{\lambda_2}(i) + V_{\lambda_4}(\sqrt{3}\,i) , \\
V_{\lambda_2}(i) + V_{\lambda_3}(\tfrac35) + V_{\lambda_4}(\tfrac{\sqrt{3}}{5}\,i) - V_{\lambda_5}(\tfrac{3\,\sqrt{3}}{5}), \quad 
& V_{\lambda_2}(1) - V_{\lambda_3}(\tfrac35\,i) + V_{\lambda_4}(\tfrac{\sqrt{3}}{5}) + V_{\lambda_5}(\tfrac{3\,\sqrt{3}}{5}\,i) \; . 
\end{align*}
One checks that \,$\wh{\liem}'$\, is a Lie triple system of rank \,$2$\, and dimension \,$6$\,, and therefore
(by the preceding classification of the Lie triple systems of \,$\lieg_2$\, of rank \,$2$\,) of type \,$(\Sph\times\Sph,3,3)$\,.
Hence \,$\wh{\liem}'$\, is congruent under the adjoint action to the standard Lie triple system of type \,$(\Sph\times\Sph,3,3)$\,
given in the theorem. 
\,$\liem'$\, corresponds to the diagonal in the local sphere product corresponding to \,$\wh{\liem}'$\,, and is therefore
congruent under the adjoint action to the diagonal in the standard Lie triple system of type \,$(\Sph\times\Sph,3,3)$\,,
which is the standard Lie triple system of type \,$(\PP,\vi=\tfrac\pi6,(\R,3))$\,. Therefore also \,$\liem'$\, itself
is of type \,$(\PP,\vi=\tfrac\pi6,(\R,3))$\,.

Let us finally turn our attention to the case where \,$\Delta' = \{\pm \alpha,\pm 2\alpha\}$\, holds.
From the classification of Riemannian symmetric spaces of rank \,$1$\,, we then must have \,$n_{2\alpha}'=1$\,. 
Without loss of generality we suppose \,$\liem_{2\alpha}' = V_{\lambda_6}(\R)$\,. We then have besides \,$v_0 \in \liem_{\alpha}'$\,
also \,$\liem_{\alpha}' \ni R(H,V_{\lambda_6}(1))v_0 = V_{\lambda_3}(\tfrac92\,i)+ V_{\lambda_5}(\tfrac{3\,\sqrt{3}}{2}\,i)$\,
and therefore \,$v_0' := V_{\lambda_3}(\sqrt{3}\,i) + V_{\lambda_5}(i) \in \liem_{\alpha}'$\,. Thus we have
\,$\R\,v_0 \oplus \R\,v_0' \subset \liem_{\alpha}'$\,. Below, we will show that in fact \,$\liem_{\alpha}' = \R\,v_0 \oplus
\R\,v_0'$\, holds. Therefore we have \,$\liem' = \R\,H \oplus \liem_{\alpha}' \oplus \liem_{2\alpha}' 
= \R\,H \oplus \R\,v_0 \oplus \R\,v_0' \oplus \R\,V_{\lambda_6}(1)$\,, and hence \,$\liem'$\, is of type
\,$(\PP,\vi=\tfrac\pi6,(\C,2))$\,.

For the proof of \,$\liem_{\alpha}' = \R\,v_0 \oplus \R\,v_0'$\, we let \,$v\in\liem_{\alpha}'$\, be given, and suppose that
\,$v$\, is orthogonal to \,$v_0$\, and \,$v_0'$\,. Then we are to show \,$v=0$\,. We write 
\,$v = \sum_{k=2}^5 V_{\lambda_k}(c_k)$\, with \,$c_2,\dotsc,c_5 \in \C$\, as before. Then because \,$v$\, is orthogonal to \,$v_0$\,,
we have as in the treatment of the case \,$\Delta' = \{\pm \alpha\}$\, (see Equation~\eqref{eq:G2:rk1:alpha:Re24})
$$ \RE(c_2) = \RE(c_4) = 0 \;, $$
and an analogous calculation based on the facts that \,$v$\, is orthogonal to \,$v_0'$\, and that the \,$\liea$-component of 
\,$R(H,v)v_0'$\, is proportional to \,$H$\, gives us
$$ \IM(c_3) = \IM(c_5) = 0 \; . $$
These two equations show that we have \,$v = V_{\lambda_2}(it_2) + V_{\lambda_3}(t_3) + V_{\lambda_4}(it_4) + V_{\lambda_5}(t_5)$\,
with \,$t_2,\dotsc,t_5\in\R$\,. Equation~\eqref{eq:G2:rk1:v35} now gives
\begin{equation}
\label{eq:G2:rk1:alpha2alpha:t5}
t_5 = -\sqrt{3}\,t_3 
\end{equation}
and an calculation analogous to the one leading to Equation~\eqref{eq:G2:rk1:v35}, but based on the fact that the
\,$V_{\lambda_1}(\C)$-component of \,$R(H,v)v_0'$\, vanishes, gives
\begin{equation}
\label{eq:G2:rk1:alpha2alpha:t2}
t_2 = -\sqrt{3}\,t_4 \; . 
\end{equation}
By plugging Equations~\eqref{eq:G2:rk1:alpha2alpha:t5} and \eqref{eq:G2:rk1:alpha2alpha:t2} into Equation~\eqref{eq:G2:rk1:lengthrel}
we obtain \,$3\,(-\sqrt{3}\,t_4)^2 + t_3^2 = t_4^2 + 3\,(-\sqrt{3}\,t_3)^2$\, and therefore
\begin{equation}
\label{eq:G2:rk1:alpha2alpha:t4}
t_4 = \eps\,t_3
\end{equation}
with \,$\eps \in\{\pm 1\}$\,. By plugging Equations~\eqref{eq:G2:rk1:alpha2alpha:t5}, \eqref{eq:G2:rk1:alpha2alpha:t2}
and \eqref{eq:G2:rk1:alpha2alpha:t4} into Equation~\eqref{eq:G2:rk1:c2345rel} we now obtain
$$ \tfrac{\sqrt{3}}{2}\,(i\,\sqrt{3}\,\eps\,t_3\,t_3 + i\,\eps\,t_3\,\sqrt{3}\,t_3) = - t_3\,i\,\eps\,t_3 $$
and therefore \,$3\,t_3^2 = -t_3^2$\,, hence \,$t_3=0$\,. Equations ~\eqref{eq:G2:rk1:alpha2alpha:t5}, \eqref{eq:G2:rk1:alpha2alpha:t2}
and \eqref{eq:G2:rk1:alpha2alpha:t4} now imply also \,$t_2=t_4=t_5=0$\, and therefore \,$v=0$\,.

This completes the classification of the Lie triple systems in \,$\lieg_2$\,. 
\hfill $\Box$

\subsection[Totally geodesic submanifolds in \,$G_2$\,]{Totally geodesic submanifolds in \,$\boldsymbol{G_2}$\,}
\label{SSe:G2:tgsub}

Once again, we describe totally geodesic isometric embeddings for the maximal Lie triple systems of \,$G_2$\, to determine
the global isometry type of the totally geodesic submanifolds of \,$G_2$\,. 
We obtain the results of the following
table, using the same notations for the isometry types as in Section~\ref{SSe:EIII:tgsub}:
\begin{center}
\begin{tabular}{|c|c|c|}
\hline
type of Lie triple system & corresponding global isometry type & properties \\
\hline
$(\mathrm{Geo},\vi=t)$ & \,$\R$\, or \,$\Sph^1$\, & \\
$(\Sph,\vi=0,\ell)$ & \,$\Sph^\ell_{r=1}$\, & \\
$(\Sph,\vi=\arctan(\tfrac{1}{3\sqrt{3}}),\ell)$\, & $\Sph^\ell_{r=\tfrac23\sqrt{21}}$\, & \,$\ell=3$\,: maximal \\
$(\Sph,\vi=\tfrac\pi6,\ell)$ & \,$\Sph^\ell_{r=1/\sqrt{3}}$\, & \,$\ell=3$\,: Helgason sphere \\
$(\PP,\vi=\tfrac\pi6,(\K,\ell))$ & \,$\KP^\ell_{\vkap=3/4}$\, & \\
\hline
$(\Sph\times\Sph,\ell,\ell')$ & $(\Sph^\ell_{r=1} \times \Sph^{\ell'}_{r=1/\sqrt{3}})/\{\pm \id\}$\, & \,$\ell=\ell'=3$\,: meridian, maximal \\
$(\mathrm{AI})$ & \,$\SU(3)/\SO(3)_{\mathrm{srr}=\sqrt{3}}$\, & \\
$(\mathrm{A}_2)$ & \,$\SU(3)_{\mathrm{srr}=\sqrt{3}}$\, & maximal \\
$(\mathrm{G})$ & \,$G_2/\SO(4)_{\mathrm{srr}=1}$\, & polar, maximal \\
\hline
\end{tabular}
\end{center}

\paragraph{Type \,$\boldsymbol{(\textrm{G})}$\,.}
The totally geodesic embedding corresponding to this type is the Cartan embedding \,$f: G_2/\SO(4) \to G_2$\, of the
Riemannian symmetric space \,$G_2/\SO(4)$\,. 

We describe the Cartan embedding for the general situation of a Riemannian symmetric space \,$M = G/K$\,. Let \,$\sigma: G \to G$\,
be the involutive automorphism which describes the symmetric structure of \,$M$\,. Then the map
$$ f: G/K \to G,\; g\cdot K \mapsto \sigma(g)\cdot g^{-1} $$
is called the \emph{Cartan map} of \,$M$\,. Because of \,$\Fix(\sigma)_0 \subset K \subset \Fix(\sigma)$\,,
\,$f$\, is a well-defined covering map onto its image; moreover \,$f$\, turns out to be totally geodesic. If \,$M$\, is a
``bottom space'', i.e.~there exists no non-trivial symmetric covering map with total space \,$M$\,, we have
\,$K = \Fix(\sigma)$\, and therefore \,$f$\, is a totally geodesic embedding in this case. Then \,$f$\, is called the
\emph{Cartan embedding} of \,$M$\,. 

\paragraph{Type \,$\boldsymbol{(\Sph\times\Sph,\ell,\ell')}$\, and the types of rank 1 with isotropy angle \,$\boldsymbol{\in\{0,\tfrac\pi6\}}$\,.}
For the construction of these types we consider the skew-field of quaternions \,$\HH$\, and the division algebra of the 
octonions \,$\OO$\,. \,$\OO$\, can be realized as \,$\OO = \HH \oplus \HH$\,, where the octonion multiplication is for any
\,$x,y \in \OO$\,, say \,$x=(x_1,x_2)$\, and \,$y=(y_1,y_2)$\, with \,$x_i,y_i \in \HH$\,, given by the equation
$$ x\cdot y = (x_1\,y_1 - \overline{y_2}\,x_2 \,,\, x_2\,\overline{y_1} + y_2\,x_1) \; . $$
In this setting, the symplectic group \,$\Sp(1)$\, is realized as the space of unit quaternions
with the quaternion multiplication as group action (hence \,$\Sp(1)$\, is isometric to a \,$3$-sphere),
and the Lie group \,$G_2$\, is realized as the automorphism group of \,$\OO$\,, i.e.
$$ G_2 = \Menge{g \in \GL(\OO)}{\forall x,y\in \OO: g(x\cdot y) = g(x)\cdot g(y)}\; . $$
In this setting a group homomorphism \,$\Phi: \Sp(1)\times \Sp(1) \to G_2$\, has been described by Yokota in 
\cite{Yokota:invol1-1990}, Section~1.3: For any \,$g_1,g_2 \in \Sp(1)$\,, \,$\Phi(g_1,g_2)$\, is given by
$$ \forall x = (x_1,x_2) \in \OO \;: \; \Phi(g_1,g_2)x = (g_1\,x_1\,g_1^{-1} \,,\, g_2\,x_2\,g_1^{-1}) \; . $$
\,$\Phi$\, is in particular a totally geodesic map; one easily sees that \,$\ker(\Phi) = \{\pm(1,1)\}$\, holds,
and therefore \,$\Phi$\, is a two-fold covering map onto its image. The image is therefore a 6-dimensional totally geodesic
submanifold of \,$G_2$\, of rank \,$2$\, which is isometric to \,$(\Sp(1)\times\Sp(1))/\{\pm (1,1)\} \cong
(\Sph^3_{r=1} \times \Sph^3_{r=1/\sqrt{3}})/\{\pm (1,1)\}$\,, and which turns out to be of type \,$(\Sph\times\Sph,3,3)$\,.

The totally geodesic submanifolds of type \,$(\Sph\times\Sph,\ell,\ell')$\, correspond to the submanifolds
\,$(\Sph^\ell_{r=1} \times \Sph^{\ell'}_{r=1/\sqrt{3}})/\{\pm (1,1)\}$\, in this product, the totally geodesic submanifolds
of type \,$(\Sph,\vi=0,\ell)$\, resp.~\,$(\Sph,\vi=\tfrac\pi6,\ell')$\, correspond to the factors 
\,$\Sph^\ell_{r=1}$\, resp.~\,$\Sph^{\ell'}_{r=1/\sqrt{3}}$\, in that product, and the totally geodesic submanifolds
of type \,$(\PP,\vi=\tfrac\pi6,(\R,\ell))$\, correspond to the diagonal \,$\Menge{\pm(x,\tfrac{1}{\sqrt{3}}x)}{x\in\Sph^\ell_{r=1}}$\,
in that product.

\paragraph{Types \,$\boldsymbol{(\mathrm{A}_2)}$\, and \,$\boldsymbol{(\mathrm{AI})}$\,.}
We again realize \,$G_2$\, as the automorphism group of \,$\OO$\,. We fix an imaginary unit octonion \,$i$\, of \,$\OO$\,, and consider
the subgroup \,$H := \Menge{g\in G_2}{g(i)=i}$\, of \,$G_2$\,. \,$H$\, is isomorphic to \,$\SU(3)$\,; as totally
geodesic submanifold of \,$G_2$\,, this subgroup is of type \,$(\mathrm{A}_2)$\,. Consider the splitting \,$\OO = V \oplus W$\, of \,$\OO$\,
with \,$V := \spn_{\R}\{1,i\} \cong \C$\, and \,$W := V^\perp$\,; \,$V$\, and \,$W$\, are complex subspaces of dimension \,$1$\, resp.~\,$3$\, 
with respect to the complex structure induced by the element \,$i\in\OO$\,. Then \,$H \cong \SU(3)$\, acts trivially on \,$V$\, and in the canonical way on \,$W\cong\C^3$\,.

Fixing a real form \,$W_{\R}$\, of \,$W$\,, we obtain the subgroup \,$H' := \Menge{g\in H}{g(W_{\R})=W_{\R}}$\,, which is isomorphic to \,$\SO(3)$\,. 
\,$H/H'$\, is a Riemannian symmetric space isomorphic to \,$\SU(3)/\SO(3)$\,, and the image of the Cartan embedding \,$H/H' \to H \subset G_2$\,
is a totally geodesic submanifold of \,$G_2$\, of type \,$(\mathrm{AI})$\,.

\paragraph{Type \,$\boldsymbol{(\Sph,\vi=\arctan(\tfrac1{3\sqrt{3}}),3)}$\,.} 
Let \,$\liem'$\, be a Lie triple system of type \,$(\Sph,\vi=\arctan(\tfrac1{3\sqrt{3}}),3)$\,.
It is apparent from the part of the proof of Theorem~\ref{T:G2:cla} which handled the classification for the case \,$\rk(\liem')=1$\,, \,$\vi_0=\arctan(\tfrac1{3\sqrt{3}})$\,
that (with respect to a suitable choice of the Cartan subalgebra \,$\liea$\, of \,$\lieg_2$\, and of the positive root system \,$\Delta_+$\, corresponding
to the root system \,$\Delta$\, of \,$\lieg_2$\, with respect to \,$\liea$\,) the unit vector \,$H := \tfrac1{\sqrt{21}}\,(9\lambda_1^\sharp + 5\lambda_2^\sharp)$\,
lies in \,$\liem'$\,, 
with respect to its Cartan subalgebra \,$\R H$\, the Lie triple system \,$\liem'$\, has only one positive root \,$\alpha$\,, which is characterized by
\,$\alpha(H)=\tfrac{\sqrt{21}}{14}$\,, hence we have \,$\|\alpha^\sharp\|^2=\tfrac{3}{28} = \tfrac{1}{r^2}$\, with \,$r := \tfrac23\,\sqrt{21}$\,.

It follows that the connected, complete totally geodesic submanifold \,$M' \subset G_2$\, corresponding to \,$\liem'$\, is a symmetric space of constant curvature
\,$\tfrac{1}{r^2}$\,, and therefore isometric either to the sphere \,$\Sph^3_r$\,, or to the real projective space \,$\RP^3_{\vkap=1/r^2}$\,. To distinguish
between these two cases, we calculate the length of a closed geodesic in \,$M'$\,. 

To do so, we use the well-known fact (see \cite{Helgason:1978}, Theorem~VII.8.5, p.~322) that the unit lattice \,$\liea_e := \Menge{v\in\liea}{\exp(v) = e}$\,
is generated by the vectors \,$2\,X_\lambda$\,, where we put \,$X_\lambda := \tfrac{2\pi}{\|\lambda^\sharp\|^2}\,\lambda^\sharp \in \liea$\,, and \,$\lambda$\,
runs through all the roots of \,$G_2$\,. In this specific situation, \,$\liea_e$\, is generated by the vectors \,$2\,X_{\lambda_2} = \tfrac{4\pi}{3}\,\lambda_2^\sharp$\,
and \,$2\,X_{\lambda_5} = \tfrac{4\pi}{3}\,\lambda_5^\sharp$\,.

The length of the geodesic \,$\gamma$\, tangent to \,$H$\, equals the smallest \,$t > 0$\, so that \,$tH \in \liea_e$\, holds, i.e.~so that there exist
\,$k,\ell\in\Z$\, with \,$tH = k\cdot \tfrac{4\pi}{3}\,\lambda_2^\sharp + \ell\cdot \tfrac{4\pi}{3}\,\lambda_5^\sharp$\,. Because we have 
\,$H = \tfrac1{\sqrt{21}}(2\lambda_2^\sharp + 3\lambda_5^\sharp)$\,, that equation leads to the conditions
$$ k = 2\cdot \tfrac{3}{4\pi\cdot\sqrt{21}}\, t \qmq{and} \ell = 3\cdot \tfrac{3}{4\pi\cdot\sqrt{21}}\, t \; . $$
Therefore, the smallest \,$t>0$\, such that \,$k,\ell\in\Z$\, holds, is \,$t = \tfrac{4\pi\cdot{\sqrt{21}}}{3} = 2\pi r$\,, and hence the 
geodesic \,$\gamma$\, is closed and has the length \,$2\pi r$\,. It follows that the totally geodesic submanifold \,$M'$\, is isometric to the
sphere \,$\Sph^3_r$\,. 

\subsection[Totally geodesic submanifolds in \,$G_2/\SO(4)$\,]{Totally geodesic submanifolds in \,$\boldsymbol{G_2/\SO(4)}$\,}
\label{SSe:G2:G}

Finally we derive from the classification of the Lie triple systems resp.~totally geodesic submanifolds in \,$G_2$\, the same classification
in the totally geodesic submanifold \,$G_2/\SO(4)$\, of \,$G_2$\,.

For this purpose, we consider the Lie group \,$G_2$\, as a Riemannian symmetric space in the same way as in the Sections~\ref{SSe:G2:geometry} and 
\ref{SSe:G2:lts}, and use the names for the types of Lie triple systems of \,$\lieg_2$\, as introduced in Theorem~\ref{T:G2:cla}. 

Further, we let \,$\liem_1$\, be a Lie triple system of \,$\lieg_2$\, of type \,$(\mathrm{G})$\,, i.e.~\,$\liem_1$\, corresponds to a totally geodesic submanifold
which is locally isometric to \,$G_2/\SO(4)$\,. 

\begin{Theorem}
\label{T:G2:G:cla}
Exactly the following types of Lie triple systems of \,$G_2$\, have representatives which are contained in \,$\liem_1$\,:
\begin{itemize}
\item \,$(\mathrm{Geo},\vi=t)$\, with \,$t \in [0,\tfrac\pi6]$\,
\item \,$(\Sph,\vi=0,2)$\,
\item \,$(\Sph,\vi=\arctan(\tfrac1{3\sqrt{3}}),2)$\,
\item \,$(\Sph,\vi=\tfrac\pi6,2)$\,
\item \,$(\PP,\vi=\tfrac\pi6,(\K,2))$\, with \,$\K\in\{\R,\C\}$\,
\item \,$(\mathrm{AI})$\,
\item \,$(\Sph\times\Sph,\ell,\ell')$\, with \,$\ell,\ell' \leq 2$\,
\end{itemize}
Among these, the Lie triple systems which are maximal in \,$\liem_1$\, are: \,$(\Sph,\vi=\arctan(\tfrac1{3\sqrt{3}}),2)$\,, \,$(\PP,\vi=\tfrac\pi6,(\C,2))$\,, \,$(\mathrm{AI})$\, 
and \,$(\Sph\times\Sph,2,2)$\,.
\end{Theorem}

\beweis
Again similar to the proofs of Theorems~\ref{EIII:SO5:cla} and \ref{EIII:DIII:cla}.
\beweisende

\begin{Remark}
\label{R:G2:G:CN}
The maximal totally geodesic submanifolds of \,$G_2/\SO(4)$\, of type \,$(\Sph,\vi=\arctan(\tfrac{1}{3\sqrt{3}}),2)$\,, which are isometric
to a 2-sphere of radius \,$\tfrac23\,\sqrt{21}$\,, are missing from the classification by \textsc{Chen} and \textsc{Nagano} 
in Table~VIII of \cite{Chen/Nagano:totges2-1978}. They are in a similar ``skew'' position in \,$G_2/\SO(4)$\, as in \,$G_2$\,, compare Remark~\ref{R:G2:CN}.
\end{Remark}

We can infer the isometry type of the totally geodesic submanifolds corresponding to the Lie triple systems of \,$G_2/\SO(4)$\, from the corresponding information
on the totally geodesic submanifolds of \,$G_2$\,, given in Section~\ref{SSe:G2:tgsub}:

\begin{center}
\begin{tabular}{|c|c|c|}
\hline
type of Lie triple system & corresponding global isometry type & properties \\
\hline
$(\mathrm{Geo},\vi=t)$ & \,$\R$\, or \,$\Sph^1$\, & \\
$(\Sph,\vi=0,2)$ & \,$\Sph^2_{r=1}$\, & \\
$(\Sph,\vi=\arctan(\tfrac{1}{3\sqrt{3}}),2)$\, & $\Sph^2_{r=\tfrac23\sqrt{21}}$\, & maximal \\
$(\Sph,\vi=\tfrac\pi6,2)$ & \,$\Sph^2_{r=1/\sqrt{3}}$\, & Helgason sphere \\
$(\PP,\vi=\tfrac\pi6,(\K,2))$ & \,$\KP^2_{\vkap=3/4}$\, & \,$\K=\C$\,: maximal \\
\hline
$(\mathrm{AI})$ & \,$\SU(3)/\SO(3)_{\mathrm{srr}=\sqrt{3}}$\, & maximal \\
$(\Sph\times\Sph,\ell,\ell')$ & $(\Sph^\ell_{r=1} \times \Sph^{\ell'}_{r=1/\sqrt{3}})/\{\pm \id\}$\, & \,$\ell=\ell'=2$\,: polar, meridian, maximal \\
\hline
\end{tabular}
\end{center}

\section{Summary}
\label{Se:summary}

In the following table, we list the global isometry types of the maximal totally geodesic submanifolds of all the irreducible,
simply connected Riemannian symmetric spaces \,$M$\, of rank \,$2$\,, thereby combining information from the papers
\cite{Klein:2007-claQ}, \cite{Klein:2007-tgG2} and \cite{Klein:2007-Satake} (Section~6), as well as the present paper.

We once again use the notations from Section~\ref{SSe:EIII:tgsub} for describing the scaling factor of the invariant Riemannian metric on
the symmetric spaces involved. For the three infinite families of Grassmann manifolds \,$G_2^+(\R^n)$\,, \,$G_2(\C^n)$\, and \,$G_2(\HH^n)$\,,
we also use the notation \,${}_{\textrm{srr}=1*}$\, to denote the invariant Riemannian metric scaled in such a way that the shortest root
occurring \emph{for large \,$n$\,} has length \,$1$\,, disregarding the fact that this root might vanish for certain small values of \,$n$\,.

\begin{center}
\begin{longtable}{|c|l|}
\hline
$M$ & maximal totally geodesic submanifolds \\
\hline
\endhead
\hline
\endfoot
\,$G_2^+(\R^{n+2})_{\textrm{srr}=1*} $\, & \,$G_2^+(\R^{n+1})_{\textrm{srr}=1*}$\,, \,$(\Sph^\ell_{r=1} \times \Sph^{\ell'}_{r=1})/\Z_2$\, with \,$\ell+\ell' = n$\, \\
& \emph{for \,$n \geq 4$\, even:} \,$\CP^{n/2}_{\vkap=1/2}$\, \\
& \emph{for \,$n=2$\,:} \,$\CP^1_{\vkap=1/2} \times \RP^1_{\vkap=1/2}$\, \\
& \emph{for \,$n=3$\,:} \,$\Sph^2_{r=\sqrt{5}}$\, \\
\hline
\,$G_2(\C^{n+2})_{\textrm{srr}=1*}$\, & \,$\CP^n_{\vkap=1}$\,,  \,$G_2(\R^{n+2})_{\textrm{srr}=1*}$\,, \,$G_2(\C^{n+1})_{\textrm{srr}=1*}$\, \\
& $\CP_{\vkap=1}^{\ell} \times \CP_{\vkap=1}^{\ell'}$\, with \,$\ell+\ell'=n$\, \\
& \emph{for \,$n$\, even:} \,$\HP^{n/2}_{\vkap=1/2}$\, \\
& \emph{for \,$n=2$\,:} \,$G_2^+(\R^5)$\,, \,$(\Sph^3_{r=1/\sqrt{2}}\times\Sph^1_{r=1/\sqrt{2}})/\Z_2$\, \\
& \emph{for \,$n=4$\,:} \,$\CP^2_{\vkap=1/5}$\, \\
\hline
\enlargethispage{2em}
\,$G_2(\HH^{n+2})_{\textrm{srr}=1*}$\, & \,$\HP^n_{\vkap=1}$\,, \,$G_2(\HH^{n+1})_{\textrm{srr}=1*}$\,, \,$G_2(\C^{n+2})_{\textrm{srr}=1*}$\, \\
& $\HP_{\vkap=1}^{\ell} \times \HP_{\vkap=1}^{\ell'}$\, with \,$\ell+\ell'=n$\, \\
& \emph{for \,$n=2$\,:} \,$(\Sph^5_{r=1/\sqrt{2}}\times\Sph^1_{r=1/\sqrt{2}})/\Z_2$\,, \,$\Sp(2)_{\textrm{srr}=\sqrt{2}}$\, \\
& \emph{for \,$n=4$\,:} \,$\Sph^3_{r=2\sqrt{5}}$\, \\
& \emph{for \,$n=5$\,:} \,$\HP^2_{\vkap=1/5}$\, \\
\hline
\,$\SU(3)/\SO(3)_{\textrm{srr=1}}$ & \,$\RP^2_{\vkap=1/4}$\,, \,$(\Sph^2_{r=1} \times \Sph^1_{r=\sqrt{3}})/\Z_2$\, \\
\hline
\,$\SU(6)/\Sp(3)_{\textrm{srr=1}}$ & \,$\HP^2_{\vkap=1/4}$\,, \,$\CP^3_{\vkap=1/4}$\,, \,$\SU(3)_{\textrm{srr}=1}$\,, \,$(\Sph^5_{r=1}\times\Sph^1_{r=\sqrt{3}})/\Z_2$\,\\
\hline
\,$\SO(10)/\Ug(5)_{\textrm{srr=1}}$\, & \,$\CP^4_{\vkap=1}$\,, \,$\CP^3_{\vkap=1} \times \CP^1_{\vkap=1}$\,, \,$G_2^+(\R^8)_{\textrm{srr}=\sqrt{2}}$\,, \,$G_2(\C^5)_{\textrm{srr}=1}$\,,
\,$\SO(5)_{\textrm{srr}=1}$\, \\
\hline
\,$E_6/(\Ug(1)\cdot\Spin(10))_{\textrm{srr=1}}$\, & \,$\OP^2_{\vkap=1/2}$\,, \,$\CP^5_{\vkap=1} \times \CP^1_{\vkap=1}$\,, \,$G_2^+(\R^{10})_{\textrm{srr}=\sqrt{2}}$\,, \\
& \,$G_2(\C^6)_{\textrm{srr}=1}$\,, \,$(G_2(\HH^4)/\Z_2)_{\textrm{srr}=1}$\,, \,$\SO(10)/\Ug(5)_{\textrm{srr}=1}$\, \\
\hline
\,$(E_6/F_4)_{\textrm{srr=1}}$\, & \,$\OP^2_{\vkap=1/4}$\,, \,$\HP^3_{\vkap=1/4}$\,, \,$((\SU(6)/\Sp(3))/\Z_3)_{\textrm{srr}=1}$\,, \,$(\Sph^9_{r=1}\times\Sph^1_{r=\sqrt{3}})/\Z_4$\, \\
\hline
\,$G_2/\SO(4)_{\textrm{srr=1}}$\, & \,$\Sph^2_{r=\tfrac23\,\sqrt{21}}$\,, \,$\CP^2_{\vkap=3/4}$\,, \,$\SU(3)/\SO(3)_{\textrm{srr}=\sqrt{3}}$\,, \,$(\Sph^2_{r=1}\times\Sph^2_{r=1/\sqrt{3}})/\Z_2$\, \\
\hline
\,$\SU(3)_{\textrm{srr=1}}$\, & \,$\CP^2_{\vkap=1/4}$\,, \,$\RP^3_{\vkap=1/4}$\,, \,$\SU(3)/\SO(3)_{\textrm{srr}=1}$\,, \,$(\Sph^3_{r=1}\times\Sph^1_{r=\sqrt{3}})/\Z_2$\, \\
\hline
\,$\Sp(2)_{\textrm{srr=1}}$\, & \,$\Sph^3_{r=\sqrt{5}}$\,, \,$\HP^1_{\vkap=1/2}$\,, \,$(\Sph^3_{r=\sqrt{2}} \times\Sph^3_{r=\sqrt{2}})/\Z_2$\,,
\,$G_2^+(\R^5)_{\textrm{srr}=1}$\, \\
\hline
\,$(G_2)_{\textrm{srr=1}}$\, & \,$\Sph^3_{r=\tfrac23\,\sqrt{21}}$\,, \,$(\Sph^3_{r=1}\times\Sph^3_{r=1/\sqrt{3}})/\Z_2$\,, \,$\SU(3)_{\textrm{srr}=\sqrt{3}}$\,, \,$G_2/\SO(4)_{\textrm{srr}=1}$\, \\
\end{longtable}
\end{center}

\enlargethispage{2em}

\end{document}